\documentclass[3p]{elsarticle}

\usepackage{mydef}

\usepackage{url}
\usepackage{nicefrac}
\usepackage{bm}
\usepackage{algpseudocode}
\usepackage{algorithm}
\usepackage{multirow}              		
\usepackage{mathrsfs}
\usepackage{amsmath}
\usepackage{dsfont}
\usepackage{cancel}
\usepackage{graphicx}		
\usepackage{amssymb}
\usepackage{adjustbox}

\newcounter{qcounter}[section]

\newtheorem{thm}{Theorem}[section]
\newtheorem{cor}[thm]{Corollary}
\newtheorem{lem}[thm]{Lemma}
\newtheorem{prop}[thm]{Proposition}

\newtheorem{exmp}[thm]{Example}

\newtheorem{rem}[thm]{Remark}

\newproof{pf}{Proof}

\newcommand\numberthis{\addtocounter{equation}{1}\tag{\theequation}}

\graphicspath{{./},{./figures/}}

\title{Explicit cost bounds of stochastic Galerkin approximations for parameterized {PDEs} with random coefficients{\small \tnoteref{fn1}}}

\author[UTK]{N.C. Dexter}
\author[ORNL]{C.G. Webster}
\author[ORNL]{G. Zhang}

\address[UTK]{Department of Mathematics, University of Tennessee, Knoxville, TN 37996.}
\address[ORNL]{Department of Computational and Applied Mathematics, Oak Ridge National Laboratory, Oak Ridge, TN 37831.}

\begin{document}

\begin{abstract}
This work analyzes the overall computational complexity of the stochastic Galerkin finite element method (SGFEM) for approximating the solution of parameterized elliptic partial differential equations with both affine and non-affine random coefficients.
To compute the fully discrete solution, such approaches employ a Galerkin projection in both the deterministic and stochastic domains, produced here by a combination of finite elements and a global orthogonal basis, defined on an isotopic total degree index set, respectively. 
To account for the sparsity of the resulting system, we present a rigorous cost analysis that considers the total number of coupled finite element systems that must be simultaneously solved in the SGFEM.  
However, to maintain sparsity as the coefficient becomes increasingly nonlinear in the parameterization, it is necessary to also approximate the coefficient by an additional orthogonal expansion.
In this case we prove a rigorous complexity estimate 
for the number of 
floating point operations (FLOPs) 
required per matrix-vector multiplication of the coupled system.
Based on such complexity estimates we also develop explicit cost bounds in terms of FLOPs 
to solve the stochastic Galerkin (SG) systems to 
a prescribed tolerance, 
which are used to compare with the minimal complexity estimates of a stochastic collocation finite element method (SCFEM), shown in our previous work \cite{Galindo2015}.
Finally, computational evidence complements the theoretical estimates and supports our conclusion that, in the case that the coefficient is 
affine, the coupled SG system can be solved more efficiently than the decoupled SC systems.  However, as the coefficient becomes more nonlinear, it becomes prohibitively expensive to obtain an approximation with the SGFEM.
\end{abstract}


\tnotetext[fn1]{This material is based upon work supported in part by the U.S.~Air Force of Scientific Research under grant number 1854-V521-12; by the U.S.~Department of Energy, Office of Science, Office of Advanced Scientific Computing Research, Applied Mathematics program under contract numbers ERKJ259, and ERKJE45
and by the Laboratory Directed Research and Development program at the Oak Ridge National Laboratory, which is operated by UT-Battelle, LLC, for the U.S.~Department of Energy under Contract DE-AC05-00OR22725.}

\maketitle

\section{Introduction}
\label{sec:intro}

Nowadays, stochastic polynomial methods are widely used alternatives to Monte Carlo methods (see, e.g., \cite{Fishman_96}) for predicting the solution to physical and engineering problems described by parameterized partial differential equations (PDEs) with a finite number of random variables.  In the last decade, two classes of such methods have been proposed that often feature much faster convergence rates: {\em intrusive} stochastic Galerkin (SG) methods and {\em non-intrusive} stochastic collocation (SC) methods.  
Both approaches typically employ a Galerkin projection in the physical domain, produced here by finite elements, and the resulting fully discrete approximations only differ in their choice of multivariate polynomials for the discretization in the stochastic domain. 
For details about the relations between these methods see~\cite{GW14b, Gunzburger:2014dy, ACTA_GWZ, SpectralUQ, Nobile2009}, and for computational comparisons between the SG and SC methods see, e.g., \cite{BNTT_comp, Elman2011, ACTA_GWZ}.  

The focus of this paper is to provide explicit cost bounds for applying the stochastic Galerkin finite element method (SGFEM) to the solution of an elliptic PDE, with stochastic diffusion coefficient parameterized by finitely many random variables. 
In particular, we focus on the cost of constructing isotropic total degree SG approximations when the coefficient has both affine and non-affine dependence on the parameters. 
Under very basic assumptions on the coefficient, the solution to this problem has been shown to have analytic regularity in the random variables (see \cite{Tran2014QO}). 
As a result, SG approximations that employ a global orthogonal basis have been shown to be optimal projections in the $L^2$ sense, converging sub-exponentially with respect to the cardinality of the polynomial subspace \cite{Todor2007}. 
However, the computational cost of solving the coupled SG system does not grow linearly in the cardinality of the given subspace.
Therefore, the convergence estimates do not indicate the total complexity of obtaining the approximation for a prescribed tolerance.

When the diffusion coefficient can be written as a sum of separable functions of the physical and random parameters, the coupled SG system can be written as a sum of Kronecker products of SG matrices and finite element stiffness matrices. 
For every SG matrix, each nonzero element leads to a nonzero block of the coupled SG system, where the size of the block equals the size of the finite element stiffness matrix. 
To solve the SG system, one must simultaneously solve all the coupled finite element problems. 
In the case that the coefficient is affine in the parameters, 
   the number of nonzeros in each SG matrix 
   is of order $\mathcal{O}(\Mp)$ \cite{Ernst2010}, where $\Mp$ is the cardinality of the isotropic total degree polynomial subspace of order $p\in\N$. 
Thus, a matrix-vector product involving the coupled SG system requires $\mathcal{O}(\Jh \Mp)$ floating point operations (FLOPs), where $\Jh$ is the number of physical degrees of freedom. 
Therefore, the work of solving the coupled SG system when employing an iterative method, e.g., conjugate gradient (CG), is of the order $\mathcal{O}(\Jh \Mp \NiterSG)$
where $\NiterSG$ is the number of iterations 
required to achieve a prescribed accuracy of the fully discrete approximation
\cite{BNTT_comp,Ernst2010,Ullmann2010}.

On the other hand, when the diffusion coefficient is a general non-affine function of the random parameters, the cost of obtaining an approximation with the SGFEM is not as obvious as before.
In this setting we consider two cases, namely, the coefficient is: (1) a polynomial with respect to the random variables, and; (2) a transcendental function with respect to the random variables.
In the first case, as we increase the order of the polynomial, the block-sparsity of the SG system decreases, resulting in a SG system that {\em incrementally} becomes block-dense 
\cite{Eiermann2007,Ernst2010,SpectralUQ,Matthies2005,Ullmann2010,Ullmann2012}. 
In the second case, a separable representation can be guaranteed with the use of an orthogonal expansion \cite{WienerPC,XiuKarniadakisgPC}, such that, substituting the expansion into the discretized PDE recovers the Kronecker product structure.
However, when the expansion is not truncated, the SG system is known to be {\em entirely} block-dense \cite{Ernst2010,Matthies2005}.  
Without a priori knowledge on the exact sparsity of the SG matrices 
in this case, it was estimated that the complexity of matrix-vector multiplications of the SG Kronecker product system is between
$\mathcal{O}(\Jh \Mp^2)$ and $\mathcal{O}(\Jh \Mp^3)$ \cite{Ullmann2010}. 
As such, it is impossible to make a conclusive statement about the computational cost, and, more importantly, does not account for the two cases above, i.e., when the coefficient is possibly a truncated polynomial of fixed total degree $r\in\N$ such that $1 \leq r < \infty$.
In these cases, the work of solving the coupled SG system with an iterative method is given by $\mathcal{O}(\Jh \Mpr \NiterSG)$,
where $\Mpr$ is the total number of $\mathcal{O}(\Jh)$ finite element problems that must be simultaneously solved.

The key challenge of estimating the cost of solving the SG system when the coefficient is a (truncated) polynomial of finite order is to provide bounds on the block-sparsity of the matrix, i.e., nonzeros of the SG system. 
To achieve this, 
we provide a rigorous counting argument, which can be seen as a generalization of results from \cite{Ernst2010}, 
for the exact sparsity of the SG matrices 
for an arbitrary order orthogonal expansion of a non-affine coefficient.
As a result, 
we are able to provide bounds for $\Mpr$
of the order
$\mathcal{O}(\Mp \Mr \min\{2^r, M_{\lceil r/2\rceil}\})$,
where $\Mr$ is the cardinality of the total degree polynomial subspace used in an orthogonal expansion of order $r$ of the coefficient. This result provides sharper estimates than the bounds in the case of the full orthogonal expansion from \cite{Ullmann2012} since it depends on the truncation order $r$, and allows us to estimate the total complexity of solving the coupled system for general non-affine coefficients.
Since the counting argument for the sparsity of the SG system 
relies only on the SG discretization of an elliptic operator in terms of orthogonal polynomials,
we note that this argument can be reused to estimate the complexity of solving similarly defined PDEs with this method.

In addition, we also develop explicit cost bounds in terms of FLOPs to solve the SG system. 
Our approach relies on $\eps$-complexity analysis, wherein we balance the errors arising from the approximation with the SGFEM and the iterative solver, e.g., CG, so as to ensure the solution to the fully discrete approximation achieves a given tolerance of $\eps>0$.
With this result, we are able to provide a direct comparison with $\eps$-complexity estimates for the stochastic collocation finite element method (SCFEM) in our previous work \cite{Galindo2015}. 
Finally, we present numerical results in agreement with the theoretical work estimates for both the SGFEM and SCFEM all cases described above.

An outline of the paper is as follows. 
In \S \ref{sec:prob}, we provide a discussion on the model problem, and requirements on the diffusion coefficient. 
In \S \ref{sec:SG}, 
we define the parameterized finite element and SG approximations, 
derive the SG system, and provide examples of the resulting linear systems that arise from the SG discretization with various coefficients. We then define the cost of solving the SG system and discuss preconditioning strategies. 
In \S \ref{sec:complexity_analysis}, we derive the exact number of coupled finite element problems in the SG system and bounds on the sparsity in the non-affine case, and present explicit cost bounds of the SGFEM. We also discuss the conditioning of the system in the non-affine case in order to provide a comparison with similar results from \cite{Powell2009}.
In \S \ref{sec:SC_comparison}, we briefly describe the SCFEM, and provide theoretical comparison with results from \cite{Galindo2015} in terms of minimum work to reach a given tolerance, both in the affine and non-affine cases.
Finally, in \S \ref{sec:num_ex}, we present illustrative numerical examples corroborating our theoretical results.

\section{Problem setting}
\label{sec:prob}

We consider the simultaneous solution of the parameterized linear elliptic PDE:
\begin{align}
\label{eq:model_problem}
\left\{\begin{array}{rll} 
-\nabla \cdot \left( a(x,\y) \nabla u(x,\y) \right) \hspace{-0.25cm} &= f(x)  &\forall x\in D, \; \y \in \Gamma \\
                                                u(x,\y) \hspace{-0.25cm} &= 0  &\forall x\in \partial D, \; \y \in \Gamma
\end{array}\right.
\end{align}
where $f\in L^2(D)$ is a fixed function of $x$, $D\subset \R^d$, $d=1,2,3$, is a bounded Lipschitz domain, and $\y(\omega)=(y_1(\omega),\ldots,y_N(\omega)):\Omega\to\Gamma = \prod_{i=1}^N \Gamma_i \subseteq \R^N$ is a random vector with $\omega\in\Omega$ and $\Omega$ the set of outcomes. In this setting we assume the components of $\y$ have a joint probability density function $\varrho : \Gamma \to \R_+$, with $\varrho(\y) = \prod_{i=1}^N \varrho_i(y_i)$ 
known directly through, e.g., truncations of correlated random fields \cite{loeve1978probability} in $(\Gamma,\mathcal{B}(\Gamma),\varrho(\y)d\y)$, where $\mathcal{B}(\y)$ denotes the Borel $\sigma$-algebra on $\Gamma$ and $\varrho(\y)d\y$ is the probability measure of $\y$.
We further assume that $\varrho_i$ is an even weight function for each $i=1,\ldots,N$. 
We require the following assumptions related to the continuity, coercivity, and holomorphic dependence of the coefficient $a(x,\y)$. Namely:
{\em 
\begin{list}{\em (A\arabic{qcounter})~}{\usecounter{qcounter}}
\item There exist constants $0<a_{\min} \leq a_{\max}<\infty$ such that for all $x\in \overline{D}$ and $\y\in \Gamma$, 
$
a_{\min} \leq a(x,\y) \leq a_{\max}.
$
\item The complex continuation of $a(x,\y)$, denoted $a^*: \mathbb{C}^N \to L^\infty$, is a $L^\infty(D)$-valued holomorphic function on $\mathbb{C}^N$.
\end{list}
}

The holomorphic dependence on $\y$ of the coefficient $a(x,\y)$ holds in many examples, including polynomial, exponential, and trigonometric functions of the variables $y_1,\ldots,y_N$ shown below. 

\begin{exmp}[\rm The affine case]
\label{exmp:affine_coefficient}
We consider an affine function of the random parameters, e.g.,
\begin{align}
\label{eq:affine_coefficient}
a(x,\y) = a_0(x) + \sum_{k = 1}^N y_k b_k(x), \;\;\; x \in \overline{D}, \; \y \in \Gamma,
\end{align}
where $a_0$, $\{b_k\}_{k=1}^N \subset L^2(D)$ are such that $a(x,\y)$ satisfies $\Aone$. Such examples include general \KL expansions \cite{loeve1978probability} or piecewise constant random fields. 
\end{exmp}

\begin{exmp}[\rm The non-affine, polynomial case]
\label{exmp:polynomial_coefficient}
We consider a non-affine, polynomial function of the random parameters, e.g.,
\begin{align}
\label{eq:polynomial_coefficient}
a(x,\y) = a_0(x) + \sum_{1 \leq |\bmalpha| \leq \rbar}  \y^{\bmalpha} c_{\bmalpha} (x), \;\;\; x \in \overline{D}, \; \y \in \Gamma,
\end{align}
where $\bmalpha = (\alpha_1,\ldots,\alpha_N)$ is a multi-index, $|\bmalpha| = \alpha_1 + \cdots + \alpha_N$, $\y^{\bmalpha} = y_1^{\alpha_1} \cdots y_N^{\alpha_N}$, $\rbar < \infty$ is the polynomial order of $a(x,\y)$, and $a_0, \{c_{\bmalpha}\}_{|\bmalpha| \leq \rbar} \subset L^2(D)$ are such that $a(x,\y)$ satisfies $\Aone$. 
Examples include fixed-order Taylor or orthogonal expansions of a general random field.
\end{exmp}

\begin{exmp}[\rm The non-affine, transcendental case]
\label{exmp:transcendental_coefficient}
We consider a non-affine, transcendental function of the random parameters, e.g.,
\begin{align}
\label{eq:transcendental_coefficient}
a(x,\y) = a_0(x) + g(x,\y), \;\;\; x \in \overline{D}, \; \y \in \Gamma,
\end{align}
where $a_0, g \subset L^2(D)$, and $g(x,\y)$ is a general transcendental function of $x$ and $\y$, such that $a(x,\y)$ satisfies $\Aone$. 
Examples of $g(x,\y)$ include the sine, logarithm, or exponential functions of \eqref{eq:affine_coefficient} or \eqref{eq:polynomial_coefficient}.
\end{exmp}

Let $L^2_\varrho(\Gamma)$ be the space of square integrable functions with respect to the measure $\varrho(\y) d\y$ and $L^\infty_\varrho(\Gamma)$ 
be the space of essentially bounded functions, with the norm 
\begin{align*}
\|u\|_{L^\infty_{\varrho}(\Gamma)}:= \text{ess} \sup_{\y\in\Gamma} |u(\y)|,
\end{align*}
where the essential suppremum is taken with respect to the weight $\varrho$. By $H^{-1}(D)$ we denote the dual of $H_0^1(D)$, the space of square integrable functions in $D$ having zero trace on the boundary and square integrable distributional derivatives. 
We will often use the abbreviation $\Htworho$ to denote the space
\begin{align*}
L^2_\varrho(\Gamma ; H_0^1(D)): = \left\{ u : \overline{D} \times \Gamma \to \R : u \text{ strongly measurable and } \int_\Gamma \| u \|_{H_0^1(D)}^2 \varrho(\y) d\y < \infty \right\},
\end{align*}
and $\Hinfrho$ to denote the space
\begin{align*}
L^\infty_\varrho(\Gamma; H_0^1(D)): = \left\{ u : \overline{D} \times \Gamma \to \R : u \text{ strongly measurable and } \textnormal{ess} \sup_{\y\in\Gamma} \|u(\cdot, \y)\|_{H_0^1(D)} < \infty \right\}. 
\end{align*}
For the space $H_0^1(D)$ we have the energy norm $\| v \|_{H_0^1(D)} = \| \nabla v \|_{L^2(D)}$, hence $\Htworho$ is a Hilbert space with norm $\| v \|_{\Htworho}^2 = \int_\Gamma \| v \|_{H_0^1(D)}^2 \varrho d\y$. 
The {\em stochastic weak form} of problem \eqref{eq:model_problem} is given by: {\em find $u\in \Htworho$ such that $\forall v\in \Htworho$} 
\begin{equation}
\begin{aligned}
\label{eq:weak_problem}
\int_\Gamma \mathcal{B}[u,v](\y) \varrho(\y)\; d\y = \int_\Gamma F(v) \varrho(\y)\; d\y,
\end{aligned}
\end{equation}
where
\begin{equation}
\begin{aligned}
\label{eq:deterministic_weak_problem_BF}
\mathcal{B}[u,v](\y) = \int_D a(x,\y) \nabla u(x,\y) \cdot \nabla v(x,\y) dx, \;\;\;\; F(v) = \int_D f(x) v(x,\y) dx.
\end{aligned}
\end{equation}
For convenience, we will often use the abbreviation $\mathcal{B}(\y) = \mathcal{B}[\cdot,\cdot](\y)$ and suppress the dependence on $x\in D$ in writing $a(\y) = a(\cdot,\y)$ and $u(\y) = u(\cdot,\y)$. 
It follows from $\Aone$ that $\mathcal{B}(\y)$ is a symmetric, uniformly coercive, and continuous bilinear operator on $H_0^1(D)$, parameterized by $\y\in\Gamma$, and $\mathcal{B}(\y)$ induces the norm 
\begin{align}
\label{eq:Bnorm}
\|u\|^2_{\mathcal{B}(\y)} := \int_D a(x,\y) |\nabla u|^2 dx. 
\end{align}
Assumption $\Aone$ and the Lax-Milgram lemma also ensure the existence and uniqueness of the solution $u$ to \eqref{eq:weak_problem} in $\Htworho$.

\begin{figure}[!htb]
\begin{center}
\includegraphics[width=4.27cm]{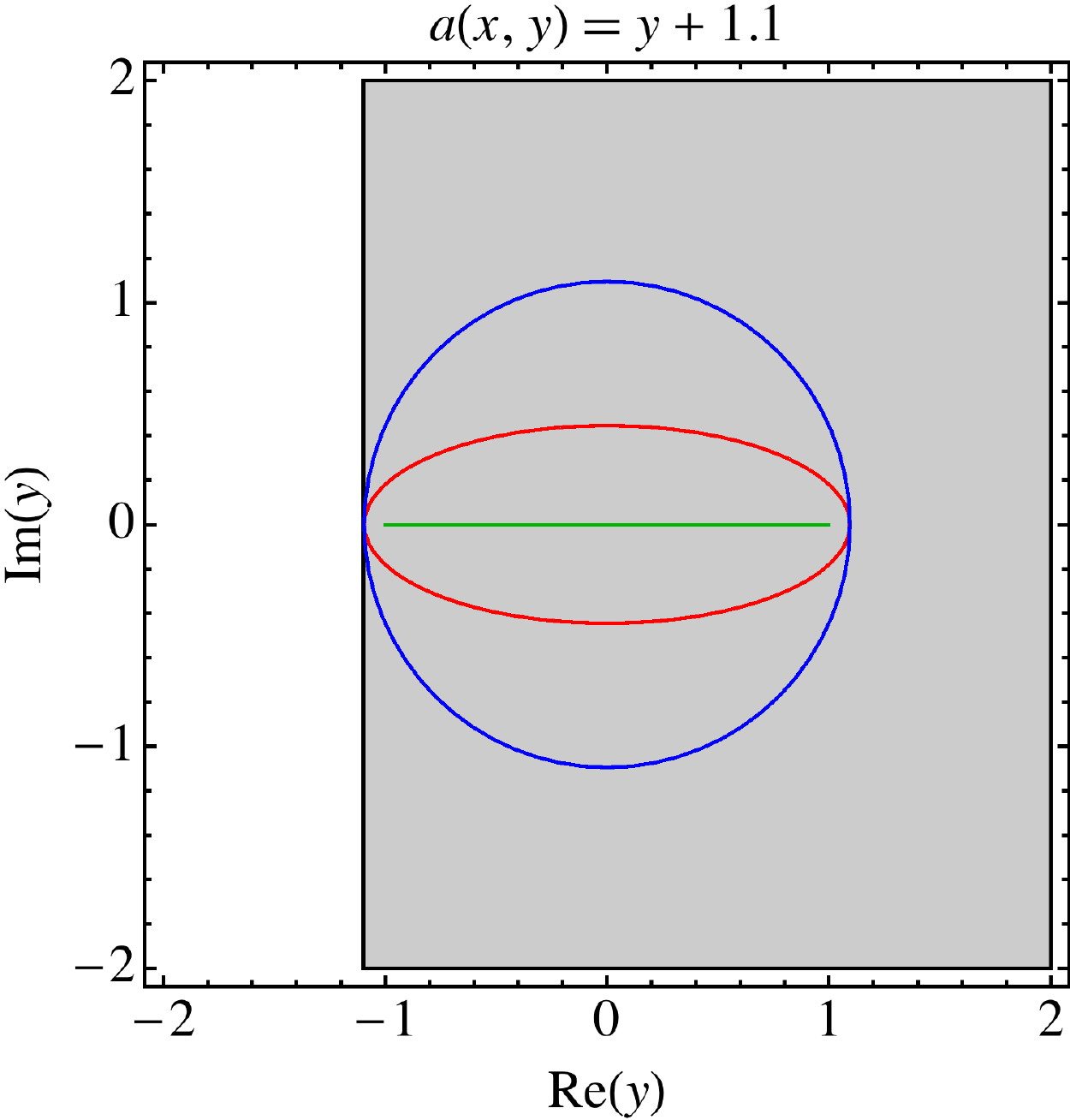}
\includegraphics[width=4.27cm]{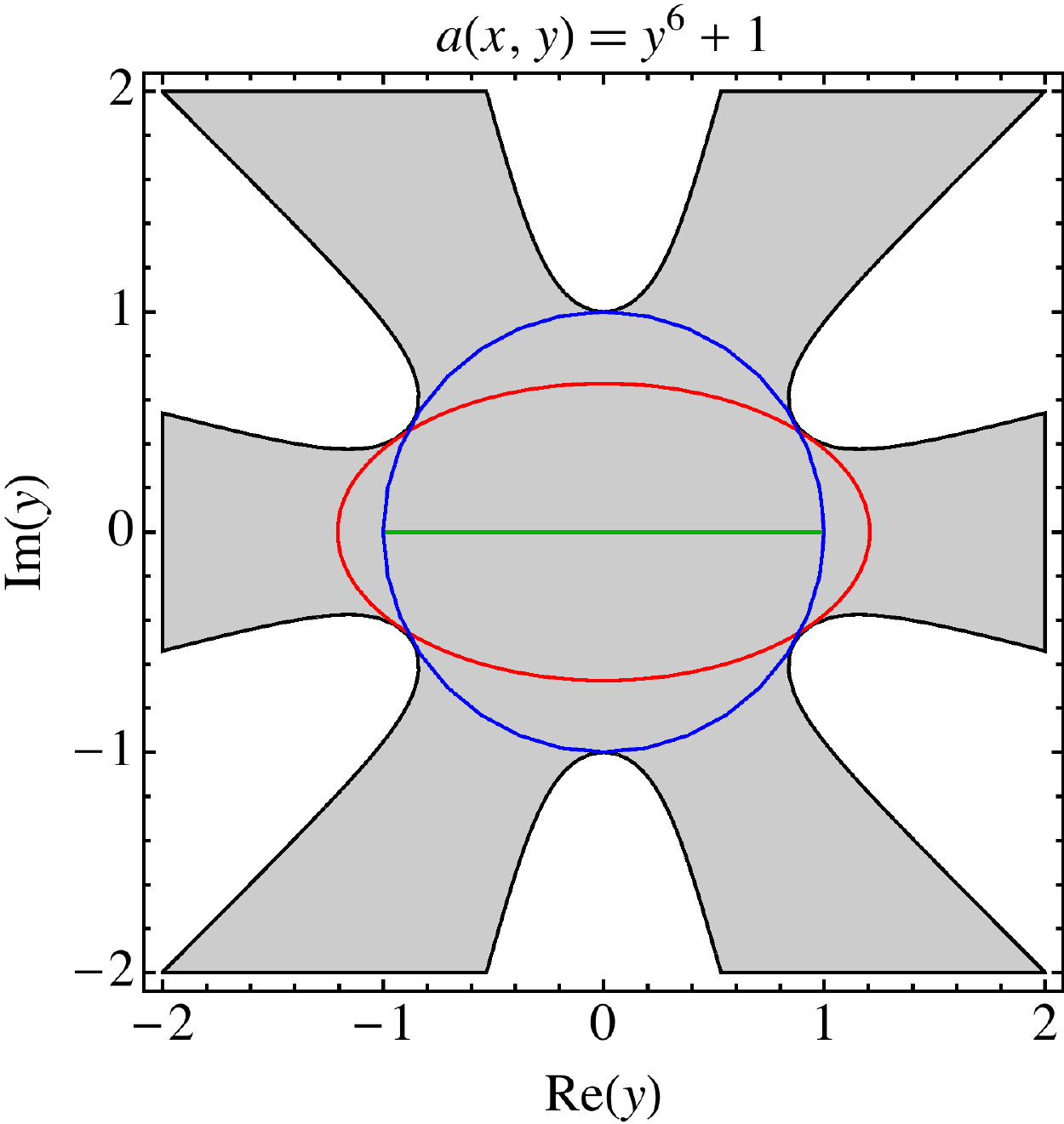}
\includegraphics[width=4.27cm]{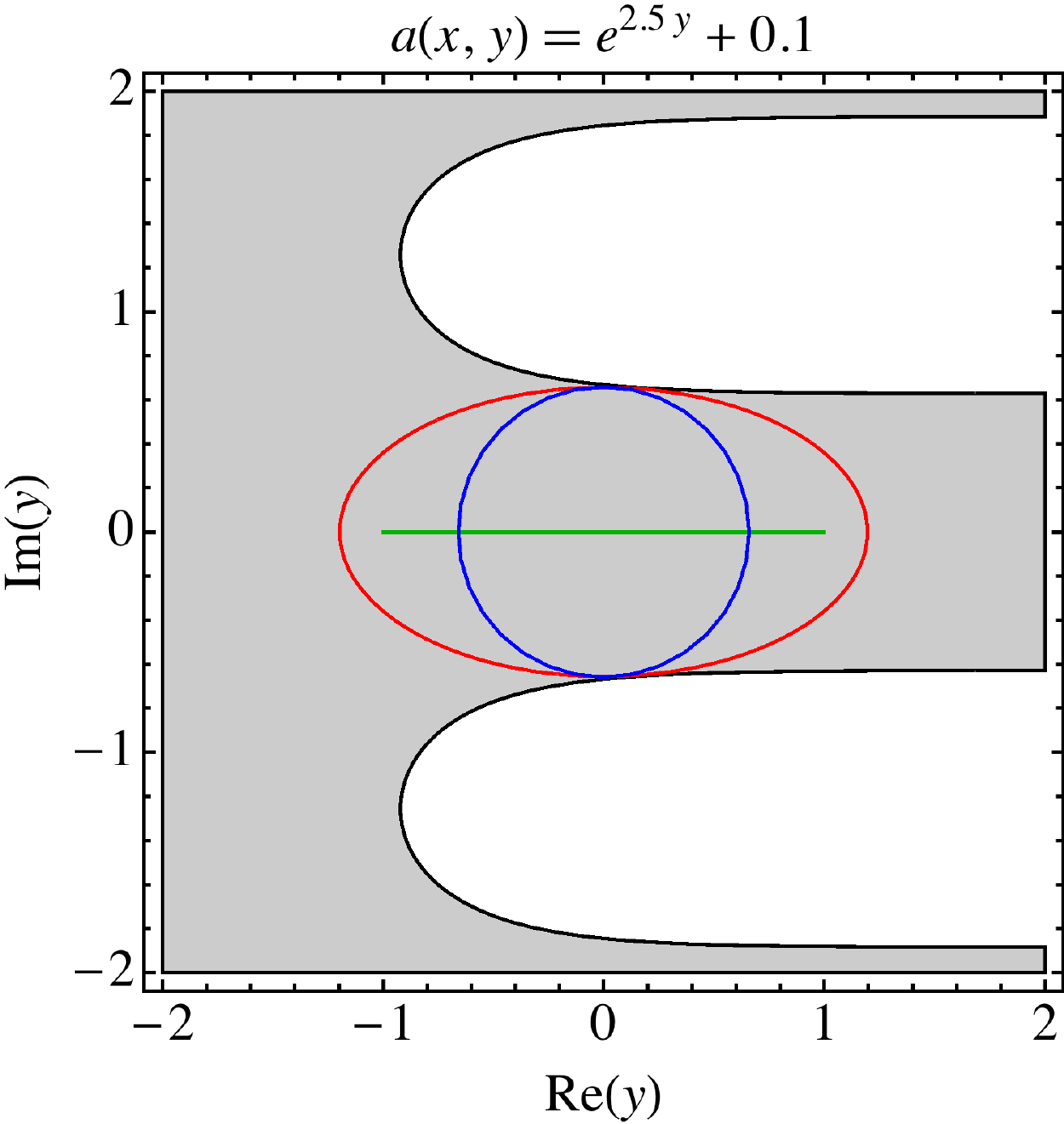}
\caption{Domains of uniform ellipticity for some one-dimensional coefficients $a(x,y)$ are indicated by the gray regions in each plot. The blue and red curves represent the maximal discs and ellipses, respectively, that can be contained in those domains, and the green lines represent the interval $\Gamma =[-1,1]$.}
\label{fig:uniform_ellipticity_regions}
\end{center}
\end{figure}

The convergence of the global stochastic polynomial methods used to approximate \eqref{eq:model_problem} exploits the uniform ellipticity of the coefficient $a(\y)$ and depends on the regularity of $u(\y)$ with respect to $\y$. By $\text{Re}(z)$ and $\text{Im}(z)$ we denote the real and imaginary parts of $z\in \mathbb{C}$, and for $0< \delta < a_{\min}$ we define
\begin{align}
\label{eq:domain_uniform_ellipticity}
U(a,\delta) = \{\z \in \C^N : \text{Re}(a(x,\z)) \geq \delta, \forall x\in \overline{D} \}.
\end{align}
If $U(a,\delta) \neq \emptyset$ for some $0<\delta<a_{\min}$, we say that $a(x,\z)$ is uniformly elliptic on the set $U(a,\delta)$ and we refer to $U(a,\delta)$ as its domain of uniform ellipticity. For $\bs{\gamma} = (\gamma_1,\ldots,\gamma_N)$ with $\gamma_i>1$ $\forall i$ we denote the polyellipse
\begin{align*}
\mathcal{E}_{\bs\gamma} = \bigotimes_{1\leq i \leq N} \left\{ z_i \in \mathbb{C} : \text{Re}(z_i) \leq \frac{\gamma_i + \gamma_i^{-1}}{2} \cos\phi, \; \text{Im}(z_i) \leq \frac{\gamma_i - \gamma_i^{-1}}{2} \sin\phi, \; \phi\in [0,2\pi) \right\}.
\end{align*}
In \cite{Tran2014QO} it was shown that if $a(\y)$ satisfies $\Aone$ and $\Atwo$, then for any $0 < \delta < a_{\min}$ there exists a $\bs{\gamma} = (\gamma_1,\ldots,\gamma_N)$ with $\gamma_i > 1$ $\forall i$ such that $\mathcal{E}_{\bs\gamma} \subset U(a,\delta)$.
We can also similarly define 
the polydisc $\mathcal{D}_{\bs\gamma} = \bigotimes_{1\leq i \leq N} \{ z_i \in \C : |z_i| \leq \gamma_i \}$, though, for arbitrary $0<\delta<a_{\min}$, it is not always possible to find a $\bs{\gamma}$ with $\gamma_i > 1$ $\forall i$ such that $\mathcal{D}_{\bs\gamma} \subset U(a,\delta)$.
Figure \ref{fig:uniform_ellipticity_regions} provides an illustration of this fact for various one-dimensional coefficients $a(\y)$, $\y\in \mathbb{C}$.
Note that in the case of the 6th degree polynomial and exponential random variables, no disc of radius $\gamma>1$ containing $\Gamma = [-1,1]$ can fit in the region. 
The following theorem, proved in \cite{Tran2014QO}, shows the regularity of the solution $u$ with respect to the parameterization.

\begin{thm}
\label{thm:regularity}
When the coefficient $a(x,\y)$ satisfies $\Aone$ and $\Atwo$, so that for some $0 < \delta < a_{\min}$ and $\bs{\gamma} = (\gamma_1,\ldots,\gamma_N)$ with $\gamma_i > 1$ $\forall i$ we have $\mathcal{E}_{\bs\gamma} \subset U(a,\delta)$, then the function $\z \mapsto u(\z)$ from \eqref{eq:model_problem} is holomorphic in an open neighborhood of $\mathcal{E}_\bs{\gamma}$. 
\end{thm}

This result states that a direct consequence of the uniform ellipticity of the function $a(x,\y)$ on the polyellipse $\mathcal{E}_{\bs\gamma} \subset U(a,\delta)$ is that the solution $u$ of \eqref{eq:model_problem} has analytically smooth dependence on the parameterization $\y$.
Theorem \ref{thm:regularity} is the key in motivating the construction of global stochastic Galerkin (SG) approximations to the solution $u$ of \eqref{eq:model_problem}, to be described in the following sections. 

\section{Stochastic Galerkin finite element method}
\label{sec:SG}

In this section we define the SGFEM for constructing fully discrete approximations to the solution $u$ of problem \eqref{eq:model_problem}. This discretization employs mixed Galerkin projections in the spatial and parameter domains. In particular we rely on the finite element method for the spatial discretization, described in \S \ref{subsec:fem}, and the stochastic Galerkin method for the parameter discretization, described in \S \ref{subsec:stochastic_galerkin_approximation}. In \S \ref{subsec:representations_of_a} we describe the linear systems that result from the SG discretization when Examples \ref{exmp:affine_coefficient}, \ref{exmp:polynomial_coefficient}, and \ref{exmp:transcendental_coefficient} are used in problem \eqref{eq:model_problem}. We then conclude in \S \ref{subsec:solving_SG} with a discussion of the cost of solving the SG systems.

\subsection{Parameterized finite element approximation}
\label{subsec:fem}

We briefly define the finite element method for obtaining a discretization of $u$ from \eqref{eq:model_problem} over the spatial domain $D$.
Let $\mathcal{T}_h$, be a triangulation of $D$ with maximum mesh size $h>0$, and $V_h(D)\subset H_0^1(D)$ a finite element space of piecewise continuous polynomials on $\mathcal{T}_h$ parameterized by $h\to 0$. Let $\{\phi_j(x)\}_{j=1}^{\Jh}$ denote a finite basis of $V_h(D)$ of dimension $\Jh$. We can write the {\em semi-discrete} problem as: {\em find $\uh(\y)\in V_h(D)$ such that $\forall v\in V_h(D)$}
\begin{align}
\label{eq:semi_discrete_SFEM}
\mathcal{B}[\uh(\y),v](\y) = F(v),
\end{align}
where $\mathcal{B}[\cdot,\cdot](\y)$ and $F(\cdot)$ are defined in \eqref{eq:deterministic_weak_problem_BF}.
For almost every $\y\in\Gamma$, problem \eqref{eq:semi_discrete_SFEM} admits a unique solution of the form
$\uh(x,\y) = \sum_{j=1}^{\Jh} u_j(\y) \phi_j(x)$.
We discretize problem \eqref{eq:semi_discrete_SFEM} by defining, for $i,j = 1,\ldots,\Jh$,
\begin{align}
\label{eq:A_ij}
\mathbf{[A]}_{i,j}(\y) & = \mathcal{B}[ \phi_j ,\phi_i ] (\y), \qquad
\mathbf{F}_{i} = F(\phi_i).
\end{align}
The coefficients $\mathbf{u}_{h}(\y)=[u_1(\y),u_2(\y),\ldots, u_{\Jh}(\y)]^{\text{T}}$ of $\uh(x,\y)$ are determined by solving the linear system 
\begin{align}
\label{eq:param_FEM_alg}
\mathbf{A}(\y) \mathbf{u}_{h}(\y) = \mathbf{F}, 
\end{align}
at fixed realizations of $\y\in\Gamma$. Here $\mathbf{A}(\y)$ is symmetric and positive-definite so that \eqref{eq:param_FEM_alg} can be solved by iterative methods such as the conjugate gradient (CG) method. 

\subsection{Stochastic Galerkin approximation with an orthogonal basis}
\label{subsec:stochastic_galerkin_approximation}

Based on the smoothness of the solution $u$ to \eqref{eq:model_problem}, characterized by Theorem \ref{thm:regularity}, we now consider the construction of approximations to $u$ in terms of global polynomials. 
Let $\Lambdap \subset \N_0^N$ be a finite set of multi-indices, e.g., having dimension $\#\Lambdap < \infty$, and define the space of polynomials $\mathcal{P}_{\Lambdap}(\Gamma) = \text{span}\{ \y^\nub : \nub \in \Lambdap\}$.
A general global polynomial approximation problem can be framed in terms of solving for the $\#\Lambdap$ stochastic degrees of freedom (SDOF) $\{u_\p\}_{\p\in\Lambdap}$.
When an interpolatory approach is used, the resulting systems of equations are decoupled finite element systems. 
When Galerkin projection with an orthogonal basis is used, the finite element systems are fully coupled, and must be solved simultaneously. 
Some isotropic examples of such index sets include 
\begin{equation}
\begin{aligned}
\label{eq:index_set_examples} 
\LambdapTP = \left\{ \p\in\N^N_0 : \max_{1\leq i \leq N} p_i \leq p \right\}, \;\;\;\;
\LambdapTD = \left\{ \p\in\N^N_0 : \sum_{n=1}^N p_n \leq p \right\}, \\
\LambdapSM = \left\{ \p\in\N^N_0 : \sum_{n=1}^N f(p_n) \leq f(p) \right\}, 
\;\;\;\; f(p) = \left\{ 
    \begin{array}{rl} 
    0, & p = 0 \\ 
    1, & p = 1 \\ 
    \lceil \log_2(p) \rceil, & p \geq 2 
    \end{array} \right. 
\end{aligned}
\end{equation}
corresponding to the Tensor Products (TP), Total Degree (TD), and Smolyak (SM) polynomial spaces $\mathcal{P}_{\LambdapTP} (\Gamma)$, $\mathcal{P}_{\LambdapTD} (\Gamma)$, and $\mathcal{P}_{\LambdapSM} (\Gamma)$, respectively.
When the solution $u$ exhibits an anisotropic dependence on the parameters $\y$, anisotropic weighted versions of the index sets defined in 
\eqref{eq:index_set_examples}
can be introduced to further reduce the number of SDOF needed to approximate $u$ at a desired accuracy \cite{BNTT_comp,WebsterAnisotropic}. 

\begin{rem}[\rm Best $M$-term and quasi-optimal approximations]
\label{rem:best_M_terms}
The optimal choice of $\Lambdap$ would be the set $\Lambda$ of cardinality $M$ such that the corresponding approximation provides maximum accuracy out of all sets of size $M$. Such approximations are referred to as best $M$-term approximations, and recent work has focussed on the construction of best $M$-term Taylor and Galerkin approximations \cite{Chkifa2014,Cohen2012,Cohen2010,Cohen2011,Tran2014QO}. 
These approaches construct $\Lambda$ by utilizing the largest $M$ coefficients $u_\p$ or sharp upper bounds of $u_\p$. However, in this effort we focus on analyzing the computational complexity of finding solutions to \eqref{eq:weak_problem} in $\mathcal{P}_{\Lambdap}(\Gamma)$ for a prescribed index set $\Lambdap$.
\end{rem}

For each $n=1,\ldots,N$, let $\{\psi_{p_n}(y_n)\}_{p_n=1}^\infty$ be a sequence of univariate polynomials over $\Gamma_n$, orthonormal with respect to the $L^2_{\varrho_n}(\Gamma_n)$ inner product. 
Then $\{\Psi_\p(\y)\}_{0 \leq |\p|}$ with $\Psi_\p(\y):=\prod_{n=1}^N \psi_{p_n}(y_n)$ is a sequence of multivariate polynomials over $\Gamma$, orthonormal with respect to the $L^2_\varrho(\Gamma)$ inner product. 
In the case that $\varrho= \frac{1}{2}$ for each $n=1,\ldots,N$, $\{\psi_{p_n}\}_{p_n=1}^\infty$ and $\{\Psi_\p\}_{0 \leq |\p|}$ are the univariate and multivariate Legendre polynomials, respectively. 
Given a specific choice of index set $\Lambdap$, 
it follows that $\{\Psi_\p\}_{\p\in \Lambdap}$ forms a basis of $\mathcal{P}_{\Lambdap}(\Gamma)$ with dimension $\Mp = \dim (\mathcal{P}_{\Lambdap}(\Gamma))=\#\Lambdap$. 
Hence, with $\{\phi_j\}_{j=1}^{\Jh}$ as in \S \ref{subsec:fem} and $\{\Psi_\p\}_{\p\in\Lambdap}$ as above, we can now write the fully discrete stochastic Galerkin (SG) approximation as 
\begin{align}
\label{eq:fully_discrete_SG_soln}
\uSG (x,\y) = \sum_{\p \in \Lambdap} \sum_{j=1}^{\Jh} u_{j,\p} \, \phi_j(x) \Psi_\p (\y),
\end{align}
whose coefficients can be found by solving the following coupled problem: 
{\em find $\uSG\in V_h(D) \otimes \mathcal{P}_{\Lambdap}(\Gamma)$ such that for all $v \in V_h(D) \otimes \mathcal{P}_{\Lambdap}(\Gamma)$ }
\begin{align}
\label{eq:fully_discrete_SG_problem}
\bbE\left[ \mathcal{B}[\uSG,v](\y) \vphantom{\big\vert} \right] = \bbE \left[ F(v) \vphantom{\big\vert} \right] ,
\end{align}
where $\mathcal{B}[\cdot,\cdot](\y)$ and $F(\cdot)$ are defined in \eqref{eq:deterministic_weak_problem_BF}.
To form the linear system of equations resulting from the SG approximation given by \eqref{eq:fully_discrete_SG_soln}, we let $\mathbf{u}_{h,\p} = [u_{1,\p},\ldots,u_{\Jh,\p}]^{\text{T}}$ be the vector of nodal values of the finite element solution corresponding to the $\p$-th stochastic mode of $u_{h,p}$, and $\mathbf{u}_{h,p} = [ \mathbf{u}_{h,\p} ]_{\p\in\Lambdap}^{\text{T}}$. Observe that when $f$ is deterministic
$\langle \Psi_\p \mathbf{F}_i \rangle = \mathbf{F}_i \delta_{\bs{0},\p}$ for all $i=1,\ldots,\Jh$,
where $\delta_{\bs{0},\p} = 1$ if $\p = \bs{0}$ and $\delta_{\bs{0},\p} = 0$ otherwise. Performing a Galerkin projection onto $\text{span}\{\Psi_{\p}\}_{\p\in\Lambdap}$ for the solution of \eqref{eq:fully_discrete_SG_problem} yields the following system: for each $\p\in\Lambdap$
\begin{align}
\label{eq:fully_discrete_SG_algebraic_problem}
\sum_{\q\in\Lambdap} \langle \Psi_\p(\y), \mathbf{A}(\y)\Psi_\q(\y)\rangle \mathbf{u}_{h,\q} = \langle \Psi_\p(\y), \mathbf{F} \rangle, 
\end{align}
which can be written algebraically as a system of {\em fully coupled finite element problems}: for each $\p\in\Lambdap$
\begin{align}
\label{eq:fully_discrete_SG_algebraic_KuF}
\sum_{\q\in \Lambdap} [\mathbf{K}]_{\p,\q} \mathbf{u}_{h,\q} = \mathbf{F}\delta_{\bs{0},\p} 
\end{align}
with $[\mathbf{K}]_{\p,\q} = \langle\Psi_\p(\y), \mathbf{A}(\y) \Psi_\q(\y) \rangle$ and $\mathbf{A}(\y)$ as given in \eqref{eq:A_ij}. 

\begin{rem}
\label{rem:K_storage_and_computation}
Typically matrix free methods are applied to solve \eqref{eq:fully_discrete_SG_algebraic_KuF} without ever explicitly forming $\mathbf{K}$ in memory, as described in \cite{Pellissetti00iterativesolution}. When the resulting system is sparse, as a result of an affine coefficient $a(x,\y)$, e.g., Example \ref{exmp:affine_coefficient}, this can lead to computationally efficient solution strategies. 
However, these implementations rely on the fact that the coefficient $a(x,\y)$ can be written as a sum of separable functions of $x$ and $\y$, e.g., $a(x,\y) = \sum_{j=1}^N b_j(x) c_j(\y)$. 
For the transcendental function $a(x,\y)$ from Example \ref{exmp:transcendental_coefficient}, this may not be the case. 
Moreover, when $\mathbf{K}$ is block-dense, matrix-vector multiplications require approximately $\mathcal{O}(\Jh \Mp^2)$ floating point operations (FLOPs), so that when iterative methods are used, the solution of the fully coupled finite element problems given in \eqref{eq:fully_discrete_SG_algebraic_KuF} becomes unfeasible.
\end{rem}

\subsection{Representations of $a(\y)$ and the corresponding matrix $\mathbf{K}$} 
\label{subsec:representations_of_a}

For a general coefficient $a(x,\y)$, the matrix $\mathbf{K}$ in \eqref{eq:fully_discrete_SG_algebraic_KuF} requires the storage of at most $\Mp^2$ block matrices of the size and sparsity of $\mathbf{A}(\y)$, i.e., $\mathcal{O}(J_h \Mp^2)$ elements. However, in several specific cases the actual block-sparsity of $\mathbf{K}$ is much less. 
We recall the coefficient from Example \ref{exmp:affine_coefficient}, where $\mathbf{K}$ can be rewritten
\begin{align*}
[\mathbf{K}]_{\p,\q} = \langle \Psi_\p(\y), \Psi_\q(\y)\rangle \mathbf{A}_0 + \sum_{k=1}^N \langle y_k \Psi_\p(\y), \Psi_\q(\y)\rangle \mathbf{A}_k,
\end{align*}
with
$[\mathbf{A}_0]_{i,j} = \int_D a_0(x) \nabla \phi_j (x) \cdot \nabla \phi_i(x) dx$ and
$[\mathbf{A}_k]_{i,j} = \int_D b_k(x) \nabla \phi_j (x) \cdot \nabla \phi_i(x) dx$.
If we let
$[\mathbf{G}_0]_{\p,\q} = \langle \Psi_\p(\y), \Psi_\q(\y) \rangle$ and $[\mathbf{G}_k]_{\p,\q} = \langle y_k \Psi_\p(\y), \Psi_\q(\y) \rangle$,
then $\mathbf{K}$ has a matrix representation, given by,
\begin{align}
\label{eq:SG_linear_kronecker_K_matrix}
\mathbf{K} = \mathbf{G}_0 \otimes \mathbf{A}_0 + \sum_{k=1}^N \mathbf{G}_k \otimes \mathbf{A}_k,
\end{align}
where $\mathbf{A} \otimes \mathbf{B}$ denotes the Kronecker product of $\mathbf{A}$ and $\mathbf{B}$.  
We note that a similar construction can be optained for any coefficient $a(x,\y)$ which can be written as a sum of separable functions of $x$ and $\y$, such as the polynomial function of Example \ref{exmp:polynomial_coefficient}. 

However, when $a(x,\y)$ is not separable in $x$ and $\y$, 
this construction is no longer valid, and the resulting matrix $\mathbf{K}$ may be block-dense if we simply carry out the Galerkin projections and compute $\mathbf{K}$ directly.
For certain special cases, e.g., when the diffusion coefficient is given by a log-transformed random field, the problem can be reformulated as a convection-diffusion problem, and the resulting system can be solved much more efficiently than the original problem \cite{Ullmann2012}. 
In general, this reformulation is not applicable.
Hence, for a general transcendental coefficient $a(x,\y)$, as given in Example \ref{exmp:transcendental_coefficient}, we project the coefficient onto an additional subspace $\mathcal{P}_{\Lambdar}(\Gamma)$, $r\in\N_0$, in order to obtain a separable representation.
To see this, define $\{\Psi_\rb(\y)\}_{0 \leq |\rb|}$ to be the (infinite) basis of orthonormal polynomials of $L^2_\varrho(\Gamma)$ as in \S \ref{subsec:stochastic_galerkin_approximation}. 
Then $a(x,\y)$ can be written as an expansion such that 
$a(x,\y) = \sum_{0 \leq |\rb|} a_\rb(x) \Psi_\rb(\y)$,
where the coefficients $a_\rb(x)= \langle a(x,\y), \Psi_\rb(\y) \rangle$.
Let $\Lambdar$ be an index set of the type described in \S \ref{subsec:stochastic_galerkin_approximation}.
Since infinite series representations are not practical in computations, we seek a truncation 
\begin{align}
\label{eq:a_r_truncation}
a^r(x,\y) := \sum_{\rb\in\Lambdar} a_\rb(x) \Psi_\rb(\y) 
\end{align}
in the subspace $\mathcal{P}_{\Lambdar}(\Gamma)$ for some $r\in\N_0$. 
When $a^r(x,\y) \neq a(x,\y)$, e.g., in the case that the projection order $r$ is chosen to minimize error independent of the SG discretization, we let $u_{h,p}^r$ denote the corresponding solution to the fully discrete SG approximation problem with $a(x,\y)$ replaced with $a^r(x,\y)$.
By substituting $a^r(x,\y)$ into \eqref{eq:deterministic_weak_problem_BF} we obtain
\begin{gather*}
\int_{D} \left(\sum_{\rb\in\Lambdar} a_\rb(x)\Psi_\rb(\y)\right) \; \nabla \phi_j(x) \cdot \nabla \phi_i(x) dx = \sum_{\rb\in\Lambdar} [\mathbf{A}_\rb]_{i,j} \Psi_\rb(\y), \numberthis \label{eq:operator_expansion} \\
\quad [\mathbf{A}_\rb]_{i,j} = \int_D a_\rb(x) \nabla \phi_j (x) \cdot \nabla \phi_i(x) dx. \numberthis \label{eq:A_kij}
\end{gather*}
Equation \eqref{eq:operator_expansion} represents an expansion of the stochastic finite element stiffness matrix $\mathbf{A}(\y)$ and equation \eqref{eq:A_kij} represents the $\rb$-th mode of the expansion. 
Let $\mathbf{u}_{h,\p}^r = [u^r_{1,\p},\ldots,u^r_{\Jh,\p}]^{\text{T}}$ denote the vector of nodal values of the finite element solution corresponding to the $\p$-th stochastic mode of $u_{h,p}^r$, and $\mathbf{u}_{h,p}^r = [\mathbf{u}_{h,\p}^{r}]_{\p\in\Lambdap}^{\text{T}}$. 
We substitute the expansion of $\mathbf{A}(\y)$ into the Galerkin equations \eqref{eq:fully_discrete_SG_algebraic_problem}, to obtain the coupled system: for each $\p\in\Lambdap$
\begin{align}
\label{eq:G_rpq}
\sum_{\rb\in\Lambdar} \sum_{\q\in\Lambdap} [\mathbf{G}_\rb]_{\p,\q}  \mathbf{A}_\rb \mathbf{u}_{h,\q}^r = \langle \Psi_\p, \mathbf{F} \rangle,  \qquad [\mathbf{G}_\rb]_{\p,\q} = \langle \Psi_{\p} \Psi_{\q} \Psi_\rb \rangle.
\end{align}
Alternatively, similar to \eqref{eq:SG_linear_kronecker_K_matrix}, we may define 
$\Kr = \sum_{\rb\in\Lambdar} \mathbf{G}_\rb\otimes\mathbf{A}_\rb$
to again obtain the coupled system of finite element problems: for all $\p\in\Lambdap$
\begin{align}
\label{eq:fully_discrete_SG_algebraic_sep_KuF}
\sum_{\q\in\Lambdap} [\Kr]_{\p,\q} \mathbf{u}_{h,\q}^r = \mathbf{F} \delta_{\bs{0},\p} \qquad \forall \p\in \Lambdap.
\end{align}
Note that in forming $\Kr$, we now need only store the matrices $\{\mathbf{G}_\rb\}_{\rb\in\Lambdar}$ and $\{\mathbf{A}_\rb\}_{\rb\in\Lambdar}$, so that the efficient, matrix-free, solution strategies discussed in Remark \ref{rem:K_storage_and_computation} can be applied. 

\begin{rem}[\rm Projection and well-posedness]
\label{rem:projection_well_posed}
In the case that the coefficient is a transcendental function of the random variables, as in Example \ref{exmp:transcendental_coefficient}, 
   there does not exist a $r\in\N_0$ such that the projection \eqref{eq:a_r_truncation} is exact. 
Due to the orthogonality of the basis, setting $r=2p$ in the construction of $\Kr$ yields an entirely block-dense system \cite{Matthies2005} that is equivalent to \eqref{eq:fully_discrete_SG_algebraic_KuF}, and
computationally infeasible to solve.
A more practical approach is to choose the expansion order $0\leq r\leq 2p$, based on a-priori estimates of the error in the solution introduced by the truncation, so that the error when using the truncated expansion does not exceed that of the SG approximation. 
In this approach however, it becomes important to consider whether the truncated projection violates the well-posedness of \eqref{eq:model_problem} by failing to satisfy assumption $\Aone$.
One way to guarantee this is to choose $\tilde{r} \leq r \leq 2 p$ such that
\begin{align}
\label{eq:r_tilde}
\tilde{r} := \min \{ r\in\N_0 : \|a - a^\nu \|_{L^\infty_\varrho(\Gamma;L^\infty(D))} \leq a_{\min}, \; \forall \nu\in\N_0, \; \nu \geq r \}. 
\end{align}
An example of this problem can be seen in Figure \ref{fig:ar_uniform_ellipticity_regions} where for the function $a(x,y) = 0.1 + \exp(2.5 y)$, uniform ellipticity of the truncated projection $a^r(x,y)$ does not hold on $\Gamma = [-1,1]$ for $r < 4$. 

\begin{figure}[!htb]
\begin{center}
\includegraphics[width=0.243\textwidth]{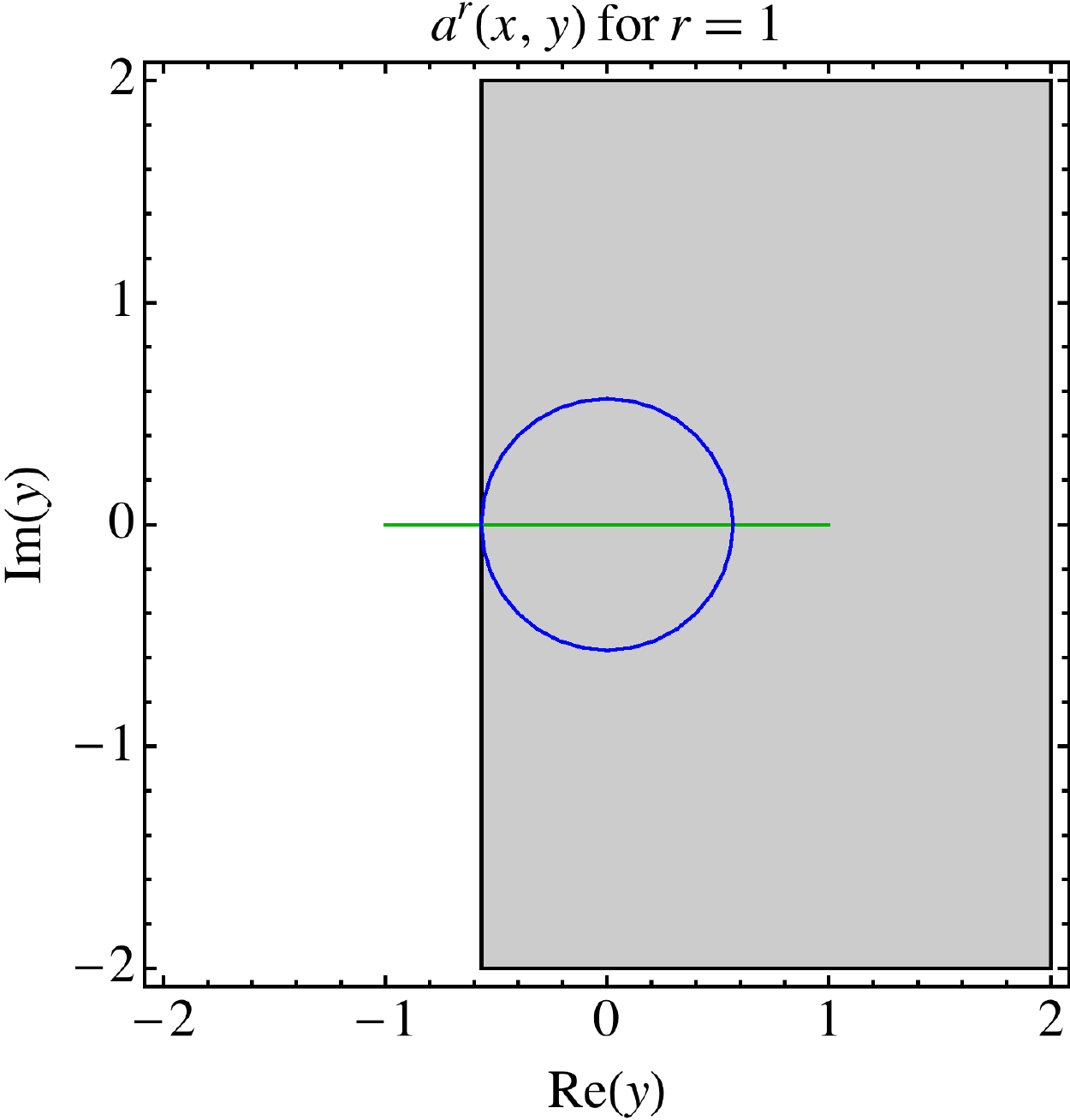}
\includegraphics[width=0.243\textwidth]{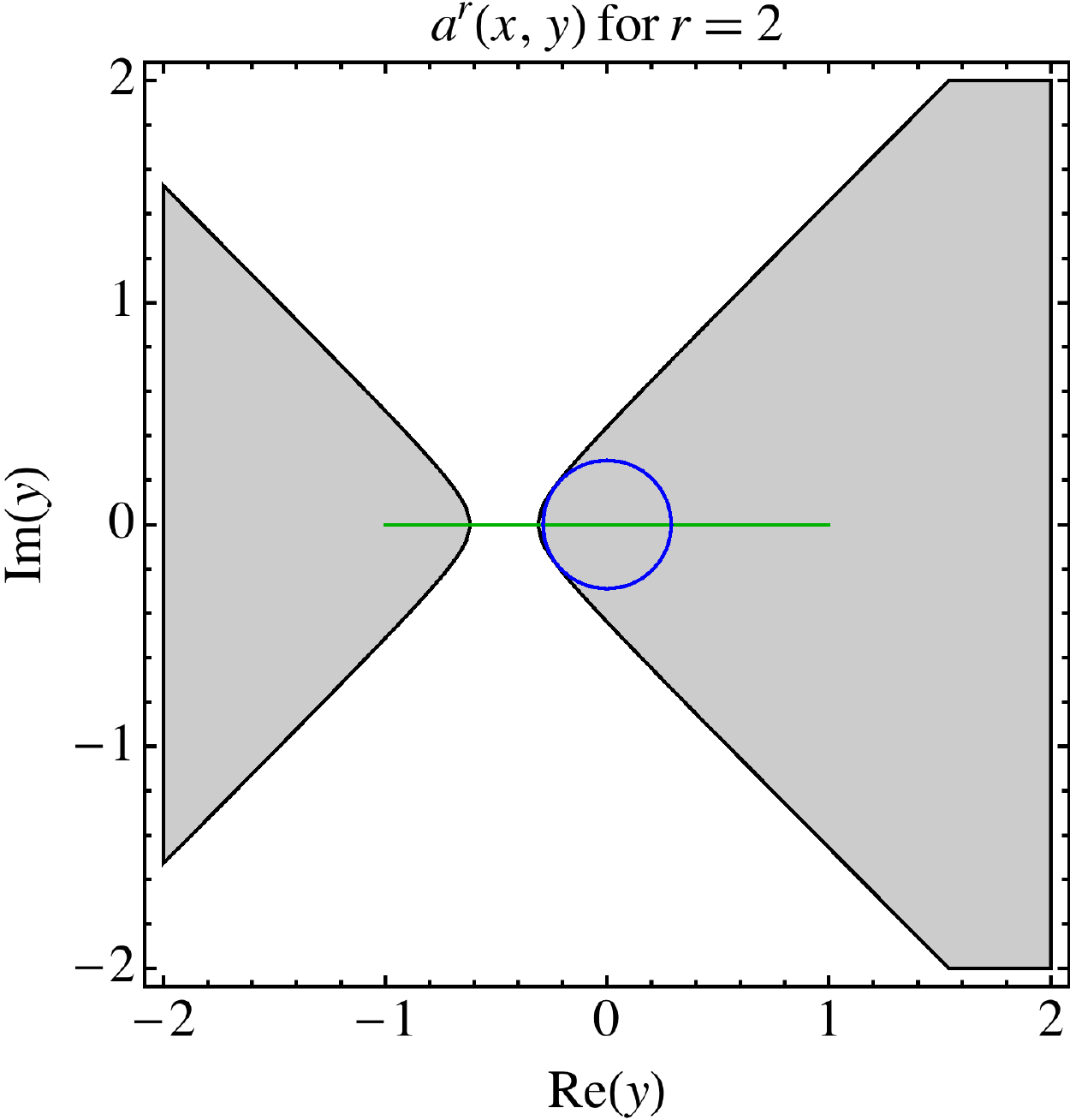}
\includegraphics[width=0.243\textwidth]{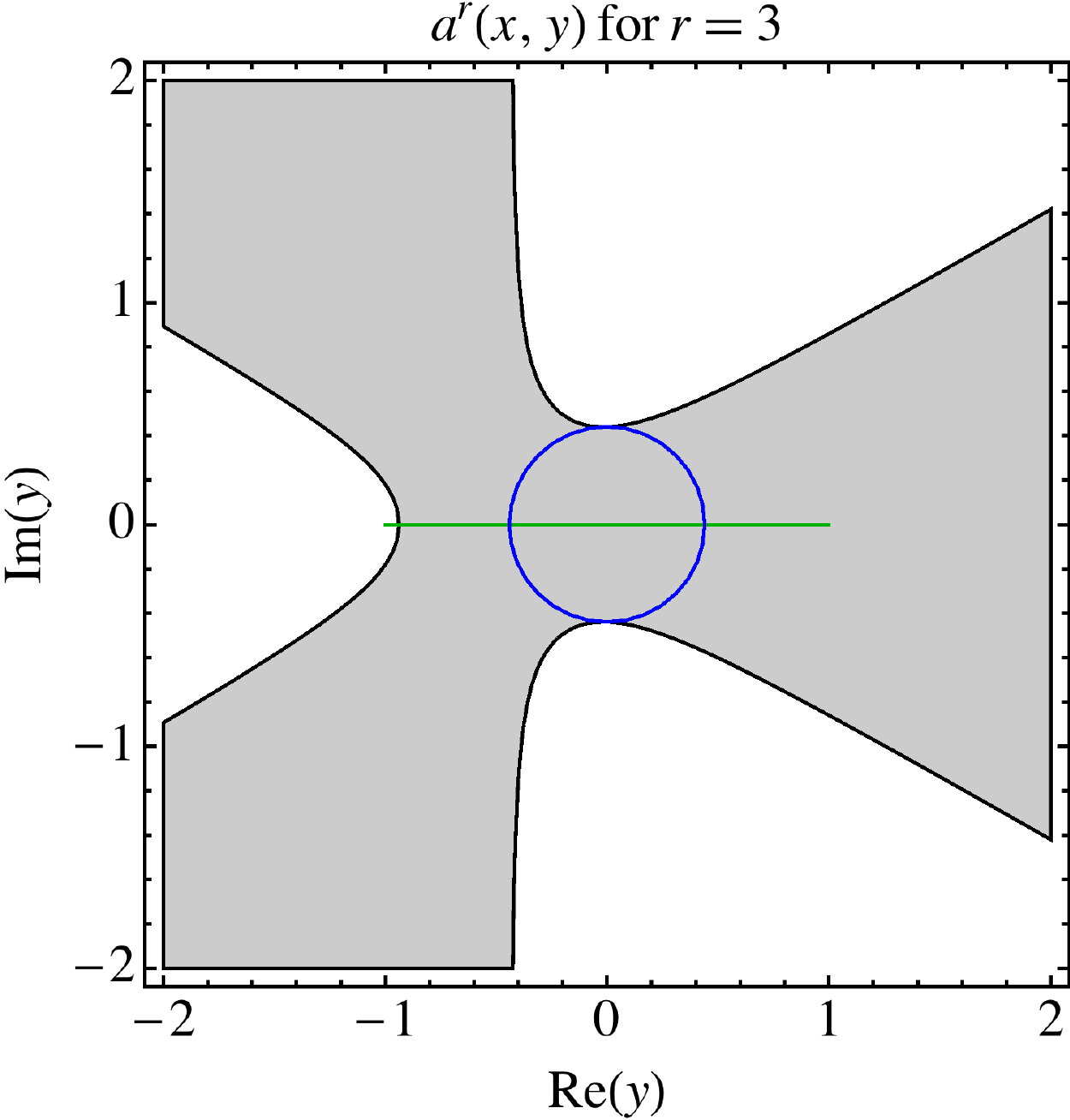}
\includegraphics[width=0.243\textwidth]{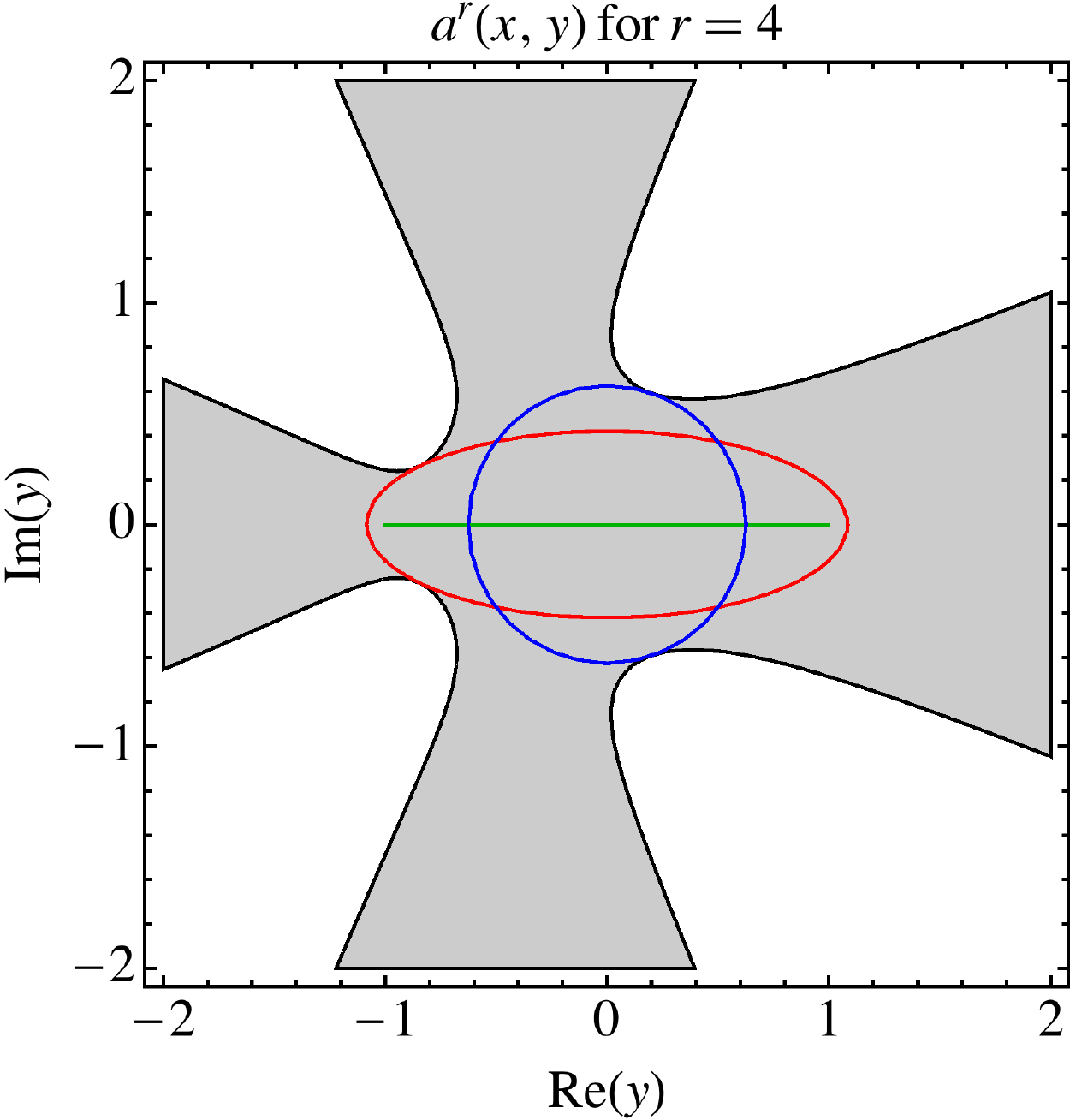}
\includegraphics[width=0.243\textwidth]{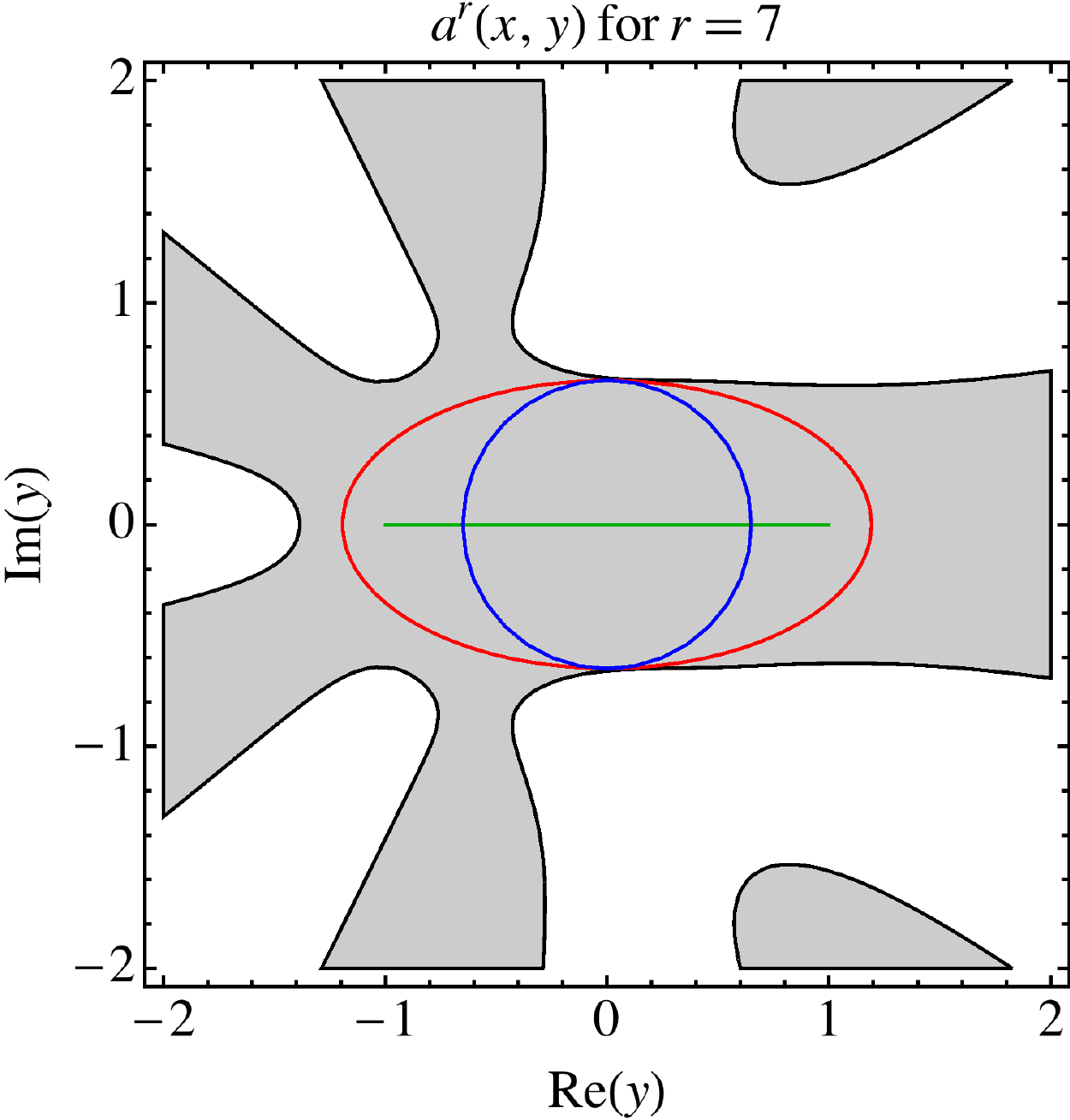}
\includegraphics[width=0.243\textwidth]{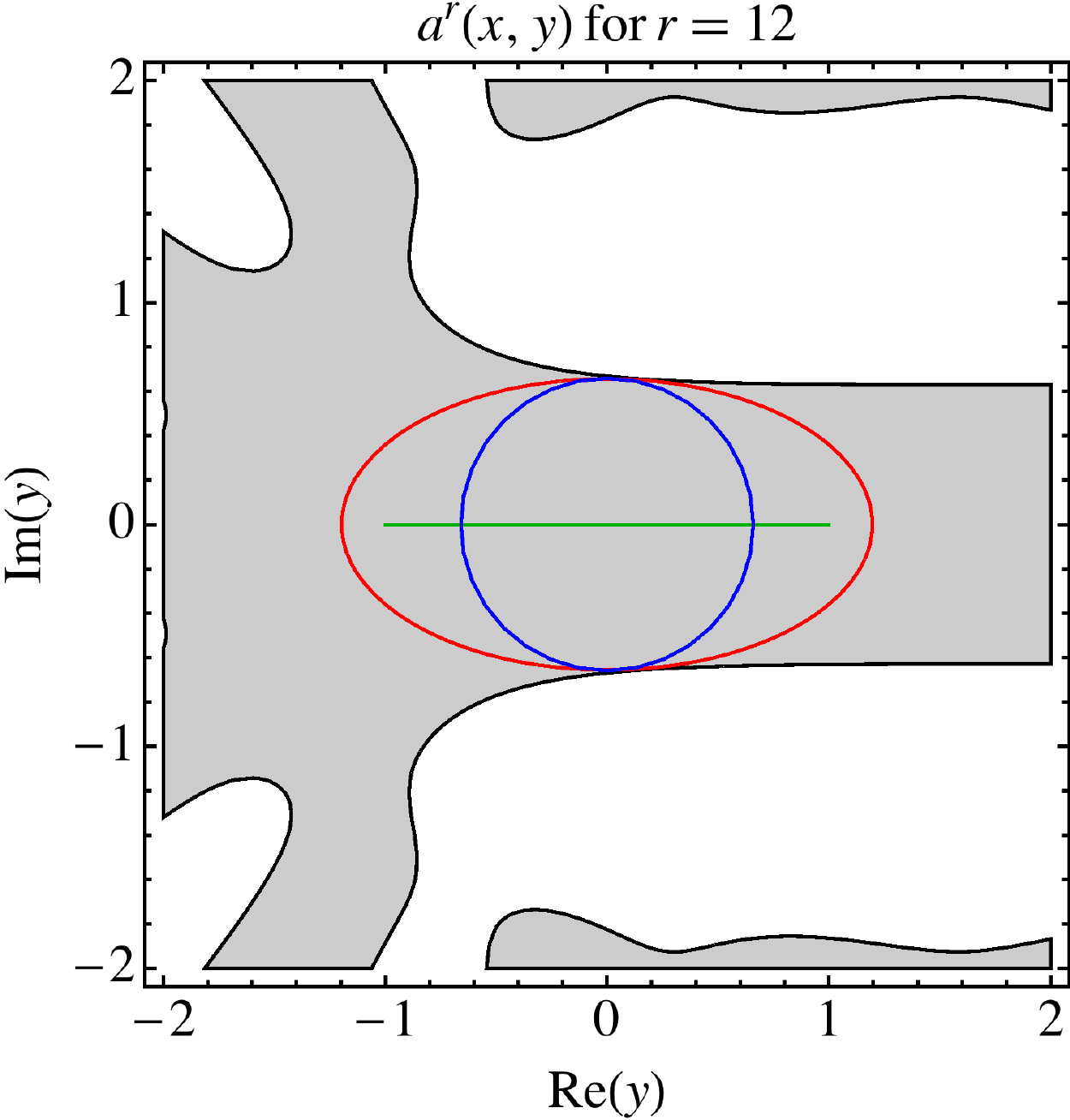}
\includegraphics[width=0.243\textwidth]{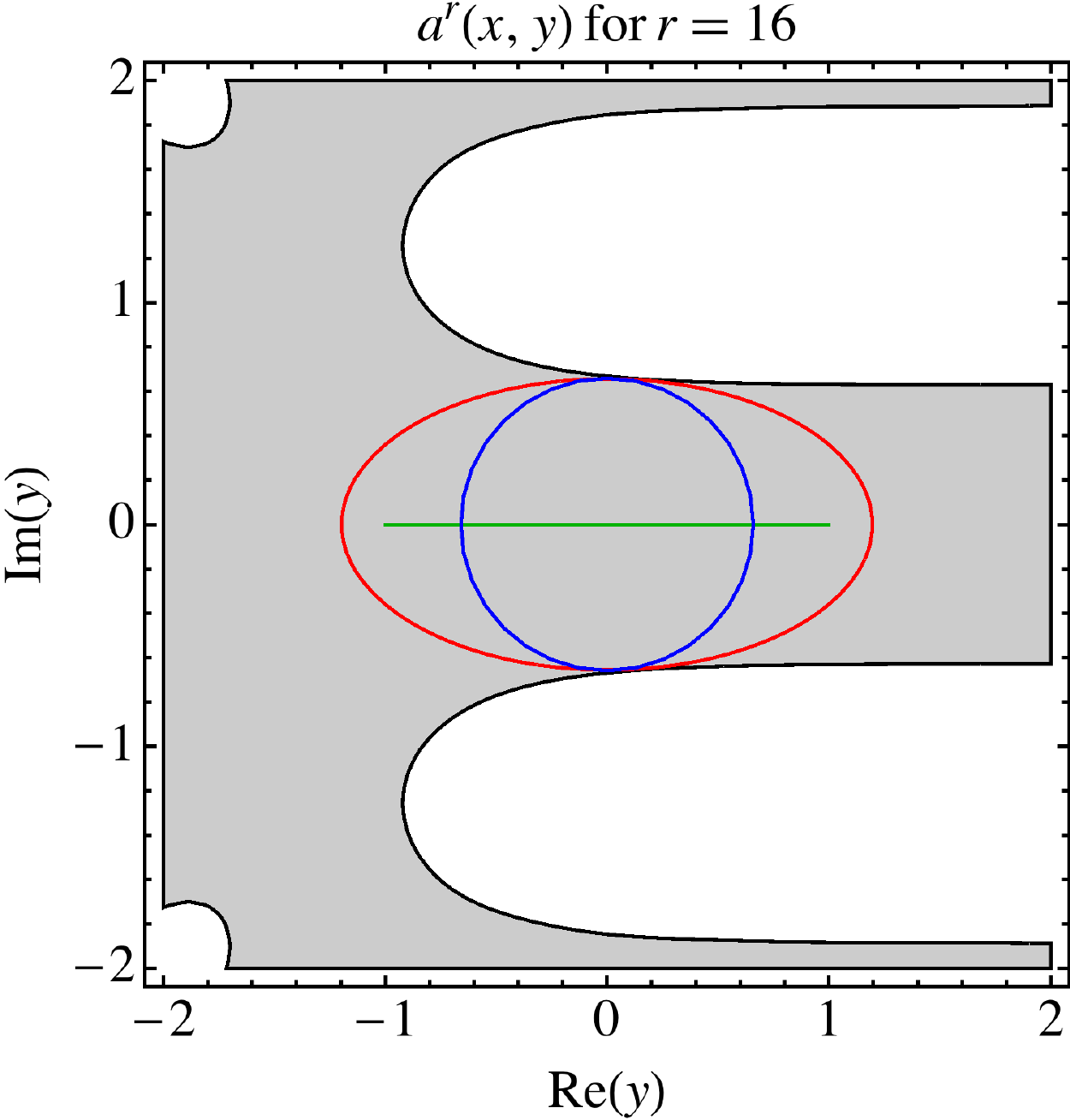}
\includegraphics[width=0.243\textwidth]{exponential_rv.pdf}
\caption{Domains of uniform ellipticity for the total degree orthogonal expansions of order $r$ of the one-dimensional coefficient $a(x,y)=0.1+\exp(2.5 y)$, for $y\in\Gamma = [-1,1]\subset \R^1$, are indicated by the gray regions in each plot. The last plot shows the domain of uniform ellipticity of the original function $a(x,y)$. The blue and red curves represent the maximal discs and ellipses, respectively, that can be contained in those domains, and the green lines represent the interval $\Gamma $.}
\label{fig:ar_uniform_ellipticity_regions}
\end{center}
\end{figure}
\end{rem}

\subsection{Cost of solving the generalized SG system}
\label{subsec:solving_SG}

Without loss of generality, to solve the stochastic Galerkin system \eqref{eq:fully_discrete_SG_algebraic_sep_KuF} for $\mathbf{u}_{h,p}^r$, we use the precondtioned conjugate gradient (PCG) method, 
wherein, for the unpreconditioned CG method, we have the estimate
\begin{align}
\label{eq:CG_error_SG}
\|\mathbf{u}^{r}_{h,p} - \mathbf{u}^{r,(k)}_{h,p} \|_{\Kr} \leq 2 \left( \frac{\sqrt{\kappar} - 1}{\sqrt{\kappar} + 1} \right)^k \|\mathbf{u}^{r}_{h,p} - \mathbf{u}^{r,(0)}_{h,p} \|_{\Kr}.
\end{align}
Here $\kappar$ is the condition number of $\Kr$, $\mathbf{u}^{r,(0)}_{h,p}$ is the vector of the initial guess, and $\mathbf{u}^{r,(k)}_{h,p}$ is the output of the $k$-th iteration of the CG solver.
The CG method is highly dependent on the conditioning of the system, and when $\kappar$ is large, the number of iterations needed to reduce the error in $\mathbf{u}^{r,(k)}_{h,p}$ will also be significant. Hence we introduce the mean-based block-diagonal preconditioner (see, e.g., \cite{Pellissetti00iterativesolution,Powell2009}), 
\begin{align}
\label{eq:mean_based_preconditioner}
\mathbf{P} := \mathbf{G}_\bs{0} \otimes \mathbf{A}_\bs{0},
\end{align}
with $\mathbf{A}_\bs{0}$ and $\mathbf{G}_\bs{0}$ the matrices defined in \eqref{eq:A_kij} and \eqref{eq:G_rpq} for $\rb = \bs{0}$, respectively.

For $\rb \in \Lambdar$, at every iteration of the CG method, or any iterative approach, each nonzero entry in each matrix $\mathbf{G}_\rb$ implies a matrix-vector product of the form 
$\langle \Psi_\p\Psi_\q\Psi_\rb \rangle \mathbf{A}_\rb \mathbf{p}^{(k)}_\q$, 
where $\langle \Psi_\p\Psi_\q\Psi_\rb \rangle$ is a scalar quantity. 
Let $\text{nnz}(\mathbf{A})$ denote the number of nonzeros of a matrix $\mathbf{A}$, and define
\begin{align}
\label{eq:nnz_G}
\Mpr = \sum_{\rb\in\Lambdar} \text{nnz}(\mathbf{G}_\rb) 
\end{align}
to be the total number of nonzeros in all of the matrices $\{\mathbf{G}_\rb\}_{\rb\in\Lambdar}$ at order $p$. With this in mind, an upper bound for the work in floating point operations (FLOPs) of solving \eqref{eq:fully_discrete_SG_algebraic_sep_KuF} is given by
\begin{align}
\label{eq:SG_cost_solve}
\WSGsolve \approx \mathcal{O}(\Jh) * \Mpr * \NiterSG, 
\end{align}
where the term $\mathcal{O}(\Jh)$ corresponds to the cost of a single finite element matrix-vector product, and $\NiterSG$ is the number of iterations of the CG solver without a preconditioner.
If we apply a preconditioner, in hopes to minimize $\NiterSG$, we must also account for the added cost of applying the preconditioner at each iteration. With the mean-based preconditioner from \eqref{eq:mean_based_preconditioner}, at each iteration we multiply an additional matrix of size $\Jh \Mp \times \Jh \Mp$, but the matrix consists only of $\Mp$ diagonal blocks. 
Here we assume that in finding the inverse of the preconditioning matrix $\mathbf{P}$, a sparse approximation to $\mathbf{P}^{-1}$ is used, which we will denote $\mathbf{\tilde P}^{-1}$. Such a decomposition can be found from, e.g., incomplete LU or incomplete Cholesky factorizations. 
Hence, for each iteration we require $\Mp$ additional matrix-vector products of the size, and complexity, of the original finite element system, so the work estimate in FLOPs for the case of this preconditioner is given by 
\begin{align}
\label{eq:SG_cost_solve_precon}
\WSGsolvep \approx \mathcal{O} (\Jh) * ( \Mp + \Mpr) * \NiterpreSG , 
\end{align}
where $\NiterpreSG$ is the number of iterations needed by the PCG method.
Other preconditioners, such as the Kronecker product preconditioner suggested in \cite{Ullmann2010} would require a different form of \eqref{eq:SG_cost_solve_precon}. 

\begin{figure}[!htb]
\begin{center}
\includegraphics[width=0.161\textwidth]{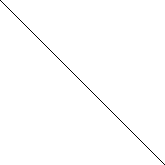}
\includegraphics[width=0.161\textwidth]{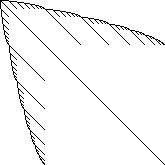}
\includegraphics[width=0.161\textwidth]{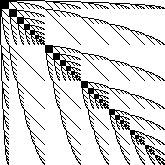}
\includegraphics[width=0.161\textwidth]{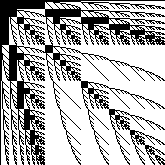}
\includegraphics[width=0.161\textwidth]{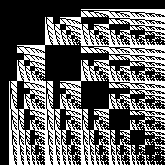}
\includegraphics[width=0.161\textwidth]{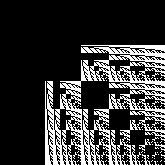}
\caption{Visualization of the number of nonzeros of a $165\times165$ SG matrix with elements $[\Kr]_{\p,\q} = \sum_{\rb\in\Lambdar} [\mathbf{G}_\rb]_{\p,\q} * \mathbf{A}_\rb$. Each pixel represents a block finite element system when using a total degree projection of the solution of fixed degreee $p=3$, and increasing the total degree of the projection of the coefficient, i.e., $r=0,1,2,3,4,5$. At $r = 6$, the matrix is entirely block-dense.}  
\label{fig:TD_sparsity}
\end{center}
\end{figure}

Figure \ref{fig:TD_sparsity} displays the effect of fixing the projection order of the solution but increasing the order of the projection of the coefficient. In order to minimize the error of the projection, such a situation would be required if the coefficient is highly nonlinear and reflects the importance of considering $\Mpr$ in the computational cost of the SGFEM.

\section{Explicit cost bounds for the SGFEM} 
\label{sec:complexity_analysis}

The primary goal of this section is to estimate the algorithmic complexity required by the SGFEM to construct an approximation 
to \eqref{eq:model_problem} within a prescribed tolerance $\eps>0$. 
We assume $a(x,\y)$ is a general non-affine coefficient, as in Examples \ref{exmp:polynomial_coefficient} and \ref{exmp:transcendental_coefficient}, satisfying assumptions $\Aone$ and $\Atwo$.
Let $a^r(x,\y)$ be the orthogonal expansion of $a(x,\y)$, given by \eqref{eq:a_r_truncation}, of total degree $r$, i.e., $a^r (x,\y)\in \mathcal{P}_{\Lambdar}(\Gamma)$ with $\Lambdar=\LambdarTD$ from \eqref{eq:index_set_examples}.
We further assume that $\tilde{r} \leq r \leq 2p$, with $\tilde{r}$ given in \eqref{eq:r_tilde}, so that $a^r(x,\y)$ also satisfies $\Aone$ and $\Atwo$.
We will focus on the complexity of solving \eqref{eq:fully_discrete_SG_algebraic_sep_KuF}, when the stochastic discretization to \eqref{eq:model_problem} 
is performed in 
$\mathcal{P}_{\Lambdap}(\Gamma)$ with $\Lambdap = \LambdapTD$ from \eqref{eq:index_set_examples}, i.e., in the space of total degree polynomials of order $p$, and the physical discretization is performed with the finite element method. 
These results are presented in the context of solving the linear system \eqref{eq:fully_discrete_SG_algebraic_sep_KuF} with a PCG method when a zero initial vector is used to seed the solver. 
The results, however, can be generalized to other methods, such as preconditioned GMRES and other Krylov subspace methods. 

The results are organized as follows. In \S \ref{subsec:SG_mat_vec_complexity} we discuss the overall complexity of the matrix-vector products associated with solving \eqref{eq:fully_discrete_SG_algebraic_sep_KuF} when using the SG matrix $\Kr = \sum_{\rb\in\Lambdar} \mathbf{G}_{\rb} \otimes \mathbf{A}_{\rb}$. Our analysis extends the results of \cite{Ernst2010} in order 
to provide a bound on the block-sparsity of the SG system $\Kr$ in the more general setting of a non-affine coefficient $a^r(x,\y)$, as given in Examples \ref{exmp:polynomial_coefficient} and \ref{exmp:transcendental_coefficient}.
In particular, for $\mathbf{G}_\rb$ given in \eqref{eq:G_rpq}, we show that $\text{nnz}(\mathbf{G}_\rb) = \mathcal{O}(\min\{2^{|\rb|}, M_{\lceil |\rb|/2 \rceil}\} \Mpsubhalfrb)$ for every $\rb\in\Lambdar$, where $M_{r} = {N+r \choose N}$ for $r\in\N_0$, 
when solving \eqref{eq:fully_discrete_SG_algebraic_sep_KuF}, so that the total complexity of the matrix-vector products with the Galerkin system is $\mathcal{O}(\Jh \Mp \Mr \min\{ 2^r , M_{\lceil r/2 \rceil}\} )$.
In \S \ref{subsec:eps_complexity_analysis}, we 
perform an $\eps$-complexity analysis to derive the explicit cost bounds of the SGFEM using PCG, in terms of FLOPs as the tolerance $\eps\to0$.
Finally, in \S\ref{subsec:SG_condition_number} we discuss issues related to the conditioning of the SG system. 

\subsection{Complexity of matrix-vector multiplications in the SG approximation}
\label{subsec:SG_mat_vec_complexity}

In this section we provide rigorous 
bounds 
on the sparsity of the SG matrix $\mathbf{K}_{r}$ from \eqref{eq:fully_discrete_SG_algebraic_sep_KuF}, for arbitrary $0 \leq r \leq 2p$, and $p\in\N_0$. 
Our main result, given by Theorem \ref{thm:nnz_Gk_general},
provides an exact count for $\text{nnz}(\mathbf{G}_{\rb})$ in the general case $|\rb|\in\N_0$ and $N\in\N$ when the integrals $[\mathbf{G}_\rb]_{\p,\q} = \langle \Psi_\rb \Psi_\p \Psi_\q \rangle$ are defined in terms of even weight functions $\varrho$. This result can be seen as an extension of estimates from \cite{Ernst2010},
where bounds on the sparsity of 
$\mathbf{G}_\rb$ were shown in the cases (i) $|\rb|= 1$ and $N\in\N$ and (ii) $N=1$ and $\rb=r\in\N_0$. 
We then provide upper bounds on the total number of nonzeros blocks $\Mpr = \sum_{\rb\in\Lambdar} \text{nnz}(\mathbf{G}_\rb)$ of the matrix $\Kr$ from \eqref{eq:fully_discrete_SG_algebraic_sep_KuF}, both in the cases that $a(x,\y)$ is a finite order polynomial as in Example \ref{exmp:polynomial_coefficient}
and the case that $a(x,\y)$ is a transcendental function of the random variables, as in Example \ref{exmp:transcendental_coefficient}.
Our first major result is summarized in the following Theorem: 

\begin{thm}
\label{thm:nnz_Gk_general}
Let $\Lambdap$ and $\Lambdar$ be the isotropic total degree index sets corresponding to the solution and coefficient, respectively, with $p,r\in\N_0$, and $0 \leq r \leq 2p$. 
If $\rb \in \Lambdar$,
and $\varrho_{i}$ are even for all $i=1,\ldots,N$, then for the matrix $\mathbf{G}_\rb$ from \eqref{eq:G_rpq} we have
\begin{align}
\label{eq:nnz_Gk_general}
\textnormal{nnz}(\mathbf{G}_\rb) & = \sum_{\ell = \lceil |\rb|/2 \rceil}^{|\rb|} c(\rb,\ell) { N + p - \ell \choose p - \ell }, \qquad
c(\rb,\ell) = \begin{cases} \#\bm{S}(\rb,\ell) & \text{$|\rb|$ even, $\ell = |\rb|/2$,} \\
2 \#\bm{S}(\rb,\ell) & \text{otherwise,}
\end{cases} 
\end{align}
with $\bm{S}(\rb,\ell) = \left\{ \bm{s}\in\N^N_0 :  |\sbb| = \ell, \; \sbb \leq \rb \right\}$,
so that $\#\bm{S}(\rb,\ell)$ is equal to the coefficient of $t^\ell$ in the polynomial 
$P_\rb (t) = \prod_{i=1}^N \sum_{j=0}^{r_i} t^j.$
Moreover, we have the following bound for $\textnormal{nnz}(\mathbf{G}_\rb)$, i.e.,
\begin{align*}
\textnormal{nnz}(\mathbf{G}_\rb) \leq 2 \min\left\{ 2^{|\rb|}, { N + \lceil |\rb|/2 \rceil \choose N } \right\} { N + p - \lceil |\rb|/2 \rceil \choose N }, \numberthis \label{eq:nnz_Gk_general_bound}
\end{align*}
so that 
\begin{align*}
\Mpr 
    \leq 2 \sum_{j=0}^r \min\left\{ 2^j, { N + \lceil j/2 \rceil \choose N } \right\} { N - 1 + j \choose N - 1} { N + p - \lceil j/2 \rceil \choose N }. \numberthis \label{eq:Mpr_general_bound} 
\end{align*}
\end{thm}

\begin{pf}
For a given $\rb\in\Lambdar$, we estimate the number of pairs $(\p,\q)\in \Lambdap\times \Lambdap$ such that $\langle \Psi_\rb \Psi_\p \Psi_\q \rangle = \prod_{i=1}^N \langle \psi_{r_i} \psi_{p_i} \psi_{q_i} \rangle \neq 0$. To do this, we extend the result of \cite[Lemma 28]{Ernst2010} to a general matrix $\mathbf{G}_\rb$ with $|\rb|\in\N$. Since $\{\Psi_\rb\}_{\rb\in\Lambdar}$
are
orthonormal with respect to the even 
weight function $\rho(\y) = \prod_{i=1}^N \rho_i(y_i)$, we see that $\langle \Psi_\rb \Psi_\p \Psi_\q \rangle \neq 0$ if and only if $(\p,\q)\in\Theta_\rb$, where 
\begin{align*}
\Theta_\rb = \{(\p,\q) \in \Lambdap \times \Lambdap : |p_i - q_i| \leq r_i \leq p_i + q_i, \\
                         \text{ and } p_i + q_i + r_i \text{ is even } \forall i = 1,\ldots,N \}.
\end{align*}
Therefore, to estimate the number of nonzeros in the matrix $\mathbf{G}_\rb$, we must estimate $\#\Theta_\rb$. However, $\Theta_\rb$ is different for each $\rb\in\Lambdar$. Even when $\rb_1,\rb_2\in\Lambdar$ are such that $|\rb_1| = |\rb_2|$, in general we do not have that $\#\Theta_{\rb_1}=\#\Theta_{\rb_2}$. 
On the other hand, if $\rb_2$ is a permutation of $\rb_1$, then it is easy to see that $\#\Theta_{\rb_1}=\#\Theta_{\rb_2}$ since $\Lambdap$ is the isotropic total degree set.
Also note that $\langle \Psi_\rb \Psi_\p \Psi_\q \rangle = \langle \Psi_\rb \Psi_\q \Psi_\p \rangle$ so if $(\p,\q)\in\Theta_\rb$ then $(\q,\p)\in\Theta_\rb$ as well. Note that $\Theta_\rb$ can be rewritten
\begin{align*}
\Theta_\rb = \{(\p,\q) \in \Lambdap \times \Lambdap : |p_i - q_i| \leq r_i \leq p_i + q_i, \\
                         \text{ and } |p_i - q_i| + r_i \text{ is even } \forall i = 1,\ldots,N \}.
\end{align*}
Hence if $(\p,\q) \in \Theta_\rb$, then we see that $(\p,\q)$ must satisfy 
\begin{list}{(\roman{qcounter})~}{\usecounter{qcounter}}
\item $|p_i - q_i| \leq r_i$ for all $i=1,\ldots,N$, 
\item $r_i \leq p_i + q_i$ for all $i=1,\ldots,N$, and 
\item $|p_i - q_i| + r_i$ is even for all $i=1,\ldots,N$.
\end{list}
Note that when $r_i=0$, we see that $p_i=q_i \leq p$, and when $r_i > 0$ we see that (i) and (iii) imply 
$|p_i - q_i| \in \{ 0, 2, 4, \ldots, r_i \}$ 
for 
$r_i$ even, and
$|p_i - q_i| \in \{ 1, 3, 5, \ldots, r_i \}$ 
for 
$r_i$ odd.
For each $i =1,\ldots, N$, let $\{k_i^{(n)}\}_{n=0}^{\lfloor r_i/2 \rfloor}$ be the sequence defined by 
\begin{align*}
k_i^{(n)} = \begin{cases} 2n +1 & r_i \text{ odd}, \\ 2n & r_i \text{ even}, \end{cases}
\end{align*}
so that fixing $|p_i - q_i| = k_i^{(n)}$ implies 
that conditions (i) and (iii) are met.

To satisfy (ii) we must have $r_i \leq p_i + q_i$ and to satisfy (i) and (iii) we must have $|p_i - q_i| = k_i^{(n)}$. 
To avoid overcounting due to symmetry, we first fix possible values of $p_i$ and consider what $q_i$ must be. 
Let $\{s_i^{(n)}\}_{n=0}^{\lfloor r_i/2 \rfloor}$ be the sequence defined by
\begin{align*}
s_i^{(n)} = \frac{r_i + k_i^{(n)}}{2},
\end{align*}
which we will refer to as the sequence of {\em starting points} for $p_i$ corresponding to $k_i^{(n)}$. 
Note that the starting points $\{s_i^{(n)}\}_{n=0}^{\lfloor r_i/2 \rfloor}$ enumerate the integers between $\lceil r_i/2 \rceil$ and $r_i$.
Picking $p_i \in \{s_i^{(n)},s_i^{(n)}+1,\ldots,p\}$ and $q_i = p_i - k_i^{(n)}$ we have 
\begin{align*}
p_i + q_i = 2p_i - k_i^{(n)} \geq 2 s_i^{(n)} - k_i^{(n)}  = 2 \left(\frac{r_i + k_i^{(n)}}{2} \right) - k_i^{(n)} = r_i
\end{align*}
or $p_i + q_i \geq r_i$, so that (ii) is satisfied.

Since (i), (ii), and (iii) are 
satisfied by setting $p_i\in \{s_i^{(n)},s_i^{(n)}+1,\ldots,p\}$ and $q_i = p_i - k_i^{(n)}$ for a fixed $0 \leq n \leq \lfloor r_i/2 \rfloor$, we count the number of admissible pairs for these choices. 
In $N-1$ variables, the number of polynomials of total degree less than or equal to $p-p_i$ is given by
\begin{align*}
{ N - 1 + p - p_i \choose p - p_i },
\end{align*}
where ${ n \choose k } = 0$ if $n<k$ or $n,k< 0$.
To simplify notation, pick $s_i = s_{i}^{(n)}$ (one of the starting points in the $i$-th direction) and $k_i = k_i^{(n)}$ (its associated distance), where $0 \leq n \leq \lfloor r_i/2 \rfloor$ is fixed. To count the number of admissible pairs associated with the difference $k_i$ and starting point $s_i$, we compute 
\begin{align*}
\sum_{p_i = s_i}^p { N - 1 + p - p_i \choose p - p_i } = \sum_{j = 0}^{p-s_i} { N - 1 + j \choose j } = { N + p - s_i \choose p - s_i }.
\end{align*}
Define $\sbb \in\N^N_0$ with the $s_i$ as above, then $\sbb$ corresponds to a possible combination of starting points in each direction. 
To estimate the number of polynomials associated with the starting point $\sbb$, we compute 
\begin{align}
\label{eq:polynomial_bincoeff_sum}
\sum_{p_1 = s_1}^p \sum_{p_2 = s_2}^{p-p_1} \cdots \sum_{p_N = s_N}^{p-p_1-\cdots-p_{N-1}} { p - p_1 - \cdots - p_N \choose p - p_1 - \cdots - p_N } = { N + p - |\sbb| \choose p - |\sbb| },
\end{align}
where the sum easily follows by an induction argument and Pascal's rule.

Enumerating all of the pairs $(\p,\q)\in\Theta_\rb$ thus reduces to counting the number of possible combinations of starting points. Hence, in $N$ dimensions we consider all such multi-indices of the $\{s_i^{(n)}\}_{n=0}^{\lfloor r_i/2 \rfloor}$ whose components sum to some integer $\lceil |\rb|/2 \rceil \leq \ell \leq |\rb|$. 
For two multi-indicies $\sbb,\rb \in \N^N_0$, we say $\sbb \leq \rb$ if and only if $s_i \leq r_i$ for all $i=1,\ldots,N$. 
Define the set 
$\bm{S}(\rb,\ell) = \left\{ \bm{s}\in\N^N_0 :  |\sbb| = \ell, \; \sbb \leq \rb \right\}$,
which corresponds to a particular slice of the desired set of starting points.
To estimate $\#\bm{S}(\rb,\ell)$, we consider the familiar counting argument of placing $N$ bars among $\ell$ stars with the added restriction that the numer of stars in the $i$-th bin not exceed $r_i$. Such a problem can be reframed in terms of finding the coefficient $c(\rb,\ell)$ of $t^{\ell}$ in the generating function
$P_\rb (t) = \prod_{i=1}^N \sum_{j=0}^{r_i} t^j$.
Combining \eqref{eq:polynomial_bincoeff_sum} and summing over $\ell$ between $\lceil |\rb|/2 \rceil \leq \ell \leq |\rb|$ we arrive at \eqref{eq:nnz_Gk_general}, where the coefficients $c(\rb,\ell) = \#\bm{S}(\rb,\ell)$ when $|\rb|$ is even and $\ell = |\rb|/2$ (in this case the roles of $p_i$ and $q_i$ can not be reversed) and $c(\rb,\ell) = 2\#\bm{S}(\rb,\ell)$ otherwise. 

Noting that 
$\cup_{\ell = \lceil |\rb|/2 \rceil}^{|\rb|} \bm{S}(\rb,\ell)$ is a change of coordinates of a total degree index set of order $\lceil |\rb|/2 \rceil$ intersected with the hyperrectangle $\{\sbb \in\N^{N}_0 : \sbb \leq \rb\}$ yields the bound 
\begin{equation*}
\sum_{\ell = \lceil |\rb|/2 \rceil}^{|\rb|} c(\rb,\ell) \leq 2 \sum_{\ell=\lceil |\rb|/2 \rceil}^{|\rb|} \#\bm{S}(\rb,\ell) \leq 2 { N + \lceil |\rb|/2 \rceil \choose N }. 
\end{equation*}
On the other hand, from the generating function $P_\rb(t)$ we see that $c(\rb,\ell)$ is bounded by ${ |\rb| \choose \ell }$ when $|\rb|$ is even and $\ell = |\rb|/2$ and $2{|\rb| \choose \ell}$ otherwise. This follows from the fact that when $\bm{k}$ is the multi-index having $|\rb|$ ones and the rest zeros, since $\ell \leq |\rb|$, we have that $\#\bm{S}(\rb,\ell) \leq \#\bm{T}(\kb,\ell)$ where $\bm{T}(\kb,\ell) = \{ \sbb\in\N^N_0 : |\sbb| = \ell, \sbb \leq \kb \}$ and $\#\bm{T}(\kb,\ell)$ is given by the coefficient of $t^\ell$ in $P_{\kb}(t) = (1+t)^{|\rb|} = \sum_{\ell = 0}^{|\rb|} {|\rb| \choose \ell} t^\ell$ from the binomial theorem.
Then
\begin{align*}
\sum_{\ell = \lceil |\rb|/2 \rceil}^{|\rb|} c(\rb,\ell) \leq 2 \sum_{\ell = 0}^{|\rb|} c(\rb,\ell) = 2^{|\rb|+1},
\end{align*}
so that 
\begin{align*}
\allowdisplaybreaks
\text{nnz}(\mathbf{G}_\rb) = \sum_{\ell = \lceil |\rb|/2 \rceil}^{|\rb|} c(\rb,\ell) { N + p - \ell \choose p - \ell } 
      \leq 2 \min \left\{ 2^{|\rb|}, { N + \lceil |\rb|/2 \rceil \choose N }\right\} { N + p - \lceil|\rb|/2\rceil \choose N } ,
\end{align*}
showing \eqref{eq:nnz_Gk_general_bound}.
Substituting \eqref{eq:nnz_Gk_general_bound} into \eqref{eq:nnz_G} shows the bound of
$\Mpr$ from \eqref{eq:Mpr_general_bound}. \qed
\end{pf}

\begin{figure}[!htb]
\begin{center}
\includegraphics[width=5.012cm,clip=true,trim=31mm 64mm 44mm 59mm]{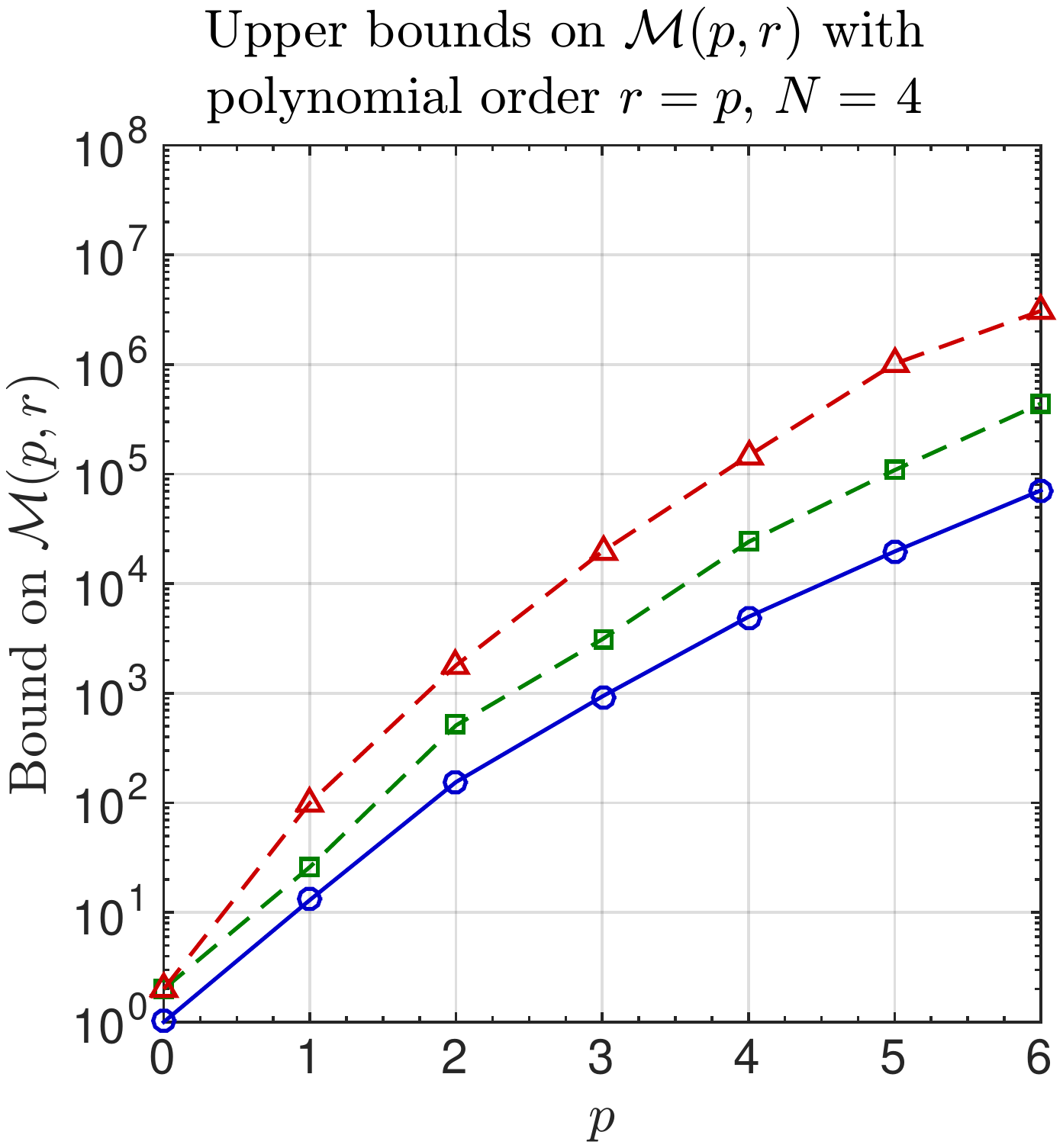}
\includegraphics[width=5.120cm,clip=true,trim=28mm 64mm 44mm 59mm]{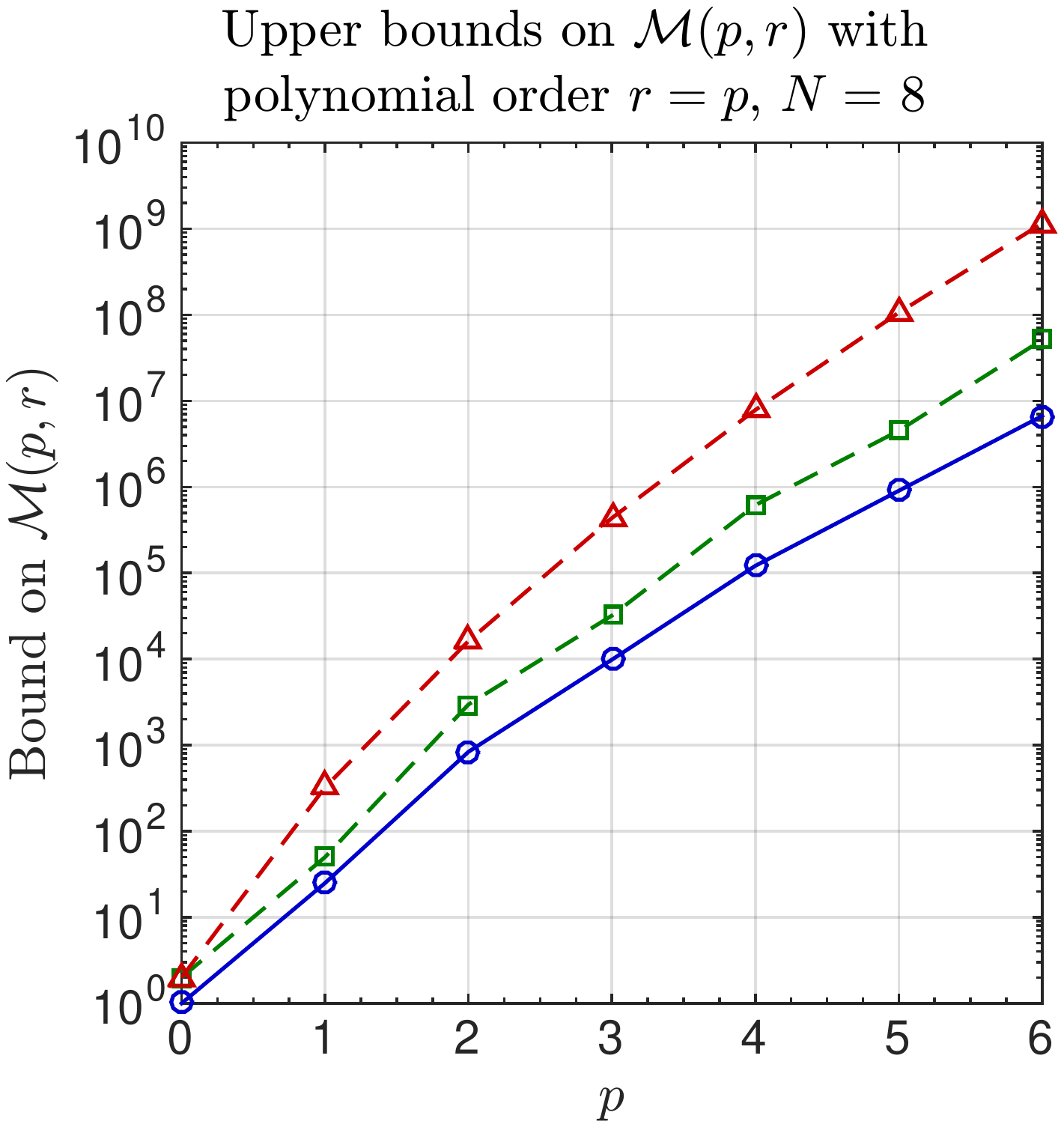}
\caption{For $r=p$ with $p$ ranging from $0,1,\ldots,6$ we plot for $N=4$ (left) and $N=8$ (right) the actual sparsity $\Mpr$ given by \eqref{eq:nnz_G} of the Galerkin system $\Kr$ from \eqref{eq:fully_discrete_SG_algebraic_sep_KuF} (blue), the bound on the sparsity from \eqref{eq:Mpr_general_bound} (green), and the bound on the sparsity from \eqref{eq:Mpr_general_r_dep_bound} (red).}
\label{fig:bound_on_sparsity}
\end{center}
\end{figure}

We note that the bound of $\Mpr$ from \eqref{eq:Mpr_general_bound} is an overestimate due to the particular form of \eqref{eq:nnz_Gk_general}, which is different for each $\rb\in\Lambdar$. 
As a consequence, we see that $\text{nnz}(\mathbf{G}_\rb) = \mathcal{O}(\min\{2^{|\rb|}, M_{\lceil |\rb|/2 \rceil}\} \Mpsubhalfrb)$ for $\rb\in\Lambdar$.
For large $N$ and small $r$, $2^{r}$ is smaller than $M_{\lceil r/2\rceil}$,
however, as $r\to\infty$ the bound $M_{\lceil r/2\rceil}$ is sharper. For the $\eps$-complexity analysis in the next section, we note that 
\begin{align*}
\mathcal{M}(p,r) = \sum_{\rb \in \Lambdar} \textnormal{nnz}(\mathbf{G}_{\rb}) 
    & \leq \sum_{\rb\in\Lambdar} 2 \min \left\{ 2^{|\rb|}, {N + \lceil |\rb| /2 \rceil \choose N} \right\} {N+p - \lceil|\rb|/2\rceil\choose N} \\ 
    & = \sum_{j = 0}^r 2 \min \left\{ 2^{j}, {N + \lceil j /2 \rceil \choose N} \right\} {N-1 + j \choose N-1} {N+p - \lceil j/2\rceil\choose N} \\ 
    & \leq 2 \min \left\{ 2^{r}, {N + \lceil r /2 \rceil \choose N} \right\} {N+p \choose N} \sum_{j = 0}^r {N - 1 + j \choose N-1} \\  
    & = 2 \min \left\{ 2^{r}, {N + \lceil r /2 \rceil \choose N} \right\} {N+p \choose N} {N + r \choose N} \numberthis \label{eq:Mpr_general_r_dep_bound} 
\end{align*} 
Figure \ref{fig:bound_on_sparsity} plots how sharply $\Mpr$ is bounded by \eqref{eq:Mpr_general_bound} and \eqref{eq:Mpr_general_r_dep_bound}.
We are also able to show that Theorem \ref{thm:nnz_Gk_general} yields a sharp result in the case $|\rb|=1$. 

\begin{cor}
\label{cor:nnz_Gk_linear}
Under the same conditions in Theorem \ref{thm:nnz_Gk_general}, when $\rb\in\Lambdar$ is such that $|\rb|=1$, we have
\begin{align}
\label{eq:nnz_Gk_linear}
\text{\em nnz}(\mathbf{G}_\rb) = \sum_{\ell = \lceil |\rb|/2 \rceil}^{|\rb|} c(\rb,\ell) { N + p - \ell \choose p - \ell} = 2 { N + p - 1 \choose p - 1}.
\end{align}
\end{cor}

Corollary \ref{cor:nnz_Gk_linear} is the result of \cite[Lemma 28]{Ernst2010}, and follows from 
the application of 
the exact formula for $\text{nnz}(\mathbf{G}_\rb)$ from \eqref{eq:nnz_Gk_general}. Here $|\rb| = \lceil |\rb|/2 \rceil = 1$ is odd and $\bm{S}(\rb,1)$ has only one element $\bm{S}(\rb,1) = \{ \sbb \in \N_0^N : |\sbb| = 1, \sbb \leq \rb \} = \{ \rb \}$. Hence $c(\rb,1) = 2\#\bm{S}(\rb,1) = 2$, and \eqref{eq:nnz_Gk_linear} is shown. We are also able to show that the formula for $\text{nnz}(\mathbf{G}_\rb)$ from \eqref{eq:nnz_Gk_general} yields a result that is sharp in the case $N=1$.

\begin{cor}
\label{cor:sharp_bounds}
Under the same conditions in Theorem \ref{thm:nnz_Gk_general}, when $N=1$ and $\rb = r \in \N_0$, we have 
\begin{list}{(\alph{qcounter})~}{\usecounter{qcounter}}
\item in case $r=2k$, $k\in\N_0$,
\begin{align}
\textnormal{nnz}(\mathbf{G}_r) = \begin{cases}
(p-r+1)(r+1)+k^2, & 0\leq r\leq p, \\
(p-k+1)^2,        & p+1 \leq r \leq 2p, \\
0,                & r > 2p.
\end{cases}
\end{align}
\item in case $r=2k+1$, $k\in\N_0$,
\begin{align}
\textnormal{nnz}(\mathbf{G}_r) = \begin{cases}
(p-r+1)(r+1)+k^2+k, & 0\leq r\leq p, \\
(p-k+1)(p-k),       & p+1 \leq r \leq 2p, \\
0,                  & r > 2p.
\end{cases}
\end{align}
\end{list}
\end{cor}

Corollary \ref{cor:sharp_bounds} is the result of \cite[Lemma 25]{Ernst2010}, and its proof using Theorem 4.1 is included in the Appendix. 
In the remarks that follow, we make a distinction between the cases that $a(x,\y)$ is a polynomial of fixed degree $\rbar< \infty$, e.g., the coefficients from Examples \ref{exmp:affine_coefficient} and \ref{exmp:polynomial_coefficient}, and that $a(x,\y)$ is a transcendental function of the random variables, e.g., the coefficient from Example \ref{exmp:transcendental_coefficient}.

\begin{rem}
\label{rem:SG_affine_and_poly_nnz_complexity}
{\rm (Complexity of matrix-vector products for polynomial coefficients, see e.g., Examples \ref{exmp:affine_coefficient} and \ref{exmp:polynomial_coefficient})}
From Corollary \ref{cor:nnz_Gk_linear} and the work estimate
\eqref{eq:SG_cost_solve_precon} when using \eqref{eq:mean_based_preconditioner} as a preconditioner, 
    we see that for coefficients that are affine functions of the random variables, e.g., Example \ref{exmp:affine_coefficient},
the complexity of a single PCG iteration is of the order $\mathcal{O}(\Jh (2 \Mp +  2 N M_{p-1})) = \mathcal{O}(\Jh \Mp)$,
where $\Mp = \# \LambdapTD = {N+p\choose N}$.
On the other hand, when the coefficient $a(x,\y)$ is a polynomial function of the random variables, e.g., Example \ref{exmp:polynomial_coefficient}, having fixed order $\overline{r} \in \N$, $\rbar<\infty$, we use Theorem \ref{thm:nnz_Gk_general} to obtain a different estimate.
Since $\{\Psi_\rb \}_{\rb \in \Lambdarbar}$ is a basis for the space $\mathcal{P}_{\Lambdarbar}(\Gamma)$, there exits coefficients $\{a_\rb(x)\}_{\rb\in\Lambdarbar}$ such that $a(x,\y) = a^{\rbar}(x,\y) = \sum_{\rb\in\Lambdarbar} a_{\rb}(x) \Psi_{\rb} (\y)$. With this representation, it is clear to see that substituting $a(x,\y)$ into \eqref{eq:operator_expansion} yields $\Krbar$ from \eqref{eq:fully_discrete_SG_algebraic_sep_KuF}, and $\Krbar = \K$ from \eqref{eq:fully_discrete_SG_algebraic_KuF}. 
However, it is not clear how many of the coefficients $a_\rb(x)$ are identically zero. 
In this case, 
we can
provide an upper bound on the block-sparsity of $\Krbar$ under the assumption that
$a_\rb(x) \not\equiv 0$ $\forall \rb\in\Lambda_{\overline{r}}$.
Using the bound of \eqref{eq:Mpr_general_r_dep_bound}, the complexity of a single matrix-vector product of $\mathbf{K}_{\rbar}$ is of the order $\mathcal{O}(\Jh \Mp M_{\rbar} \min\{ 2^{\rbar}, M_{\lceil \rbar/2 \rceil}\})$. 
Thus, when $\rbar$ is fixed, $M_{\rbar} \min\{ 2^{\rbar}, M_{\lceil \rbar/2 \rceil}\}$ is a constant, and this estimate has the same asymptotic complexity as $\mathcal{O}(\Jh \Mp)$.
\end{rem}

\begin{rem}
{\rm (Complexity of matrix-vector products in the trascendental case, see e.g., Example \ref{exmp:transcendental_coefficient})}
\label{rem:SG_nonpoly_nnz_complexity}
We recall the discussion of \cite[Section 3.4]{Ullmann2010}. There, the complexity of matrix-vector products with the SG system was estimated when a full orthogonal expansion is substituted into the SG discretization. This case corresponds to fixing the expansion order $r=2p$ following Remark \ref{rem:projection_well_posed}. Assuming that $\textnormal{nnz}(\mathbf{G}_{\rb}) = \mathcal{O}(\Mp)$ or $\mathcal{O}(\Mp^2)$,  it was estimated in \cite{Ullmann2012} that the cost of matrix-vector products involving $\Kr$ is be between $\mathcal{O}(\Jh \Mp^2)$ and $\mathcal{O}(\Jh \Mp^3)$. 
However, the use of Theorem \ref{thm:nnz_Gk_general} allows us to consider the complexity in the case of truncating the expansion, where a sharper estimate can be obtained.
Let $T_r := \prod_{k=\lceil r/2 \rceil+1}^{r} \frac{N+k}{k} \ll M_{\lceil r/2 \rceil} = {N + \lceil r/2 \rceil \choose \lceil r/2 \rceil}$, which is bounded independent of $r$, i.e.,
\begin{align*}
T_r \leq \left(\frac{N + \lceil r/2 \rceil+1}{\lceil r/2 \rceil+1} \right)^{\lceil r/2 \rceil+1} \to e^N \;\;\; \text{ as } \;\;\; r \to \infty,
\end{align*}
so that $\Mr = T_r M_{\lceil r/2 \rceil} \leq e^N M_{\lceil r/2\rceil}$.
From \eqref{eq:Mpr_general_r_dep_bound}, we see that 
$\Mpr$ is of the order $\mathcal{O}(\Mp \Mr M_{\lceil r/2 \rceil})$ as $p,r\to\infty$, since $\min\{2^r, M_{\lceil r/2 \rceil} \} \to M_{\lceil r/2 \rceil}$ as $r\to\infty$.
When $r = 2p$, this implies the complexity of matrix-vector multiplications involving $\Kr$ is of the order $\mathcal{O}(\Jh \Mp \Mr M_{\lceil r/2 \rceil}) = \mathcal{O}(\Jh \Mp T_r M_{\lceil r/2 \rceil}^2)=
\mathcal{O}(\Jh M_{p}^3)$. On the other hand, when $r = p$, we see that the complexity of matrix-vector products with $\Kr$ is order $\mathcal{O}(\Jh \Mp \Mr M_{\lceil r/2 \rceil}) = \mathcal{O}(\Jh T_r^2 M_{\lceil p/2 \rceil}^3 ) = \mathcal{O}(\Jh M_{\lceil p/2 \rceil}^3)$.
\end{rem}

\subsection{$\eps$-complexity analysis of the SGFEM}
\label{subsec:eps_complexity_analysis}

An estimate of the total complexity to obtain a fully discrete approximation of tolerance $\eps>0$ with the SGFEM and PCG solver can be shown in four steps:
\begin{enumerate}
\item Estimate the maximum mesh size $h_{\max}$ and minimum polynomial order $p_{\min}$ necessary in the finite element and SG discretizations, respectively, 
\item If projection of the coefficient is necessary, estimate the minimum projection order $r_{\min}$, otherwise set $r_{\min} = \rbar$ where $\rbar<\infty$ is the order of the coefficient,
\item Estimate the minimum number of iterations $k_{\min}$ needed by the PCG solver,
\item Substitute $h_{\max}$, $p_{\min}$, $r_{\min}$, and $k_{\min}$ into the cost \eqref{eq:SG_cost_solve_precon} for $h$, $p$, $r$, and $\NiterpreSG$, respectively.
\end{enumerate}
We proceed to estimate these parameters as follows. Denote by $u^r$ the corresponding solution of \eqref{eq:model_problem} when $a^r(x,\y)$ is substituted in place of $a(x,\y)$, and let 
$\uSGCGr$ be the approximation to $\uSGr$ found by PCG. 
Then the total error for the SG approximation satisfies the following bound:
\begin{align*}
\left\|u - \uSGCGr \right\|_{\Htworho} & \leq \underbrace{\left\| u - \ur \right\|_{\Htworho}}_{\text{SG(I)}} + \underbrace{\left\| \ur - \urh \right\|_{\Htworho}}_{\text{SG(II)}} \numberthis \label{eq:SG_full_error_with_projection} 
                          + \underbrace{\left\| \urh - \uSGr \right\|_{\Htworho}}_{\text{SG(III)}} + \underbrace{\left\| \uSGr - \uSGCGr \right\|_{\Htworho}}_{\text{SG(IV)}} .
\end{align*}
In this setting SG(I) is the approximation error using a truncated expansion of $a(x,\y)$, SG(II) is the discretization error induced by the finite element method, SG(III) is the SG error coming from the orthogonal expansion \eqref{eq:fully_discrete_SG_soln}, and SG(IV) is the solver error resulting from the PCG method. 
We note that when the projection of the coefficient is exact, as discussed in Remark \ref{rem:projection_well_posed}, the approximation error SG(I) is no longer present and $\uSGr\equiv \uSG$.

We start with bounding SG(III). Without loss of generality, it is reasonable to assume that since $u^r$ has a holomorphic dependence on $\z\in\C^N$ in an open neighborhood of the polyellips $\mathcal{E}_{\bs{\gamma}}$ from Theorem \ref{thm:regularity}, then $\urh$ does as well. 
Then, the following result, whose proof is found in \cite{Todor2007}, and immediately follows from classical spectral convergence results \cite{Davis1975,Szego1975}, describes the convergence rate of the fully discrete solutions obtained by the SG method using a total degree approximation in $\mathcal{P}_{\LambdapTD}(\Gamma)$:

\begin{prop}[\rm Convergence rate for the SG method]
\label{prop:SG_TD_convergence_rate}
If Theorem \ref{thm:regularity} holds for the solution $\urh$ to \eqref{eq:semi_discrete_SFEM} with coefficient $a^r(x,\y)$, and $\uSGr$ is the solution to \eqref{eq:fully_discrete_SG_problem} with $\Lambdap$ the order $p$ total degree index set, then 
\begin{align*}
\| \urh - \uSGr \|_{\Hinfrho} \leq C_1 \exp(-C_2 p) \qquad \forall p\in \N, 
\end{align*}
for some constants $C_1, C_2 > 0$ independent of $p$.
\end{prop}

To investigate the error in SG(I), we note that since $a(x,\y)$ satisfies assumption $\Atwo$, the projection error in $\mathcal{P}_{\LambdarTD}(\Gamma)$ can be similarly estimated as 
\begin{align}
\label{eq:pointwise_trunc_rate}
\|a - a^r \|_{L^2_\varrho(\Gamma; L^\infty(D))} \leq C_3 \exp(-C_4 r) \qquad \forall r\in \N, 
\end{align}
for some constants $C_3, C_4>0$ independent of $r$. Hence, $\forall r\in\N$,
\begin{align}
\label{eq:u_a_truncation_error}
\| u - \ur \|_{\Htworho} \leq \frac{\|f\|_{H^{-1}}}{a_{\min}^2} \|a - a^r\|_{L^2_\varrho(\Gamma;L^\infty(D))} \leq \frac{\|f\|_{H^{-1}}}{a_{\min}^2} C_3 \exp(-C_4 r)  
\end{align}
providing a bound for SG(I).  For a bound of SG(II), we present the following convergence result reguarding solutions to the parameterized finite element problem, whose proof can be found in a number of standard texts on the theory of finite element methods, e.g., \cite{Atkinson2005, Johnson2012}:
\begin{lem}
\label{lem:fem_rate}
Let $\mathcal{T}_h$ be a uniform finite element mesh over $D$ with $\Jh = \mathcal{O}(h^{-d})$ degrees of freedom and $h>0$. For the elliptic PDE \eqref{eq:model_problem} and $\y\in\Gamma$, when $\ur(\y)\in H_0^1(D)\cap H^{s+1}(D)$, the error from the finite element approximation is bounded by
\begin{align*}
\|\ur(\y) - \urh(\y) \|_{H_0^1(D)} \leq \CFEM h^s, 
\end{align*}
where the constant $\CFEM>0$ is independent of $h$ and $\y$. 
\end{lem}

For the treatment of SG(IV), 
we begin by defining 
$\mathcal{B}^r(\y)$ to be the corresponding bilinear operator in \eqref{eq:deterministic_weak_problem_BF} with $a(x,\y)$ replaced with $a^r(x,\y)$. Since both $\mathcal{B}(\y)$ and $\mathcal{B}^r(\y)$ are symmetric, uniformly coercive and continuous bilinear operators on $H_0^1(D)$, there exist $\alpha,\beta > 0$ independent of $\y$ such that for every $u,v \in H^1_0(D)$
\begin{align*}
\left| \mathcal{B}^r[u,v](\y) \right| & = \left| \int_D a^r(x,\y) \nabla u \cdot \nabla v dx \right| \leq \alpha \| u \|_{H_0^1(D)} \| v \|_{H_0^1(D)}, \;\; \text{ and } \\ 
    \beta \|u \|^2_{H_0^1(D)} & \leq \int_D a^r(x,\y) | \nabla u |^2 dx = \| u \|_{\mathcal{B}^r(\y)}^2,
\end{align*}
and similarly for $\mathcal{B}_r(\y)$ with the same $\alpha,\beta$, e.g., taking $\alpha$ to be the maximum and $\beta$ to be the minimum in each case.
Recall $\mathbf{u}_{h,\p}^r = [u^r_{1,\p},\ldots,u^r_{\Jh,\p}]^{\text{T}}$, the vector of nodal values of the finite element solution corresponding to the $\p$-th stochastic mode of $\uSGr$, and $\mathbf{u}_{h,p}^r = [\mathbf{u}_{h,\p}^{r}]_{\p\in\Lambdap}^{\text{T}}$. Then we have the following estimates expressing
\begin{align}
\label{eq:SG_continuity_condition}
& \text{Continuity:} \qquad \left\| \mathbf{u}_{h,p}^r \right\|_{\Kr} = \left\| \uSGr \right\|_{\bbE[\mathcal{B}^r(\y)]} \leq \sqrt{\alpha} \left\| \uSGr \right\|_{\Htworho}, \;\;\; \text{and} \\
\label{eq:SG_ellipticity_condition}
& \text{Ellipticity:} \qquad \sqrt{\beta} \left\| \uSGr \right\|_{\Htworho} \leq \left\| \uSGr \right\|_{\bbE[\mathcal{B}^r(\y)]} = \left\| \mathbf{u}_{h,p}^r \right\|_{\Kr},
\end{align}
where $\| \mathbf{u} \|_{\Kr}^2 = (\mathbf{u})^{\text T} \Kr \mathbf{u}$ is the $\Kr$ matrix norm, and $\| u \|_{\bbE[\mathcal{B}^r(\y)]}$ is the expectation of the energy norm \eqref{eq:Bnorm}.  
Given Proposition \ref{prop:SG_TD_convergence_rate}, Lemma \ref{lem:fem_rate}, and the estimates from \eqref{eq:u_a_truncation_error}, \eqref{eq:SG_continuity_condition}, and \eqref{eq:SG_ellipticity_condition}, we can now provide the minimal projection orders $p,r\in\N$ for the SG approximation \eqref{eq:fully_discrete_SG_soln} and the coefficient \eqref{eq:a_r_truncation}, respectively, the maximum mesh size $h$ for finite element method, and the minimum number of PCG iterations $k$ necessary to ensure that the error in the SGFEM solution $\uSGCGr$ is less than the tolerance $\eps>0$.

\begin{lem}
\label{lem:SG_optimal_parameters}
Let $u\in L^2_\varrho(\Gamma;H_0^1(D)\cap H^{s+1}(D))$ be the solution to \eqref{eq:model_problem}, $\uSGr$ be the solution to \eqref{eq:fully_discrete_SG_problem} with the coefficient $a^r(x,\y)$, and $\uSGCGr$ be the approximation of $\uSGr$ found by PCG with a zero initial guess. Then, for $\eps>0$, to ensure that
$\| u - \uSGCGr \|_{\Htworho} \leq \eps$ we must choose $h \leq h_{\max}$, $r\geq r_{\min}$, $p \geq p_{\min}$, and $k \geq k_{\min}$, where:
\begin{alignat*}{2}
h_{ \max } & = \left( \frac{\varepsilon}{4 \CFEM} \right)^{\frac{1}{s}},            & \qquad   r_{ \min } & = \log\left[\left(\frac{4 C_5}{\eps}\right)^{\frac{1}{C_4}}\right], \\ 
p_{ \min } & = \log\left[\left(\frac{4 C_{1}}{\eps}\right)^{\frac{1}{C_2}}\right],  & \qquad   k_{ \min } & = \frac{ \log \left( \frac{4 C_6}{\varepsilon} \right)}{\log \left(\frac{\sqrt{\kappartilde} + 1}{\sqrt{\kappartilde} - 1} \right) } ,
\end{alignat*}
with $\CFEM>0$ the constant from Lemma \ref{lem:fem_rate}, $C_1,C_2,C_3,C_4>0$ the constants from Proposition \ref{prop:SG_TD_convergence_rate} and
\eqref{eq:pointwise_trunc_rate}, and, for $\alpha,\beta>0$ from \eqref{eq:SG_continuity_condition} and \eqref{eq:SG_ellipticity_condition}
\begin{align*}
C_5 = C_3 \frac{\|f\|_{H^{-1}}}{a_{\min}^2}, \qquad C_6 = 2 \sqrt{\frac{\alpha}{\beta}} \left\| \uSGr \right\|_{\Htworho}, 
\end{align*} 
with $\kappartilde$ is the condition number of $\mathbf{\tilde P}^{-1} \Kr$ with $\mathbf{P}$ the mean-based preconditioner from \eqref{eq:mean_based_preconditioner}. 
\end{lem}

\begin{pf}
Without loss of generality, we seek to bound the quantities SG(I)-SG(IV) from \eqref{eq:SG_full_error_with_projection} each by $\eps/4$.
For the error SG(I) we recall estimate \eqref{eq:u_a_truncation_error}
and solve for $r$. 
From Lemma \ref{lem:fem_rate}, when $u\in L^2_\varrho(\Gamma; H_0^1(D)\cap H^{s+1}(D))$ we have that 
$\| \ur - \urh \|_{\Htworho} \leq \CFEM h^{s}$  $\forall h>0,$
and from Proposition \ref{prop:SG_TD_convergence_rate} we have that 
$\| \urh - \uSGr \|_{\Hinfrho} 
\leq C_{1} \exp(-C_2 p)$ $\forall p \in \N$,
so that solving for $h$ and $p$ gives the desired maximum mesh size $h_{\max}$ and minimum polynomial order $\pmin$ to bound SG(II) and SG(III) by $\varepsilon/4$.
Let $\mathbf{u}_{h,p}^r$ and $\mathbf{u}_{h,p}^{r,(k)}$ be the coefficients of the exact SG solution $\uSGr$ and the approximate SG solution $\uSGCGr$ after $k$ PCG iterations, respectively. Then 
from \eqref{eq:CG_error_SG} and \eqref{eq:SG_ellipticity_condition} we see that
\begin{align*}
\left\| \uSGr - \uSGCGr \right\|_{\Htworho} & \leq \frac{1}{\sqrt{\beta}} \| \mathbf{u}_{h,p}^{r} - \mathbf{u}_{h,p}^{r,(k)} \|_{\Kr} 
    \leq \frac{2}{\sqrt{\beta}} \left(\frac{\sqrt{\kappartilde} - 1}{\sqrt{\kappartilde} + 1} \right)^k \| \mathbf{u}_{h,p}^{r} - \mathbf{u}_{h,p}^{r,(0)} \|_{\Kr},
\end{align*}
where $\mathbf{u}_{h,p}^{r,(0)}$ is the initial guess used in CG and $\kappartilde = \text{cond}(\mathbf{\tilde P}^{-1}\Kr)$ with mean based preconditioner $\mathbf{P}$ from \eqref{eq:mean_based_preconditioner}. 
If we use the zero vector as the initial iteration in PCG, we have from \eqref{eq:SG_continuity_condition}
\begin{align}
\label{eq:CG_err_l2_estimate}
\left\| \uSGr - \uSGCGr \right\|_{\Htworho} & \leq \frac{2}{\sqrt{\beta}} \left(\frac{\sqrt{\kappartilde} - 1}{\sqrt{\kappartilde} + 1} \right)^k \| \mathbf{u}_{h,p}^{r} \|_{\Kr}
                                                                           \leq 2 \sqrt{\frac{\alpha}{\beta}}  \left(\frac{\sqrt{\kappartilde} - 1}{\sqrt{\kappartilde} + 1} \right)^k \left\| \uSGr \right\|_{\Htworho}.
\end{align}
Solving for $k$ gives the minimum number of iterations $k_{\min}$ required to ensure SG(IV) is bounded by $\varepsilon/4$.
\qed
\end{pf}

Given the necessary parameters from Lemma \ref{lem:SG_optimal_parameters} to achieve $\|u - \uSGCGr \|_{\Htworho} \leq \eps$, and the estimates on the computational complexity of one iteration in the PCG method from \S \ref{subsec:SG_mat_vec_complexity}, we provide a bound on the minimal number of FLOPs required by the SGFEM when approximating \eqref{eq:model_problem}. 
We split these results into the cases that the stochastic coefficient $a(x,\y)$ from \eqref{eq:model_problem} is:
\begin{list}{(\roman{qcounter})~}{\usecounter{qcounter}}
\item an affine function of the random parameters, e.g., $a(x,\y)\in\mathcal{P}_{\Lambdarbar}(\Gamma)$ with $\rbar=1$, as in Example \ref{exmp:affine_coefficient}, 
\item a non-affine polynomial of the random parameters, e.g., $a(x,\y)\in\mathcal{P}_{\Lambdarbar}(\Gamma)$ for some $1<\rbar<\infty$, as in Example \ref{exmp:polynomial_coefficient}, 
\item a non-affine, transcendental function of the random parameters, e.g., $a(x,\y)\not\in\mathcal{P}_{\Lambdar}(\Gamma)$ for any $r\in\N$, as in Example \ref{exmp:transcendental_coefficient}, so that $r$ must be chosen to satisfy $r\geq r_{\min}$ from Lemma \ref{lem:SG_optimal_parameters}.
\end{list}
The results are summarized in Theorems \ref{thm:SG_opt_work_poly} and \ref{thm:SG_opt_work_nonpoly} next.
\begin{thm}
\label{thm:SG_opt_work_poly}
Let $u\in L^2_\varrho(\Gamma; H_0^1(D)\cap H^{s+1}(D))$ be the solution to \eqref{eq:model_problem}, and
$\overline{r}$ be the smallest natural number such that $a(x,\y)\in \mathcal{P}_{\Lambda_{\overline{r}}}(\Gamma)$. 
When $\overline{r}=1$, the minimum work \eqref{eq:SG_cost_solve_precon} of solving \eqref{eq:fully_discrete_SG_algebraic_sep_KuF} with PCG to a tolerance $\eps>0$ 
can be bounded by
\begin{align}
\label{eq:SG_lin_work_estimate}
\WSGaffine \leq C_7 \left(\frac{3 \CFEM}{\eps} \right)^{\frac{d}{s}}  
2 e^N (1+N) \left( 1+ \log\left[ \left( \frac{3 C_1}{\eps} \right)^{\frac{1}{C_2 N}} \right] \right)^N
\left( \frac{ \log \left( \frac{3 C_6}{\varepsilon} \right)}{\log \left(\frac{\sqrt{\kappabtilde} + 1}{\sqrt{\kappabtilde} - 1} \right) } \right),
\end{align}
and when $\overline{r}>1$, the minimum work \eqref{eq:SG_cost_solve_precon} of solving \eqref{eq:fully_discrete_SG_algebraic_sep_KuF} with PCG to a tolerance $\eps>0$ 
can be bounded by
\begin{equation}\label{eq:SG_poly_work_estimate}
\begin{aligned}
& \WSGpoly \leq C_7 \left( \frac{ 3 \CFEM}{\eps} \right)^{\frac{d}{s}} 
\left(\frac{ \log \left( \frac{3 C_6}{\varepsilon} \right)}{\log \left(\frac{\sqrt{\kappabtilde} + 1}{\sqrt{\kappabtilde} - 1} \right) }\right)
    e^N \left[ \left( 1 + \log \left[ \left( \frac{3 C_1}{\eps} \right)^{\frac{1}{C_2 N}} \right] \right)^N \right. \\
      &\hspace{1cm} \left. + 2 \sum_{j=0}^{\rbar} {N - 1 + j \choose N - 1} \min \left\{ 2^j, {N + \lceil j/2 \rceil \choose N} \right\} \left( 1 - \frac{\lceil j/2 \rceil}{N} + \log \left[ \left( \frac{3 C_1}{\eps} \right)^{\frac{1}{C_2 N}} \right] \right)^N \right],
\end{aligned}
\end{equation}
with $\CFEM,C_1,C_2,C_6$ as in Lemma \ref{lem:SG_optimal_parameters}, $C_7>0$ independent of $\eps$, and $\kappabtilde$ the condition number of the preconditioned system $\mathbf{\tilde P}^{-1} \mathbf{K}_{\rbar} = \mathbf{\tilde P}^{-1} \mathbf{K}$, using the mean-based preconditioner from \eqref{eq:mean_based_preconditioner}. 
\end{thm}

\begin{pf}
When $a(x,\y)\in\mathcal{P}_{\Lambdarbar}(\Gamma)$ we do not need to consider SG(I) from \eqref{eq:SG_full_error_with_projection}, and bound SG(II), SG(III), and SG(IV) by $\varepsilon/3$. 
Hence, to minimize the error of the SG discretization, we choose $p\geq \pmin = \log[(3 C_{1}/\varepsilon)^{1/C_2}]$ which differs from the $\pmin$ stated in Lemma \ref{lem:SG_optimal_parameters}. 
For a uniform triangulation $\mathcal{T}_h$, $\Jh = \mathcal{O}(h^{-d})$ so that 
\begin{align}
\label{eq:Jhmax}
J_{h_{\max}} = C_7 \left[ \left( \frac{\eps}{3 \CFEM} \right)^{\frac{1}{s}} \right]^{-d} = C_7 \left( \frac{3 \CFEM}{\eps} \right)^{\frac{d}{s}} 
\end{align}
for some constant $C_7>0$ depending on the connectivity of the finite element mesh, but independent of $\varepsilon$.
In the case that $\rbar=1$, we substitute $\pmin$ into 
\eqref{eq:nnz_Gk_general} for the matrices $\mathbf{G}_{\rb}$ having $0\leq |\rb| \leq 1$, and 
apply Stirling's approximation to obtain 
\begin{align*}
M_{\pmin} + \Mpminone 
        = 2{N + \pmin \choose N} + 2N {N + \pmin - 1 \choose N} 
        \leq 2 e^N (1+N) \left( 1 + \log\left[\left(\frac{3 C_{1}}{\varepsilon}\right)^{\frac{1}{C_2N}}\right] \right)^N ,
\end{align*}
Similarly, when $\rbar>1$ we use the bound from \eqref{eq:Mpr_general_bound} and Stirling's approximation to obtain
\begin{align*}
& M_{\pmin} + \Mpminrbar 
     \leq e^N \left[ \left(1 + \log\left[\left( \frac{3 C_1}{\varepsilon} \right)^{\frac{1}{C_2 N}} \right] \right)^{N} \right. \\ 
     & \qquad\qquad\qquad \left. + 2 \sum_{j=0}^{\rbar} \min \left\{ 2^j, {N + \lceil j/2 \rceil \choose N} \right\} {N - 1 + j \choose N - 1} \left( 1 - \frac{\lceil j/2 \rceil}{N} + \log\left[ \left( \frac{3 C_1}{\varepsilon} \right)^{\frac{1}{C_2 N}} \right] \right)^N \right] .
\end{align*}
Substituting $J_{h_{\max}}$ for $J_h$, $k_{\min}$ for $\NiterSG$ from Lemma \ref{lem:SG_optimal_parameters}, and the bounds for $M_{\pmin} + N_K^{(\pmin,\rbar)}$ into the work estimate \eqref{eq:SG_cost_solve_precon}, in the cases $\rbar=1$ and $\rbar>1$ above, we obtain the desired results.
\qed
\end{pf}

\begin{thm}
\label{thm:SG_opt_work_nonpoly}
Let $u\in L^2_\varrho(\Gamma; H_0^1(D)\cap H^{s+1}(D))$ be the solution to \eqref{eq:model_problem}, 
and suppose that $a(x,\y) \not \in \mathcal{P}_{\Lambdar}(\Gamma)$ for any $r\in\N$. In this case $r$ must be chosen to satisfy $r\geq \rmin$ from Lemma \ref{lem:SG_optimal_parameters}. 
Then the minimum work  \eqref{eq:SG_cost_solve_precon} of solving 
\eqref{eq:fully_discrete_SG_algebraic_sep_KuF} with PCG to a tolerance $\eps>0$ can be bounded by 
\begin{align*}
& \WSGnonpoly \leq C_7 \left( \frac{ 4 \CFEM}{\eps} \right)^{\frac{d}{s}} 
e^N \left( 1 + \log \left[ \left( \frac{4 C_1}{\eps} \right)^{\frac{1}{C_2 N}} \right] \right)^N 
\numberthis \label{eq:SG_nonpoly_work_estimate} \\
        & \hspace{1.3cm} \left[1 + 2 e^{2N} \left( 1 + \log \left[ \left( \frac{4 C_5}{\eps} \right)^{\frac{1}{C_4 N}} \right] \right)^N \left( 1 + \frac{1}{N} + \log \left[ \left( \frac{4 C_5}{\eps} \right)^{\frac{1}{2 C_4 N}} \right] \right)^N \right]
\left( \frac{ \log \left( \frac{4 C_6}{\varepsilon} \right)}{\log \left(\frac{\sqrt{\kappartilde} + 1}{\sqrt{\kappartilde} - 1} \right) } \right), 
\end{align*}
with $\CFEM,C_1,C_2,C_4,C_5,C_6$ as in Lemma \ref{lem:SG_optimal_parameters}, $C_7>0$ independent of $\eps$, and $\kappartilde$ the condition number of the preconditioned system $\mathbf{\tilde P}^{-1} \Kr$, using the mean-based preconditioner from \eqref{eq:mean_based_preconditioner}. 
\end{thm}

\begin{pf}
In this setting $r$ must be chosen to satisfy $r\geq \rmin = \log [ ( 4 C_5/\varepsilon )^{1/C_4} ]$ from Lemma \ref{lem:SG_optimal_parameters} and, therefore, we must bound the sum from \eqref{eq:Mpr_general_bound} which now depends on $\rmin$, and hence on $\eps$. Thus, we use the bound \eqref{eq:Mpr_general_r_dep_bound} for $\Mpr$, noting that as $\eps \to 0$, $r_{\min}\to \infty$ so that 
\begin{align*}
\min\left\{ 2^{\rmin}, {N + \lceil \rmin / 2 \rceil \choose N } \right\}  = {N + \lceil \rmin / 2 \rceil \choose N } . 
\end{align*}
Substituting $\pmin$ and $\rmin$ from Lemma \ref{lem:SG_optimal_parameters} into \eqref{eq:Mpr_general_r_dep_bound} and applying Stirling's approximation, we obtain
\begin{align*}
& M_{p_{\min}} + \Mpminrmin 
\leq e^N \left( 1 + \log \left[ \left( \frac{4 C_1}{\eps} \right)^{\frac{1}{C_2 N}} \right] \right)^N 
\\
        & \hspace{2.3cm} \left[1 + 2 e^{2N} \left( 1 + \log \left[ \left( \frac{4 C_5}{\eps} \right)^{\frac{1}{C_4 N}} \right] \right)^N \left( 1 + \frac{1}{N} + \log \left[ \left( \frac{4 C_5}{\eps} \right)^{\frac{1}{2 C_4 N}} \right] \right)^N \right].
\end{align*}
As in the proof of Theorem \ref{thm:SG_opt_work_poly}, we substitute $J_{h_{\max}}$ for $\Jh$ from \eqref{eq:Jhmax}, $k_{\min}$ for $\NiterSG$, and the bound for $M_{p_{\min}} + \Mpminrmin$ with $p_{\min}$ and $k_{\min}$ from Lemma \ref{lem:SG_optimal_parameters} into the cost \eqref{eq:SG_cost_solve_precon} to complete the proof.
\qed
\end{pf}

Given Theorems \ref{thm:SG_opt_work_poly} and \ref{thm:SG_opt_work_nonpoly} we see that the work of obtaining the fully discrete approximation using the SGFEM, with PCG as a solver, is asymptotically given by:
\begin{align}
\label{eq:SG_asymptotic_estimate}
\mathcal{O} \underbrace{
    \vphantom{\left( \frac{ \log\left( \frac{1}{\eps} \right) }{ \log \left( \frac{ \sqrt{\kappartilde} + 1 }{ \sqrt{\kappartilde} - 1} \right) } \right)}
    \left( \frac{1}{\eps} \right)^{\frac{d}{s}}}_{\text{(SG.1)}} 
\underbrace{
    \vphantom{\left( \frac{ \log\left( \frac{1}{\eps} \right) }{ \log \left( \frac{ \sqrt{\kappartilde} + 1 }{ \sqrt{\kappartilde} - 1} \right) } \right)}
\left[ \log \left( \frac{1}{\eps} \right) \right]^{g(N)} }_{\text{(SG.2)}}
    \underbrace{\left( \frac{ \log\left( \frac{1}{\eps} \right) }{ \log \left( \frac{ \sqrt{\kappartilde} + 1 }{ \sqrt{\kappartilde} - 1} \right) } \right)}_{\text{(SG.3)}} ,
\end{align}
where $g(N) = N$ and $\kappartilde = \kappatilde$ if $a(x,\y)$ is an affine or non-affine, polynomial function of the random parameters of fixed order $\rbar<\infty$, e.g., Examples \ref{exmp:affine_coefficient} and \ref{exmp:polynomial_coefficient}, and $g(N) = 3N$ when $a(x,\y)$ is a non-affine, transcendental function of the random parameters, e.g., Example \ref{exmp:transcendental_coefficient}, requiring a total degree orthogonal expansion of order $r\geq r_{\min}$ depending on $\eps$.
Here, 
(SG.1), (SG.2), and (SG.3) correspond to the work required by the finite element, SG, and PCG methods, respectively. In particular, (SG.2) corresponds to the estimates for the sparsity of the Galerkin system $\Kr$ from \eqref{eq:fully_discrete_SG_algebraic_sep_KuF}, 
and represents the number of coupled finite element systems that must be solved simultaneously by the PCG method.
However, due to the bound \eqref{eq:Mpr_general_bound}, the asymptotic complexity in the cases that $a(x,\y)$ is affine or polynomial in $\y$ are the same. 
This does not imply that there is no need to consider 
the work estimates in these cases separately.
Indeed, if $a(x,\y)$ is a polynomial having the representation $\sum_{\rb\in\Lambda_{\rbar}} a_{\rb}(x) \Psi_\rb(\y)$ where $a_{\rb}(x) \neq 0$ for all $\rb\in\Lambda_{\rbar}$, then the complexity of matrix-vector multiplications with $\mathbf{K}_r$ 
is of the order $\mathcal{O}(\Jh \Mp M_{\rbar} \min\{ 2^{\rbar}, M_{\lceil \rbar/2\rceil} \})$. 
Here, the constant $M_{\rbar} \min\{ 2^{\rbar}, M_{\lceil \rbar/2\rceil}\}$ grows rapidly with $\rbar$, suggesting that higher order polynomial functions of $\y$ require additional cost.

\subsection{Conditioning of the generalized SG system}
\label{subsec:SG_condition_number}

In this section, we discuss issues related to the conditioning of the linear system that results from the SGFEM discretization. 
We first recall \cite[Theorem 10]{Ernst2010}: 
the eigenvalues of the matrices $\{\mathbf{G}_\rb\}_{\rb\in\Lambdar}$ from \eqref{eq:G_rpq} lie in the interval $[\xi_\rb,\Xi_\rb]$, where
\begin{align}
\label{eq:G_k_eigenbounds}
\xi_\rb := \min\{\Psi_\rb(\y): \y\in\Gmlb\}, \qquad \Xi_\rb := \max\{\Psi_\rb(\y): \y\in\Gmlb\} , 
\end{align}
$\Gmlb$ is a tensor product grid of Gauss-Legendre quadrature points having $\m(\lb)=(m(l_1),\ldots,m(l_n))$ points in each direction, 
and $\lb$ is such that $m(l_n) := p + \lceil \frac{k_n + 1}{2} \rceil$, $n=1,\ldots,N$. 
Since $a^r(x,\y)$ satisfies $\Aone$, 
the analysis of \cite[Theorem 3.8]{Powell2009} shows that the eigenvalues for the preconditioned system $\mathbf{P}^{-1} \Kr$ lie in the interval  
$[1 - \underline{\tau}_r, 1 + \overline{\tau}_r]$ where 
\begin{align*}
\underline{\tau}_r = \frac{1}{a_{\min}} \sum_{\rb \in \Lambdar \atop |\rb|\neq 0} \xi_\rb \|a_\rb(x) \|_{L^\infty(D)}, \qquad \overline{\tau}_r = \frac{1}{a_{\min}} \sum_{\rb \in \Lambdar \atop |\rb|\neq 0} \Xi_\rb \|a_\rb(x) \|_{L^\infty(D)}. \numberthis \label{eq:Pinv_Kw_eigenbounds}
\end{align*}
As a result of \eqref{eq:Pinv_Kw_eigenbounds},
we see that in the case that the projection order $r$ of the coefficient $a^r(x,\y)$ depends on $\eps$, the condition number of the preconditioned system 
$\mathbf{P}^{-1} \Kr$ 
does as well through the number of terms in $\underline{\tau}_r$ and $\overline{\tau}_r$. This should come as no surprise since even in the case of the \KL expansion, the condition number of 
$\mathbf{P}^{-1} \mathbf{K}$ 
depends on the number of terms in the truncated \KL expansion which is chosen a-priori to minimize the error. 

\section{Comparison with the SCFEM}
\label{sec:SC_comparison}

In this section we compare our explicit cost bounds for the SGFEM with the complexity estimates for the SCFEM developed in \cite{Galindo2015}, when solving \eqref{eq:model_problem}.
The basic idea behind the SCFEM is to 
construct a fully discrete approximation in a subspace of $V_h(D)\otimes L^2_\varrho(\Gamma)$ by
collocating semi-discrete solutions $\uh$ from \eqref{eq:semi_discrete_SFEM} on a deterministic set of points 
to obtain 
solutions $\{\uh(\cdot,\y_k)\}_{k=1}^{M_L}\in V_h(D)$. 

\subsection{A generalized SCFEM using Lagrange interpolation}
\label{subsec:gen_SCFEM}

To construct the stochastic collocation (SC) approximation, we 
consider a class of multi-index sets  defined in terms of increasing functions $\bm{m}:\N_+^N \to \N_+^N$ and $g:\N_+^N \to \N_+$. 
By $\bm{m}$ we specify the multivariate function $\bm{m}(\lb) := (m_1(l_1),\cdots,m_N(l_N))$ where each $m_n:\N_+\to\N_+$ is an increasing function, possibly different for each $n = 1,\ldots,N$. Here the $m_n$ are referred to as growth functions, specifying how the number of points grows in the direction $n$.
Associated with $m_n$ we define the {\em left-inverse} $m_n^\dagger:\N_+ \to \N_+$ by $m_n^\dagger(q) = \min\{k \in \N_+ : m_n(k) \geq q\}$, and let $\bm{m}^\dagger(\q) = (m_1^\dagger(q_1),\ldots,m_N^\dagger(q_N))$. In this case, we note that $m_n^\dagger(m_n(k))=k$ and $m_n(m_n^\dagger(k))\geq k$ for each $k\in\N_+$ and $n=1,\ldots,N$. Given $\bm{m}$ and $g$ we can define the multi-index set 
\begin{align}
\label{eq:general_index_set}
\LambdaLmg  = \left\{ \q \in \N_+^N : g(\bm{m}^{\dagger}(\q + \bm{1})) \leq L \right\},
\end{align}
to be used in constructing polynomial approximations.
In particular, setting $m_n(j) = j$ for all $j\in\N_+$ and $n=1,\ldots,N$, and defining 
\begin{equation}
\label{eq:g_defs}
g_\text{TP} ( \p ) = \max_{1\leq n \leq N} p_n, \qquad
g_\text{TD} ( \p ) = \sum_{n=1}^N (p_n-1), \qquad
g_\text{SM} ( \p ) = \sum_{n=1}^N f(p_n),
\end{equation}
where $f(p)$ is given in \eqref{eq:index_set_examples},
and using the definition of $\LambdaLmg$ from \eqref{eq:general_index_set}, we obtain the TP, TD, and SM index sets $\LambdaLTP$, $\LambdaLTD$, and $\LambdaLSM$, respectively, given in \eqref{eq:index_set_examples}.

We introduce a sequence of one-dimensional Lagrange interpolation operators $\mathcal{U}^{m_n(l_n)}:C^0(\Gamma_n)\to \mathcal{P}_{m_n(l_n)-1}(\Gamma_n)$. 
Then for $v\in C^0(\Gamma)$ the generalized multi-dimensional approximation operator $\ILmg: C^0(\Gamma) \to \mathcal{P}_{\LambdaLmg}(\Gamma)$ is given by 
\begin{align*}
\ILmg[v](\y) 
                             = \sum_{g(\lb) \leq L} \sum_{\ib\in\{0,1\}^N} (-1)^{|\ib|} \left( \bigotimes_{n=1}^N \mathcal{U}_n^{m_n(l_n-i_n)} \right) [v](\y).
                            \numberthis \label{eq:SC_approximation_operator} 
\end{align*}
Construction of the approximation $\ILmg[v](\y)$ requires the independent evaluation of samples $v(\y)$ on a deterministic set of distinct collocation points 
$\GLmg$ 
having cardinality $M_L = \#\GLmg$. 
Applying $\ILmg[\cdot]$ from \eqref{eq:SC_approximation_operator} to the semi-discrete solution $\uh(x,\y)$ of problem \eqref{eq:semi_discrete_SFEM}, we obtain the fully discrete SC approximation
\begin{align}
\label{eq:fully_discrete_SC_soln}
\uSC(x,\y) = \ILmg[ \uh ] (x,\y). 
\end{align}

\paragraph{One-dimensional abscissas}
In this effort, we use three examples for constructing the fully discrete approximation. 
The first is that of a fully-nested rule constructed on the
Clenshaw-Curtis choice of abscissas \cite{ClenshawCurtis1960} with function $g_\text{TD}(\p)$ from \eqref{eq:g_defs} and an isotropic growth rule $\m=(m,\ldots,m)$ with $m$ given by
\begin{align}
\label{eq:SC_Smolyak_m}
m(1) = 1, \quad m(l_n)=2^{l_n-1}+1 \quad \text{for} \quad l_n>1, 
\end{align}
This is the classical Smolyak sparse-tensorization construction \cite{Smolyak1963}, and here the choice of $\m$ corresponds to a doubling growth rule that leads to a nested sequence of multi-dimensional grids, e.g., $\GLmg \subset \GLpmg$. 
On the other hand, 
we can construct a sparse-Smolyak approximation on the Gauss-Legendre abscissas corresponding to the zeros of the Legendre polynomials $\{\Psi_\p\}$, 
as defined in \S \ref{sec:SG}. When the points are grown isotropically according to the linear growth rule with $\m=(m,\ldots,m)$ and $m$ defined as 
\begin{align}
\label{eq:SC_linear_growth}
m(l_n) = l_n \qquad \text{for} \;\;\; l_n \in\N,
\end{align}
and $g_{\text{TD}}(\p)$ from \eqref{eq:g_defs}, we obtain 
a grid that is not nested.
Another construction that yields a sequence of nested grids is that based on the Leja points, defined as the sequence of points satisfying $y_{k+1} := \text{argmax}_{y\in\Gamma_n} \prod_{j=1}^k |y - y_{j}|$ (see \cite{DeMarchi2004}). 
Here we take the Leja sequence of points with $g_\text{TD}$ from \eqref{eq:g_defs} and the isotropic linear growth function $\m$ from \eqref{eq:SC_linear_growth}.

\subsection{Cost of solving the SCFEM systems}
\label{subsec:solving_SC}

To construct the fully discrete approximation with the SCFEM, we must solve $M_L$ distinct decoupled finite element systems, each dependent on a realization of the parameters $\y_k \in \GLmg$ for $k=1,\ldots,M_L$. Similar to the SGFEM, we can apply the PCG method to the solution of each system. 
Let $\Niterk$ be the number of iterations required by the CG method to solve the finite element system corresponding to $\y_k$ and $\Niterprek$ be the corresponding number of iterations when a preconditioner is used.
We are interested in choosing a suitable preconditioning strategy to decrease the total number of iterations 
$\NiterpreSC = \sum_{k=1}^{M_L} \Niterprek$
required to obtain the fully discrete approximation $\uSC$. We present a preconditioning strategy of choosing 
\begin{align}
\label{eq:level_zero_preconditioner}
\mathbf{P}_0 := \mathbf{A}(\y_1), 
\end{align}
with $\mathbf{A}(\y)$ from \eqref{eq:A_ij} the finite element stiffness matrix corresponding to the sample point $\y_1 \in \GLmg$, as the preconditioner for all of the individual finite element solutions. We refer to this choice of preconditioner as the level-zero preconditioner since it corresponds to the SC approximation at level $L=0$.

Since we apply CG to the solution of each individual finite element system, the work in floating point operations (FLOPs) required to obtain a fully discrete approximation with the SCFEM without a preconditioner is given by 
\begin{align}
\label{eq:SC_cost_solve}
\WSCsolve \approx \mathcal{O} \left( \Jh \right) * \sum_{k=1}^{M_L} \Niterk. 
\end{align}
On the other hand, the ``level-zero'' preconditioner induces an 
an additional matrix-vector product requiring $\mathcal{O}(\Jh)$ FLOPs per iteration when a sparse factorization of $\mathbf{P}_0$ is used. 
Hence the work of solving \eqref{eq:weak_problem} with PCG is given by 
\begin{align}
\label{eq:SC_cost_solve_precon}
\WSCsolvep \approx 2 * \mathcal{O} \left( \Jh \right) * \sum_{k=1}^{M_L} \Niterprek . 
\end{align}
Here, the reduction in work due to preconditioning will be seen in the number of iterations saved in each individual count $\Niterprek$ contributing to the sum.

\subsection{Comparing the explicit cost bounds of the SGFEM and SCFEM}
\label{subsec:complexity_comparison}

Given a particular ``sparse'' index set $\Lambdap$, we can find increasing functions $\bm{m}:\N_+^N \to \N_+^N$ and $g:\N_+^N \to \N_+$, and $L\in\N$ such that $\Lambdap = \LambdaLmg$ from \eqref{eq:general_index_set}. 
In this setting, we can either use Galerkin projection or construct an interpolant to obtain an approximation to $u$ in $\mathcal{P}_{\Lambdap}(\Gamma)$.
Let $u_{\Lambdap}$ denote the Galerkin projection of $u$ onto the space $\mathcal{P}_{\Lambdap}(\Gamma)$. 
Then 
we have the estimate
\begin{align*} 
\| u - u_{\Lambdap} \|_{L^2_{\varrho}(\Gamma; H_0^1(D))} \leq C_a \min_{v\in H_0^1(D)\otimes \mathcal{P}_{\Lambdap}(\Gamma)} \|u  - v \|_{L^2_{\varrho}(\Gamma; H_0^1(D))} 
\end{align*}
where $C_a > 0$ depends on the coefficient $a(x,\y)$ and the bounds from assumption $\Aone$. 
This estimate expresses optimality in the $L^2_{\varrho}(\Gamma)$ error of the Galerkin projection since $C_a$ does not grow with $\Lambdap$, and suggests that the Galerkin method is the best choice for approximating $u$ in the space $\mathcal{P}_{\Lambdap}(\Gamma)$. 
We can also define an interpolation operator $\ILmg : C^0(\Gamma)\to \mathcal{P}_{\Lambdap}(\Gamma)$,
and then we have the estimate
\begin{equation}
\begin{aligned}
\label{eq:Galerkin_and_interpolation_error} 
\| u - \ILmg [u] \|_{L^\infty_\varrho(\Gamma; H_0^1(D))} 
   & \leq (C_{\Lambda_L} + 1) \min_{v\in H_0^1(D)\otimes \mathcal{P}_{\Lambdap}(\Gamma)} \|u  - v \|_{L^\infty_\varrho(\Gamma; H_0^1(D))} \\ 
& = (C_{\Lambda_L} + 1) \| u - u_{\Lambdap} \|_{L^\infty_\varrho(\Gamma; H_0^1(D))} 
\end{aligned} 
\end{equation}
where $C_{\LambdaL}$ 
is the Lebesgue constant of $\ILmg$. 
A good interpolant will be one for which $C_{\LambdaL}$ grows moderately with $\#\LambdaLmg$. For example, it is known (see \cite{Dzjadyk1983,Galindo2015}) that for a one-dimensional Lagrange interpolation operator using a Clenshaw-Curtis rule, the Lebesgue constant is bounded by $\frac{2}{\pi} \log(m-1) + 1$, where $m$ is the number of points. 
For the SC method, we define SDOF to be the total number of points needed to construct the approximation. 
From \eqref{eq:Galerkin_and_interpolation_error}, if we only consider the number of SDOF needed to represent the solution, we expect the error for the Galerkin approximation to be much lower than the error in the interpolant. Indeed, this is reflected in our numerical results in Figures \ref{fig:BNTT_SDOF_and_cost} and \ref{fig:log_SDOF_and_cost}, and has been observed in previous comparisons \cite{BNTT_comp,Elman2011}.

However, if we are willing to change the space $\Lambdap$, e.g., adding more interpolation points to gain a more stable interpolant by changing $\m$ or changing which points are included in the set $\LambdaLmg$ by changing $g$, it might be possible to obtain an approximation with lower complexity to reach a given tolerance, despite having to solve more systems. 
Therefore, to properly compare the work involved in constructing $u_{\Lambdap}$ and $\ILmg[u]$, we consider the computational complexity of both methods, not in terms of SDOF, but in terms of floating point operations (FLOPs). For a chosen $\Lambdap$, this reduces to studying the complexity of the system resulting from Galerkin projections and the stability properties of the interpolant $\ILmg$. 

Let $\uSCCG$ denote the numerical solution to the fully discrete approximation $\uSC$ obtained with the SCFEM from \eqref{eq:fully_discrete_SC_soln} found by the PCG, and 
observe that we have a similar splitting to \eqref{eq:SG_full_error_with_projection} for the error in the approximation
\begin{align*}
\left\|u - \uSCCG \right\|_{\Htworho} 
                                                                      \leq \underbrace{\left\| u - \uh \right\|_{\Htworho}}_{\text{SC(I)}} + \underbrace{\left\| \uh - \uSC \right\|_{\Htworho}}_{\text{SC(II)}} + \underbrace{\left\| \uSC - \uSCCG \right\|_{\Htworho}}_{\text{SC(III)}} \numberthis \label{eq:SC_full_error}
\end{align*}
Note that unlike in the case of the SGFEM, the SCFEM does not require a further projection of the coefficient $a(x,\y)$, so that we do not need to consider the error $\|u - \ur \|_{\Htworho}$ from \eqref{eq:SG_full_error_with_projection}. 
In addition, we do not need to worry about well-posedness of the truncation as discussed in Remark \ref{rem:projection_well_posed}.
Similar to the complexity analysis for the SGFEM, we must choose $h\leq h_{\max}$ and $L\geq L_{\min}$ so that the errors $\|u - \uh\|_{\Htworho}$ from the finite element discretization and $\|\uh - \uSC\|_{\Htworho}$ from the SC interpolation are both bounded by $\varepsilon/3$. From this, a minimum tolerance $\tau_{\min}$ for the PCG solver
can be derived and the maximum number of PCG iterations, with a zero initial guess, can be estimated \cite{Galindo2015}. 
In what follows, we present a result, whose proof can be found in \cite[Theorem 4.7]{Galindo2015} that bounds
the number of PCG iterations in the context of the work estimate \eqref{eq:SC_cost_solve_precon}. Using this estimate we 
can compare the cost in FLOPs for the SCFEM with the SGFEM results from Theorems \ref{thm:SG_opt_work_poly} and \ref{thm:SG_opt_work_nonpoly} in the previous section.

\begin{thm}
\label{thm:SC_opt_work}
Let $u\in L^2_\varrho(\Gamma; H_0^1(D)\cap H^{s+1}(D))$ be the solution to \eqref{eq:model_problem}.
Then for $\eps>0$ arbitrary, the work of finding $\uSCCG$, the approximation to the fully discrete SC solution $\uSC$ from \eqref{eq:fully_discrete_SC_soln} found by PCG, denoted by $\WSCsolvep$, can be bounded by
\begin{align*}
\WSCsolvep \leq \; & 2 C_7 \left( \frac{3 \CFEM}{\varepsilon} \right)^{\frac{d}{s}} C_8 \left[ \log \left( \frac{3 \Csc}{\eps} \right) \right]^N \left[ C_9 + \frac{1}{\log 2} \log\log \left( \frac{3 \Csc}{\eps} \right)  \right]^{N-1} \numberthis \label{eq:SC_work_estimate} \\
& \times \frac{1}{\log \left(\frac{\sqrt{\bar{\kappa}}+1}{\sqrt{\bar{\kappa}} - 1} \right)} \left\{ \log \left( \frac{C_{10}}{\eps} \right) + C_{11} + 2N \log\log \left[ \frac{1}{rN} \log \left( \frac{3 \Csc}{\eps} \right) \right] \right\}.
\end{align*}
Here $\CFEM$ from Lemma \ref{lem:fem_rate}, $C_7$ from Theorem \ref{thm:SG_opt_work_poly}, and $C_8$, $C_9$, $C_{10}$, $C_{11}$, $\Csc$, and $r$ from \cite[Theorem 4.7]{Galindo2015} are positive constants independent of $\eps$. Moreover, we define $\bar{\kappa} = \sup_{\y\in\Gamma} \kappa(\y)$ where $\kappa(\y)$ is the condition number of the preconditioned system 
$\mathbf{\tilde P}_0^{-1} \mathbf{A}(\y)$ 
with $\mathbf{P}_0$ from \eqref{eq:level_zero_preconditioner}.
\end{thm}

Theorem \ref{thm:SC_opt_work} follows from the fact that $\NiterpreSC = \sum_{k=1}^{M_L} \Niterprek \leq N_{\text{zero}}$, where $N_{\text{zero}}$ is the number of iterations needed by the SCFEM with a PCG and a zero vector initial guess. Substituting the bound on $N_{\text{zero}}$ shown in \cite[Theorem 4.7]{Galindo2015} into the work estimate \eqref{eq:SC_cost_solve_precon} and using $J_{h_{\max}} = C_7 (3 \CFEM/\varepsilon)^{d/s}$ as in Theorem \ref{thm:SG_opt_work_poly}, puts the result in terms of FLOPs.
Given Theorem \ref{thm:SC_opt_work} we see that the work of obtaining the fully discrete approximation with the SCFEM with the PCG method is asymptotically bounded by:
\begin{align}
\label{eq:SC_asymptotic_estimate}
\mathcal{O} \underbrace{
    \vphantom{\left( \frac{ \log\left( \frac{1}{\eps} \right) }{ \log \left( \frac{ \sqrt{\kappar} + 1 }{ \sqrt{\kappar} - 1} \right) } \right)}
    \left( \frac{1}{\eps} \right)^{\frac{d}{s}}}_{\text{(SC.1)}}  
  \underbrace{
    \vphantom{\left( \frac{ \log\left( \frac{1}{\eps} \right) }{ \log \left( \frac{ \sqrt{\kappar} + 1 }{ \sqrt{\kappar} - 1} \right) } \right)}
    \left[\log \left( \frac{1}{\eps} \right)\right]^N
    \left[\log \log \left( \frac{1}{\eps} \right)\right]^{N-1}
}_{\text{(SC.2)}}
\underbrace{\left( \frac{ \log\left( \frac{1}{\eps} \right) 
            }{ \log \left( \frac{ \sqrt{\bar{\kappa}} + 1 }{ \sqrt{\bar{\kappa}} - 1} \right) } \right)}_{\text{(SC.3)}}
\end{align}
where (SC.1), (SC.2), and (SC.3) correspond to the work required by the finite element, SC interpolant, and PCG methods, respectively. 
Since the costs associated with the finite element discretization are the same for both methods, we focus only on the costs associated with the SG projection and the SC interpolation, coupled with the costs of the PCG method.
In particular, as the work required by the SG approximation from (SG.2) of \eqref{eq:SG_asymptotic_estimate} has different bounds depending on whether the coefficient is a fixed order polynomial or is given by a total degree orthogonal expansion having order depending on $\eps$, we now provide a comparison in both of these cases.

{\em Comparison in the affine and non-affine polynomial cases, e.g., Examples \ref{exmp:affine_coefficient} and \ref{exmp:polynomial_coefficient}.}
The terms (SG.2) from \eqref{eq:SG_asymptotic_estimate} and (SC.2) from \eqref{eq:SC_asymptotic_estimate} are asymptotic estimates of the number of coupled and decoupled finite element systems that must be solved by the SGFEM and SCFEM to construct the stochastic approximation, respectively. 
For the coefficients from Examples \ref{exmp:affine_coefficient} and \ref{exmp:polynomial_coefficient}, (SG.2) from \eqref{eq:SG_asymptotic_estimate} for the SGFEM is $\mathcal{O}([\log(1/\eps)]^N)$. Hence, our analysis shows that in these cases, the number of coupled finite element systems in the SGFEM matrix is smaller than the number of decoupled finite element systems required by the SCFEM by a factor of $(\log\log(1/\eps))^{N-1}$. This difference is enough to suggest that, if the condition numbers of both the preconditioned coupled system from the SGFEM and the preconditioned individual systems from the SCFEM are of the same order, then the SGFEM will outperform the SCFEM in terms of minimum work required to obtain a fully discrete approximation. 
In fact, whenever $a(x,\y)$ is a general non-affine, polynomial coefficient, e.g., Example \ref{exmp:polynomial_coefficient}, having fixed order $\rbar<\infty$, the complexity of matrix-vector products involving $\mathbf{K}$ from \eqref{eq:fully_discrete_SG_algebraic_KuF} is approximately $\mathcal{O}(\Jh \Mp M_{\rbar} \min\{ 2^{\rbar}, M_{\lceil \rbar/2 \rceil}\})$
from Remark \ref{rem:SG_affine_and_poly_nnz_complexity}. 
Hence, our analysis shows that 
when $\mathcal{O}(M_{\rbar} \min\{ 2^{\rbar}, M_{\lceil \rbar/2 \rceil}\}) < \mathcal{O}((\log\log(1/\eps))^{N-1})$, which, in limit as $\eps\to 0$, is always the case, and the condition numbers of both systems are of the same order, the SGFEM will outperform the SCFEM. 
However, in practical applications, it may require unrealistically small tolerance $\eps$ to see this when $\rbar$ is large.

\begin{figure}[!htb]
\begin{center}
\includegraphics[scale = 0.3,trim=0mm 0mm 15mm 10mm,clip=true]{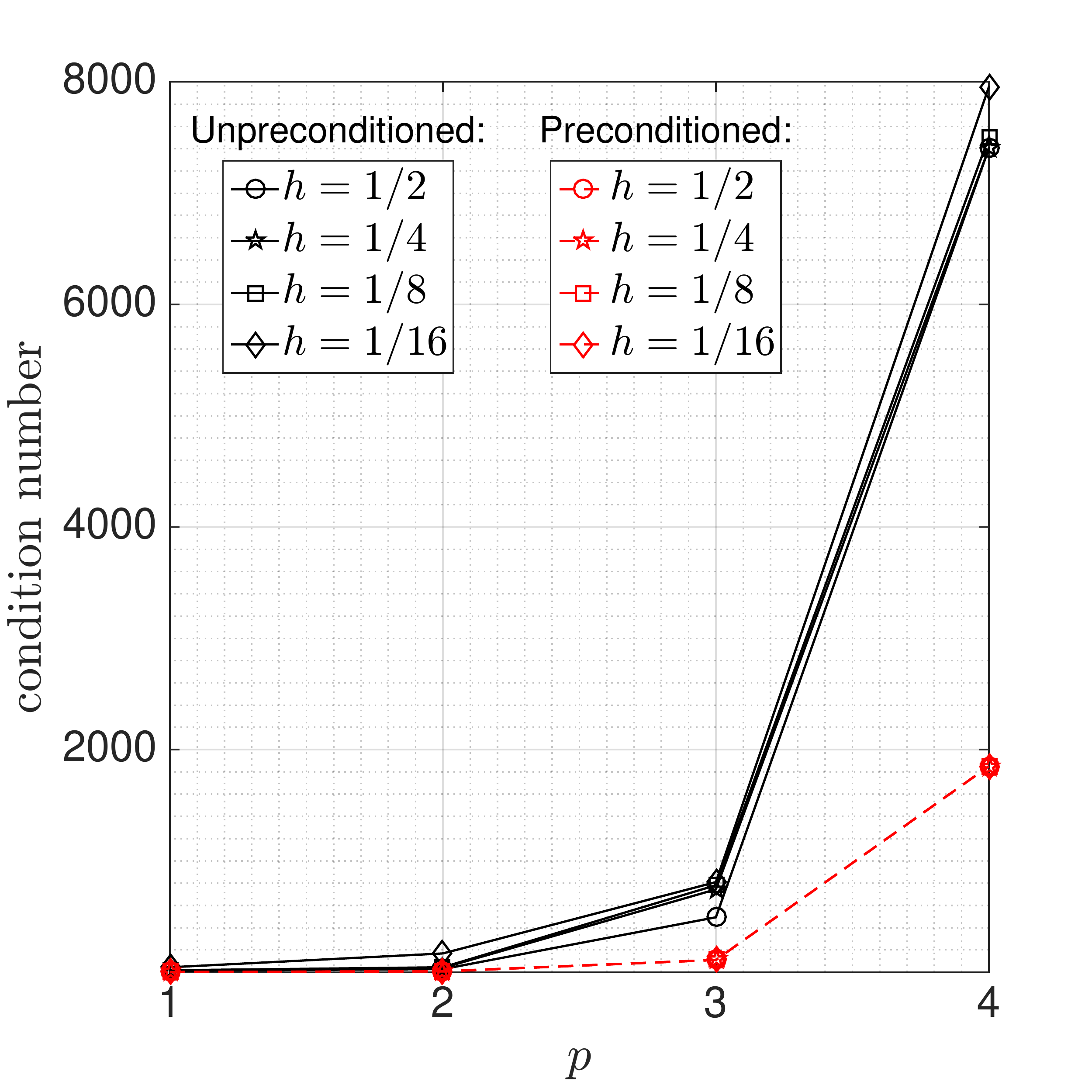}
\caption{Condition numbers of both the unpreconditioned matrix $\mathbf{K}_r$ and the preconditioned matrix $\mathbf{P}^{-1} \mathbf{K}_r$ with decreasing finite element mesh parameter $h$ and $r=p$ for the coefficient \eqref{eq:log_truncation} from \S\ref{subsec:transcendental_ex} with $N=4$ and $L_c=1/2$.}
\label{fig:cond_nums}
\end{center}
\end{figure}

{\em Comparison in the non-affine, transcendental case, e.g., Example \ref{exmp:transcendental_coefficient}}. In the case that $a(x,\y)$ is a non-affine, transcendental functon of the random parameters, e.g., Example \ref{exmp:transcendental_coefficient}, the estimate for (SG.2) is $\mathcal{O}([\log(1/\eps)]^{3N})$. Here our analysis shows that the number of coupled finite element matrices present in the SG system $\Kr$ from \eqref{eq:fully_discrete_SG_algebraic_sep_KuF} dominates the number of decoupled finite element systems needed by the SCFEM by a factor of $[\log(1/\eps)]^{2N}$. 
Note that, as in the case of the SGFEM, the term (SC.3) has a dependence on the condition numbers of the preconditioned finite element systems through the bound $\bar{\kappa} = \sup_{\y\in\Gamma} \kappa(\y)$.
For the unpreconditioned systems $\mathbf{A}(\y)$, the condition numbers can be bounded by
$\kappa(\mathbf{A}(\y)) \leq \left(C_\kappa/h\right)^2$ 
for every $\y\in\Gamma$, following from assumption $\Aone$ and the quasi-uniformity of the mesh $\mathcal{T}_h$.
However, if we use the exact inverse of $\mathbf{P}_0$ when preconditioning the SCFEM systems, 
the condition numbers are bounded independent of $h$ and $\y\in\Gamma$, since in this case $\bar{\kappa}$ 
is independent of mesh size $h$ and level $L$. 
Hence the work required by the PCG method when solving the SCFEM systems is dependent on $\varepsilon$ only through the term $\log(1/\eps)$.
On the other hand, if we use the exact inverse of $\mathbf{P}$ when preconditioning the SGFEM system,
the condition number $\kappartilde$ 
can be bounded by 
\begin{align*}
\kappartilde \leq \frac{1 + \overline{\tau}_r}{1 - \underline{\tau}_r},
\end{align*}
where $\underline{\tau}_r$ and $\overline{\tau}_r$ are defined in \eqref{eq:Pinv_Kw_eigenbounds}, hence depend on 
$\eps$ when $r$ is chosen to satisfy $r \geq r_{\min}$ from Lemma \ref{lem:SG_optimal_parameters}. 
Figure \ref{fig:cond_nums} plots the condition numbers of both the unpreconditioned matrix $\mathbf{K}_r$ and the preconditioned matrix $\mathbf{P}^{-1} \mathbf{K}_r$ with decreasing finite element mesh parameter $h$ for the coefficient $a(x,\y)$ given in \eqref{eq:log_truncation} from \S\ref{subsec:transcendental_ex} with $N=4$, $L_c = 1/2$, and letting $r=p$ with $p$ increasing. There we see that the dependence on $h$ has been removed by applying $\mathbf{P}^{-1}$, but as $p$ increases, we see a corresponding increase in the condition number.
Other preconditioners than $\mathbf{P}$ may be used to reduce the dependence on $r$, e.g. \cite{Ullmann2010}, but then their associated costs must be accounted for in the work estimate \eqref{eq:SG_cost_solve_precon} as well. 
However, even if the condition numbers of both the preconditioned coupled SG system and the preconditioned decoupled SC systems are of the same order, 
the additional work required to solve the coupled systems induced by the nonlinearity of the coefficient 
makes it difficult to see how the SGFEM can compete with the SCFEM. 

\section{Numerical examples}
\label{sec:num_ex}

In this section, we provide illustrative numerical examples comparing the complexity of the SGFEM in the three cases of Examples \ref{exmp:affine_coefficient}, \ref{exmp:polynomial_coefficient}, and \ref{exmp:transcendental_coefficient}. We then compare these results with SCFEM and the results of the theoretical complexity comparison of the previous section. We solve the model problem \eqref{eq:model_problem},
on the unit square $D = [0,1]^2$. 
For a general coefficient $a(x,\y)$ we do not know the exact solution to \eqref{eq:model_problem}. Hence we check the convergence against a ``highly enriched'' approximation, which we consider close enough to the exact one. 
To construct this ``exact'' solution $u_\text{ex}(x,\y)$, we make use of the isotropic SCFEM based on Clenshaw-Curtis abscissas using the level $L_\text{ex}$. 
We approximate the computational error for the SGFEM with orders $p=0,1,2,\ldots,p_{\max}$ and for the SCFEM with levels $L=0,1,2,\ldots,L_{\max}$ as
\begin{align}
\label{eq:error_metrics}
\| \bbE[\eps_\text{SG}] \|_{\ell^\infty} \approx \| \bbE[u_\text{ex} - \uSGCG] \|_{\ell^\infty} \qquad \text{and} \qquad \| \bbE[\eps_\text{SC}] \|_{\ell^\infty} \approx \| \bbE[u_\text{ex}-\uSCCG] \|_{\ell^\infty},
\end{align}
where $\uSGCG$ and $\uSCCG$ are the fully discrete approximations \eqref{eq:fully_discrete_SG_soln} and \eqref{eq:fully_discrete_SC_soln}, respectively, found by the PCG method, described in \S\ref{sec:SG} and \S\ref{sec:SC_comparison}. 
In \S\ref{subsec:transcendental_ex}, we measure $\|\bbE[\eps_\text{SG}] \|_{\ell^\infty} \approx \| \bbE[u_\text{ex} - \uSGCGr] \|_{\ell^\infty}$ where $\uSGCGr$ denotes the solution of  \eqref{eq:fully_discrete_SG_algebraic_sep_KuF} with
 the projected coefficient $a^r(x,\y)$.

As stated in \S\ref{subsec:solving_SG} and \S\ref{subsec:solving_SC}, we use PCG with the mean-based preconditioner for SGFEM and the level-zero preconditioner for the SCFEM. 
Hence, we believe this puts both methods at a similar starting point for comparison, if not providing a slight advantage for the SGFEM. 
With these choices, the complexity results are presented in terms of the work estimates \eqref{eq:SG_cost_solve_precon} and \eqref{eq:SC_cost_solve_precon}, respectively. 
The amount of work to reach a given error in PCG is also dependent on the tolerance used by the solver. If the tolerance is too small, we may see that the PCG method ``over-resolves'' the solution. 
To ensure that 
we do not over-resolve either solution, we set the tolerance of the solvers to be $\|\bbE[\eps_\text{SG}]\|_{\ell^\infty}/10$ and $\|\bbE[\eps_{\text{SC}}]\|_{\ell^\infty}/10$ respectively, where these quantities are first estimated for each order $p$ and level $L$ using a tolerance of $1.0\times 10^{-12}$. 
In practice, we find that this does not affect the convergence results much. 

In all three examples, we use the SG approximation constructed in terms of the orthonormal Legendre polynomials $\{\Psi_\p\}_{\p\in\Lambdap}$ for given index sets $\Lambdap$. In the presentation of the results that follow, we use the following abbreviations. For the SGFEM, we use: ``SG-TD'' to denote the approximation in the total degree subspace $\mathcal{P}_{\LambdapTD}(\Gamma)$ with $\LambdapTD$ given in \eqref{eq:index_set_examples}, and ``SG-SM'' to denote the approximation in the sparse Smolyak subspace $\mathcal{P}_{\LambdapSM}(\Gamma)$ with $\LambdapSM$ given in \eqref{eq:index_set_examples}. For the SCFEM, we use: 
``SC-GL'' and ``SC-LJ''  to denote the Smolyak approximation constructed on Gauss-Legendre abscissas and the Leja approximation constructed on Clenshaw-Curtis abscissas, both defined in terms of $g_\text{TD}$ and $\m$ given in \eqref{eq:g_defs} and \eqref{eq:SC_linear_growth}, respectively, and 
``SC-CC'' to denote the Smolyak approximation constructed on Clenshaw-Curtis abscissas with $g_\text{SM}$ and $\m$ given in \eqref{eq:g_defs} and \eqref{eq:SC_Smolyak_m}.

\subsection{Piecewise affine coefficients}
\label{subsec:linear_ex}

\begin{figure}[!htb]
\begin{center}
\includegraphics[scale = 0.7]{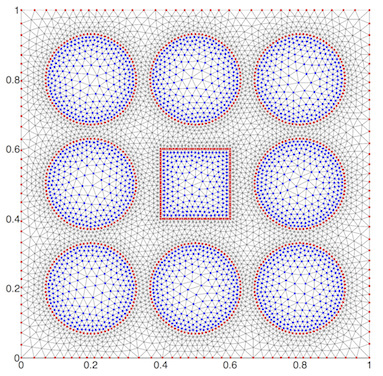}
\includegraphics[scale = 0.7]{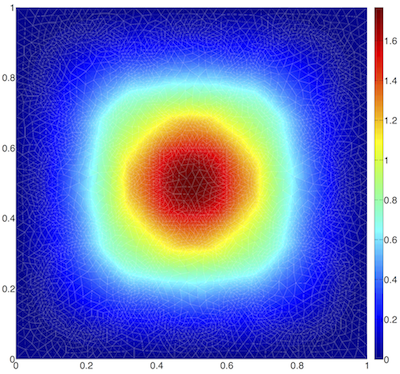}
\caption{{\bf Left:} a triangulation of the domain $D$ with circular and square inclusions. Red nodes highlight the boundary of an inclusion or the domain $D$, blue nodes highlight nodes on the interior of an inclusion. {\bf Right:} the expected value of the solution of \eqref{eq:model_problem} with stochastic conductivity coefficient \eqref{eq:subdom_problem}.}
\label{fig:BNTT_geometry_and_soln}
\end{center}
\end{figure}

One common example in engineering and the physical sciences is that of isotropic thermal diffusion problem with a stochastic conductivity coefficient. Consider a partitioning of $D=[0,1]^2$ into 8 circular inclusions arrayed about 1 square inclusion as in Figure \ref{fig:BNTT_geometry_and_soln}. 
We present the following example from \cite{BNTT_comp}, where the coefficient was given by
\begin{align}
\label{eq:subdom_problem}
a(x,\y) = b_0(x) + \sum_{n=1}^8 y_n \chi_n(x), 
\end{align}
with $b_0 = 1$, and $y_n \sim \mathcal{U}(-0.99,-0.2)$.
Here, $\chi_n$ are indicator functions corresponding to the 8 circular inclusions of radius $r=0.13$. In this example, we also set the forcing term to be 
\begin{align}
\label{eq:subdom_forcing}
f(x) = 100\chi_F(x), 
\end{align}
where $F = [0.4,0.6]^2$, is the square inclusion centered in $D$ with side length $0.2$. Figure \ref{fig:BNTT_geometry_and_soln} shows the expected value of the solution to this problem.
To solve \eqref{eq:model_problem} with the coefficient \eqref{eq:subdom_problem} and forcing function \eqref{eq:subdom_forcing}, we use a piecewise linear finite element basis in the deterministic space over a nonuniform mesh $\mathcal{T}_h$. Here, the nodes of $\mathcal{T}_h$ are adapted to the geometry of our problem, that is, we fix the nodes that lie on the boundaries of the inclusions in our domain. From this fixed boundary data, we then use the {\tt distmesh} MATLAB program \cite{Persson2004} to generate a non-degenerate triangulation that adequately resolves the details of our subdomain geometry. We further specify the subsets of the total set of nodes that belong to each geometric inclusion, and to the boundaries of the inclusion, so that the interface conditions for the coefficient may be correctly applied. The final mesh consists of 10,604 elements, 5,377 total nodes, and 5,229 unknowns.

\begin{figure}[!htb]
\begin{center}
\includegraphics[clip=true,trim=02mm 02mm 16mm 12mm,width=0.326\textwidth]{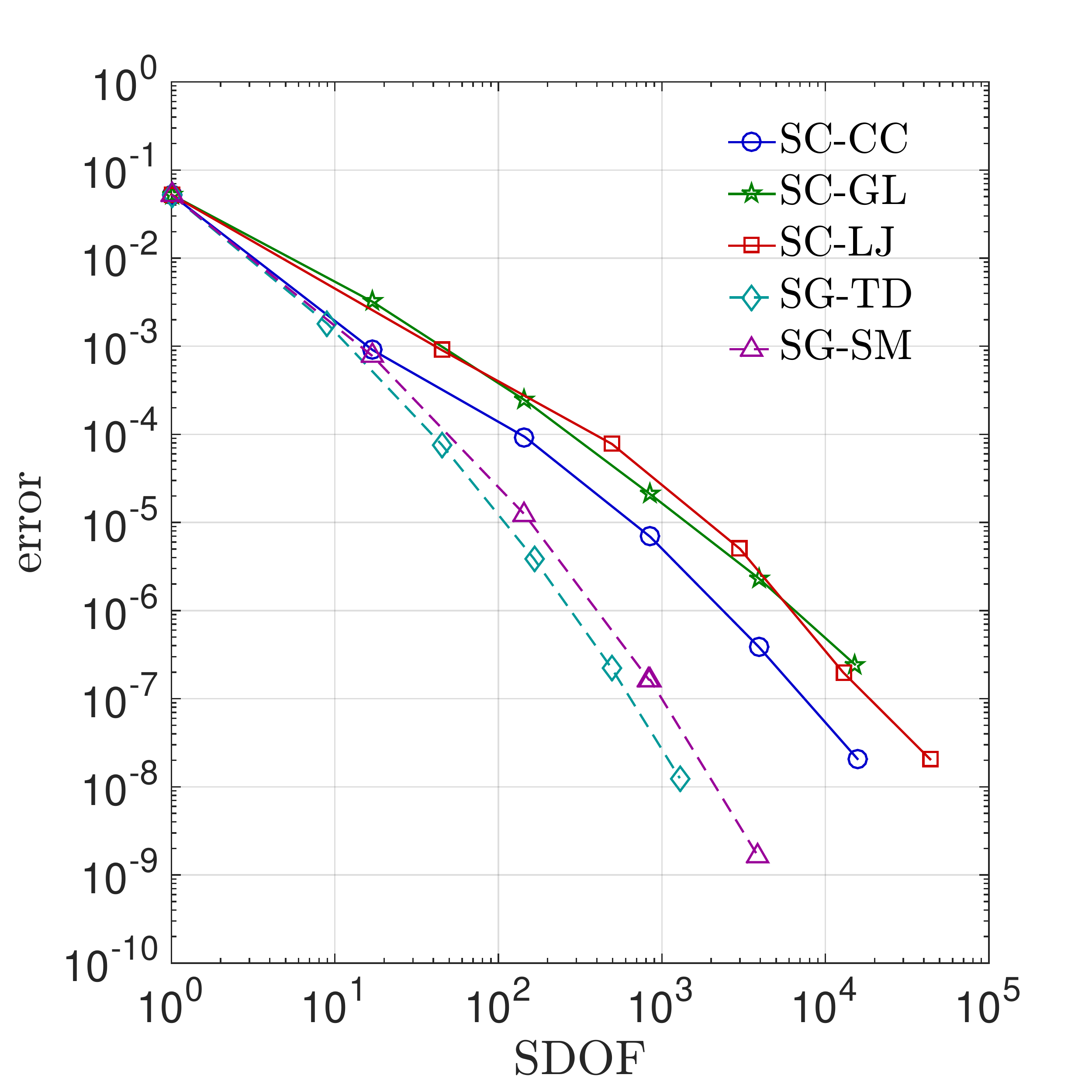}
\includegraphics[clip=true,trim=02mm 02mm 16mm 12mm,width=0.326\textwidth]{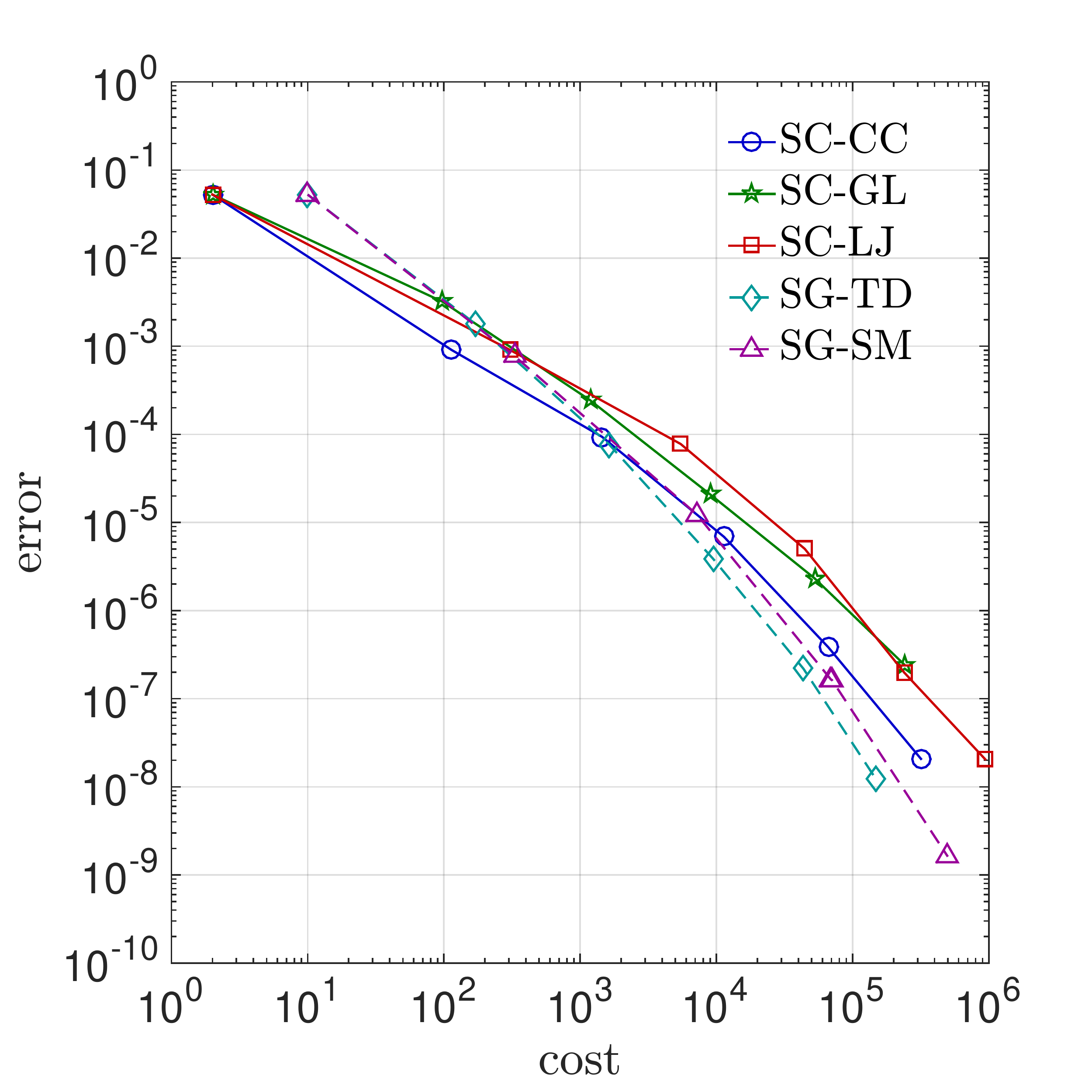}
\caption{{\bf Left:} Error versus SDOF in solving problem \eqref{eq:model_problem} with coefficient \eqref{eq:subdom_problem} and forcing \eqref{eq:subdom_forcing}. 
    {\bf Right:} Error versus computational cost with the work estimates given in \eqref{eq:SG_cost_solve_precon} and \eqref{eq:SC_cost_solve_precon} based on total number of matrix-vector products used by the CG method. 
}
\label{fig:BNTT_SDOF_and_cost}
\end{center}
\end{figure}

The coefficient \eqref{eq:subdom_problem} is an example of a coefficient $a(x,\y)$ having affine dependence on the parameters, e.g., Example \ref{exmp:affine_coefficient}.
Figure \ref{fig:BNTT_SDOF_and_cost} displays the convergence of the stochastic Galerkin and collocation methods against the total number of SDOF. 
For the SGFEM we take the SDOF to be the cardinality of the set $\Lambdap$ used in constructing the fully discrete approximation $\uSG$ from \eqref{eq:fully_discrete_SG_soln} by solving \eqref{eq:fully_discrete_SG_problem}, and for the SC method we take the SDOF to be the number of points $\#\GLmg$ corresponding to an index set $\LambdaLmg$ used in constructing the fully discrete approximation $\uSC$ from \eqref{eq:fully_discrete_SC_soln}.
From the discussion of \S \ref{subsec:complexity_comparison}, we expect to see that the approximation obtained with the SGFEM requires fewer SDOF than the SCFEM to achieve the same error, and this is indeed the observed result. For example, both the SG-TD and SC-LJ approximations require the same number of SDOF, but the error of the SC-LJ approximation is much higher.
This, of course, is a consequence of the estimate \eqref{eq:Galerkin_and_interpolation_error}, where the errors of the SC approximations are bounded above by their respective Lebesgue constants against the best-approximation error in the space $\mathcal{P}_{\LambdaLmg}(\Gamma)$.

Figure \ref{fig:BNTT_SDOF_and_cost} also displays the convergence of both methods in terms of error versus the total computational cost of solving the system with the work estimates of \eqref{eq:SG_cost_solve_precon} and \eqref{eq:SC_cost_solve_precon}, respectively.
Here, we compute the error in $\|\bbE[\eps_\text{SG}]\|_{\ell^\infty}$ and $\|\bbE[\eps_\text{SC}]\|_{\ell^\infty}$ as given in (60) and measure the cost as the number of $\mathcal{O}(J_h)$ matrix vector products required by both methods which are explicitly counted as
$$N_{\textnormal{iter}}^{\textnormal{pSG}} * \left(M_p + \sum_{\rb\in\Lambdar} \textnormal{nnz}(\mathbf{G}_\rb) \right) =
N_{\textnormal{iter}}^{\textnormal{pSG}} * (M_p + \mathcal{M}(p,r))$$
in the code for the SGFEM and $2 * \sum_{k=1}^{M_L} N_{\textnormal{iter}}^{\textnormal{p}(k)}$ in the code for the SCFEM. 
Our analysis shows the work corresponding to the SG discretization for SG-TD is asymptotically bounded by $\mathcal{O}([\log(1/\eps)]^N)$ while the analysis from \cite{Galindo2015} shows that the work corresponding to the SC discretization for SC-CC is asymptotically bounded by $\mathcal{O}([\log ( 1/\eps)]^N [\log \log ( 1/\eps )]^{N-1})$. 
This closely matches the results of the numerical experiments in Figure \ref{fig:BNTT_SDOF_and_cost}, where it can be seen that for polynomial order $p\geq 2$, the SG-TD approximation yields the best results with the least computational cost for this problem.

\subsection{Polynomial coefficients}
\label{subsec:polynomial_ex}

The next example we present is that of a polynomial function of the random paramters $\y$, e.g.
the coefficient from Example \ref{exmp:polynomial_coefficient}. We consider the following function 
\begin{align*}
a(x,\y)  = 5 + \sum_{|\rb|\leq \rbar} e^{-1.5 |\rb|} \varsigma_\rb(x) \y^\rb \numberthis \label{eq:num_polynomial_coefficient}, \;\;\;
\varsigma_\rb(x) = \left\{ \begin{array}{rl} \sin\left( |\rb|\pi x_1 \right)\cos\left( |\rb|\pi x_2 \right) & \text{ if $|\rb|$ is even,} \\ \cos\left( |\rb|\pi x_1 \right)\sin\left( |\rb|\pi x_2 \right) & \text{ if $|\rb|$ is odd,} \end{array} \right.
\end{align*}
with $y_n \sim \mathcal{U}(-1,1)$ for all $n=1,\ldots,N$ and forcing term $f(x) = 1$ $\forall x\in\overline{D}$.
For the results that follow we fix $N=4$ and study the convergence of the SGFEM and SCFEM in the cases $\rbar = 1,3,7$ in \eqref{eq:num_polynomial_coefficient}. 
As in \S \ref{subsec:linear_ex}, we set the finite element space for the spatial discretization to be the span of piecewise linear polynomials, but here we use a uniform triangulation of $D$ with $4,934$ elements and $2,340$ spatial unknowns. 

\begin{figure}[!htb]
\begin{center}
\includegraphics[clip=true,trim=03mm 03mm 16mm 12mm,width=0.326\textwidth]{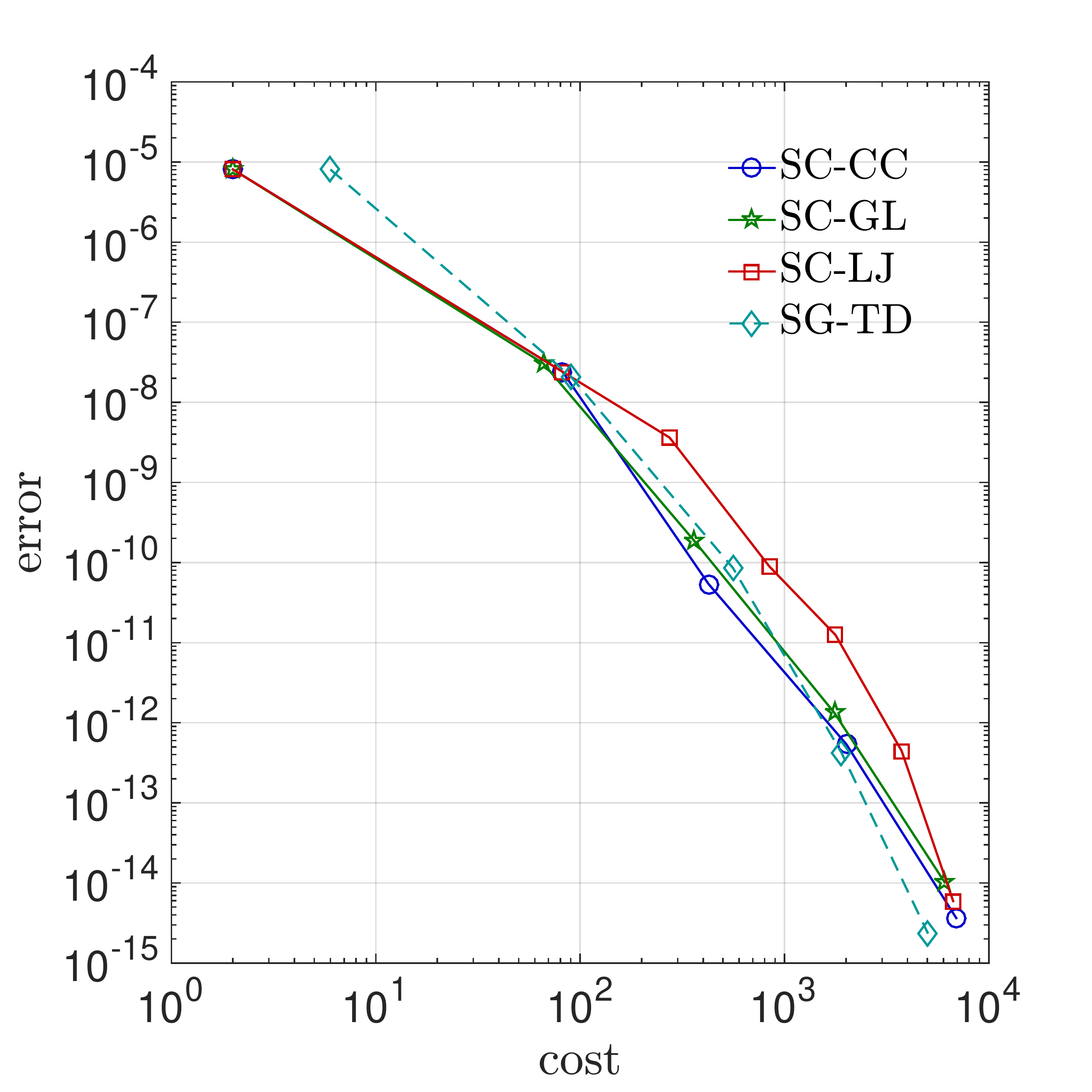}
\includegraphics[clip=true,trim=03mm 03mm 16mm 12mm,width=0.326\textwidth]{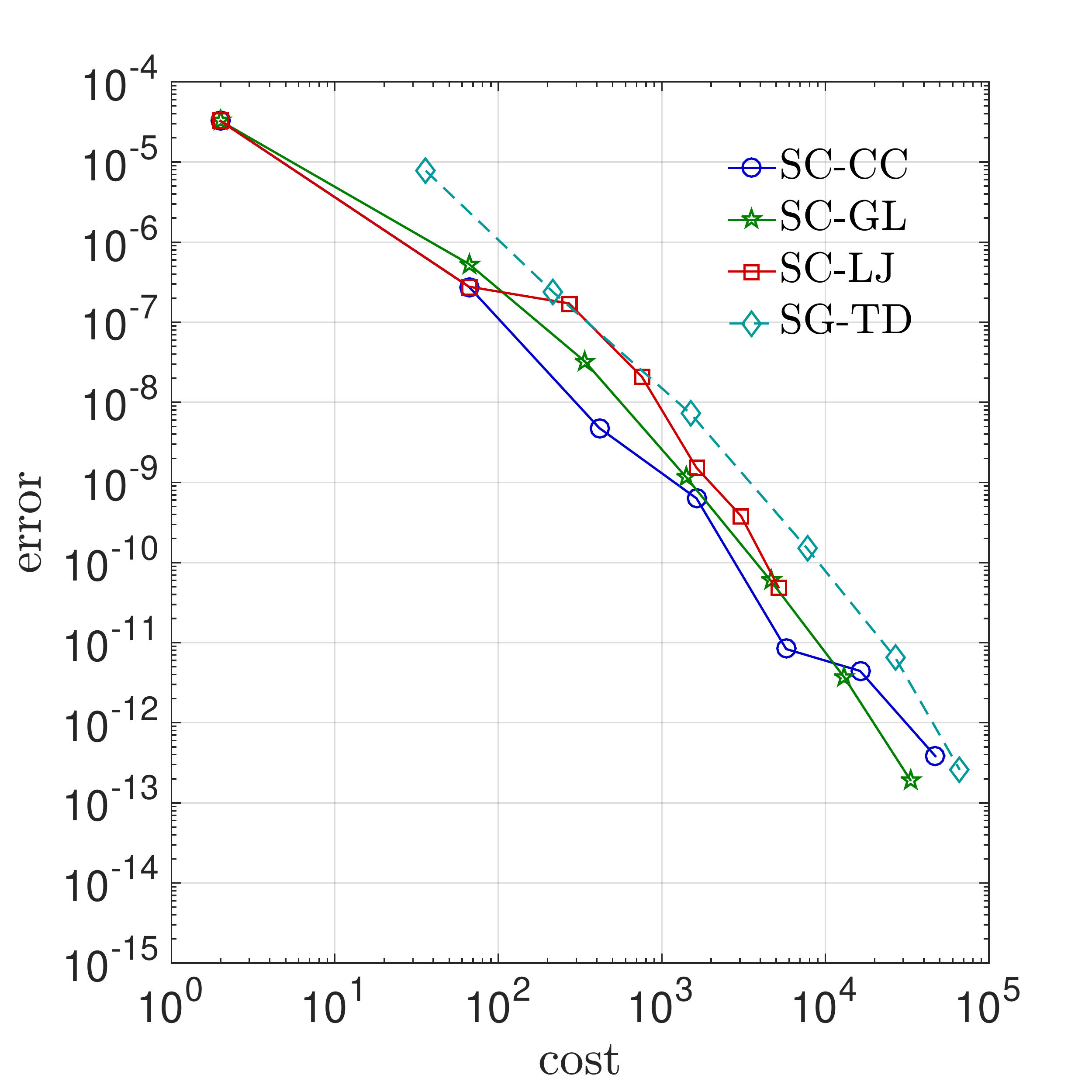}
\includegraphics[clip=true,trim=03mm 03mm 16mm 12mm,width=0.326\textwidth]{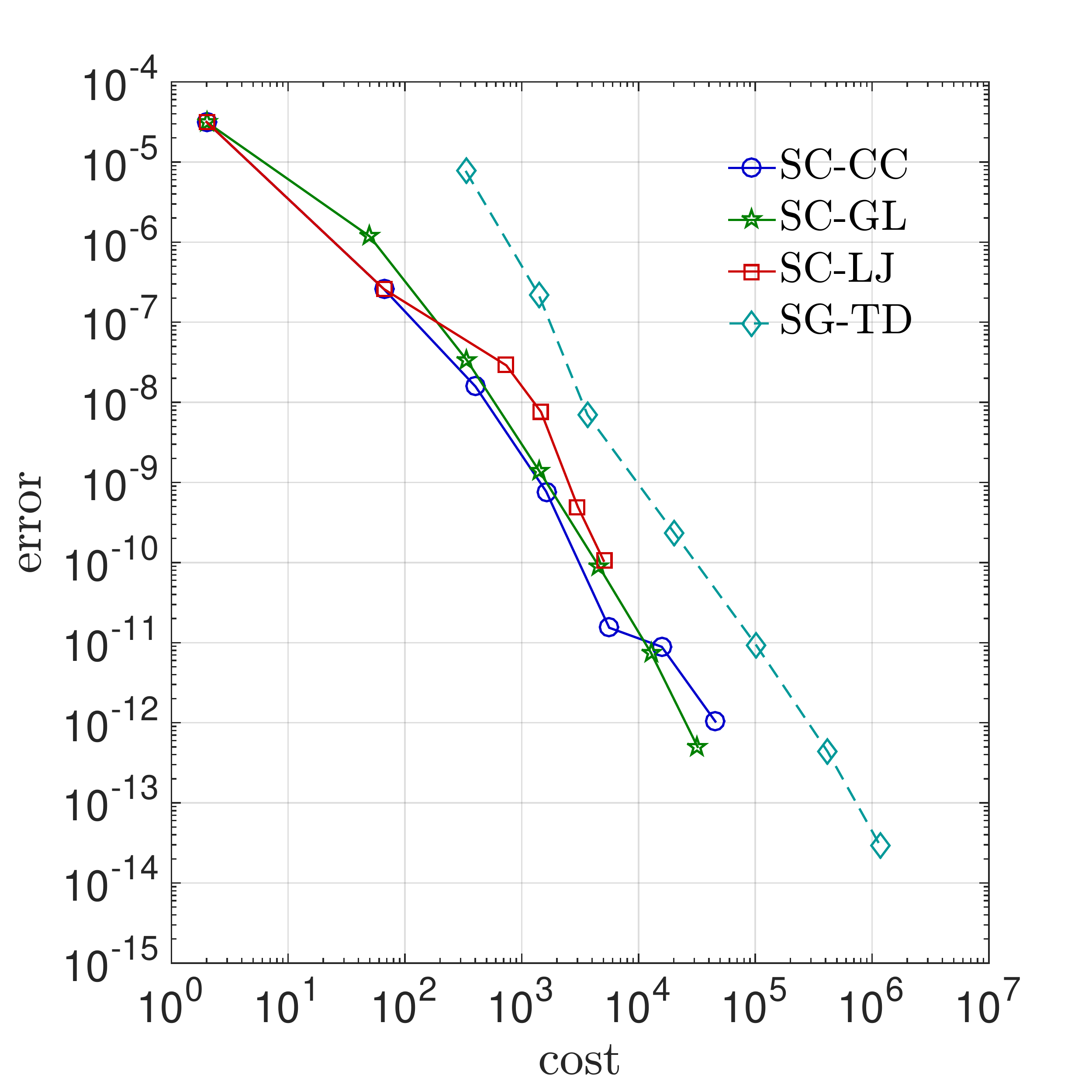}
\caption{
Error versus cost for solving problem \eqref{eq:model_problem} with coefficient \eqref{eq:num_polynomial_coefficient} having $\rbar=1$ (left), $\rbar=3$ (middle), and $\rbar=7$ (right), with forcing $f(x)=1$. The cost, given in \eqref{eq:SG_cost_solve_precon} and \eqref{eq:SC_cost_solve_precon}, is based on total number of matrix-vector products used by the PCG method. 
}
\label{fig:polynomial_cost}
\end{center}
\end{figure}

Figure \ref{fig:polynomial_cost} displays the convergence of the SGFEM and SCFEM in terms of error versus the total computational cost of solving the system with the work estimates of \eqref{eq:SG_cost_solve_precon} and \eqref{eq:SC_cost_solve_precon}.
Here, we compute the error in $\|\bbE[\eps_\text{SG}]\|_{\ell^\infty}$ and $\|\bbE[\eps_\text{SC}]\|_{\ell^\infty}$ as given in \eqref{eq:error_metrics}.
As we increase the order $\rbar$ in \eqref{eq:num_polynomial_coefficient}, we see that the work for the SGFEM increases, corresponding to the decreasing sparsity of the matrix $\mathbf{K}$ from \eqref{eq:fully_discrete_SG_algebraic_KuF}. 
Here, the work of matrix-vector multiplications with $\mathbf{K}$ are of the order $\mathcal{O}(J_h M_p M_{\rbar} \min \{2^{\rbar}, M_{\lceil \rbar/2 \rceil} \})$, where $M_{\rbar} \min \{2^{\rbar}, M_{\lceil \rbar/2 \rceil}\}$ is a large constant that grows rapidly with $\rbar$. As a result, we see that for $\rbar = 1$, the SGFEM outperforms the other methods for $p\geq 4$. However, for $\rbar = 3,7$, the extra work of the matrix-vector multiplications of the coupled SG system
dominates the overall convergence.
We also observe that the convergence rate of the SGFEM does not change in these cases, as discussed in the comparison in \S\ref{subsec:complexity_comparison}.

\subsection{Transcendental coefficients}
\label{subsec:transcendental_ex}

The next example we present is that of a random coefficient defined in terms of the truncated \KL expansion of the function $\log(a(x,\y) - a_{\min})$, for $a_{\min}>0$. This example represents a commonly used transcendental function of the physical and random parameters, e.g., Example \ref{exmp:transcendental_coefficient}, and is often presented in the context of enforcing the positivity of $a(x,\y)$ required in assumption $\Aone$. Coefficients of this type are commonly found in groundwater flow models. For these models, the permeability can exhibit large variance within each layer of sediment, and as a result are better represented on a logarithmic scale. 
We recall the problem of solving \eqref{eq:model_problem} with a coefficient having one-dimensional (layered) spatial dependence and a deterministic load $f(x_1,x_2,\omega) = 2 \cos(x_1)\sin(x_2)$ as studied in \cite{WebsterIsotropic,WebsterAnisotropic}, where $a(x,\y)$ was given by
\begin{gather*}
         \log(a(x,\omega)-0.5)  = 1 + y_1(\omega)\left(\frac{\sqrt{\pi}L}{2}\right)^{1/2} + \sum_{n=2}^N \zeta_n \varphi_n(x) y_n(\omega), \numberthis \label{eq:log_truncation} \\
   \zeta_n : = (\sqrt{\pi} L)^{1/2} \exp\left( \nicefrac{-\left( \left\lfloor \frac{n}{2} \right\rfloor \pi L\right)^2}{8} \right), \text{ for } n>1, 
\;\;\; \;\;\varphi_n(x) := \left\{ \begin{array}{rl} \sin\left( \left\lfloor\frac{n}{2}\right\rfloor \pi x_1/L_p \right),  \text{ if $n$ is even,} \\ \cos\left(\left\lfloor\frac{n}{2}\right\rfloor \pi x_1/L_p \right),  \text{ if $n$ is odd.} \end{array} \right.
\end{gather*}
Here, $\{y_n(\omega)\}_{n=1}^\infty$ are independent random variables uniformly distributed in $[-\sqrt{3},\sqrt{3}]$ with zero mean and unit variance. For $x_1\in[0,b]$, let $L_c$ be a desired physical correlation length for the random field $a(x,\y)$, chosen so that the random variables $a(x_1,\omega)$ and $a(x_1',\omega)$ become essentially uncorrelated for $|x_1-x_1'|\gg L_c$. Also, let $L_p = \max\{b,2L_c\}$ and $L=L_c/L_p$.
Expression \eqref{eq:log_truncation} represents a possible truncation of a one-dimensional random field with stationary covariance,
\begin{align*}
\cov[\log(a-0.5)](x_1,x_2) & 
                              = \exp\left(\frac{-(x_1-x_2)^2}{L_c^2}\right) .
\end{align*}
Direct integration with the coefficient $a(x,\y)$ from \eqref{eq:log_truncation} yields a fully-block dense linear system $\mathbf{K}$ from \eqref{eq:fully_discrete_SG_algebraic_KuF} that is computationally infeasible to solve
\cite{Ernst2010,Matthies2005,Ullmann2010,Ullmann2012}. 
The purpose of this example is to highlight the difficulties of obtaining a fully discrete approximation with the SGFEM in this case.

As in the previous example in \S \ref{subsec:polynomial_ex}, we set the finite element space for the spatial discretization to be the span of piecewise linear polynomials and use a uniform triangulation of $D$ with $4,934$ elements and $2,340$ spatial unknowns. 
For the results that follow, we fix the truncation length $N=9$ and correlation length $L_c = 1/64$ in \eqref{eq:log_truncation}. 
To maintain sparsity of the SG system,
we use the strategy of projecting the coefficient $a(x,\y)$ from \eqref{eq:log_truncation} onto the space $\mathcal{P}_{\Lambdar}(\Gamma)$, as in \eqref{eq:a_r_truncation}, where $\Lambdar = \LambdarTD$ for the SG-TD approximation, and $\Lambdar = \LambdarSM$ for the SG-SM approximation, obtaining the matrix $\Kr$ from \eqref{eq:fully_discrete_SG_algebraic_sep_KuF}.
We then increase $r$ while $p$ is fixed until the error in the solution stagnates, in practice finding that, for this problem, $r=p$ is sufficient to guarantee the error of the projection does not exceed that of the solution, while maintaining sparsity of the linear system.

\begin{figure}[!htb]
\begin{center}
\includegraphics[clip=true,trim=02mm 02mm 16mm 12mm,width=0.326\textwidth]{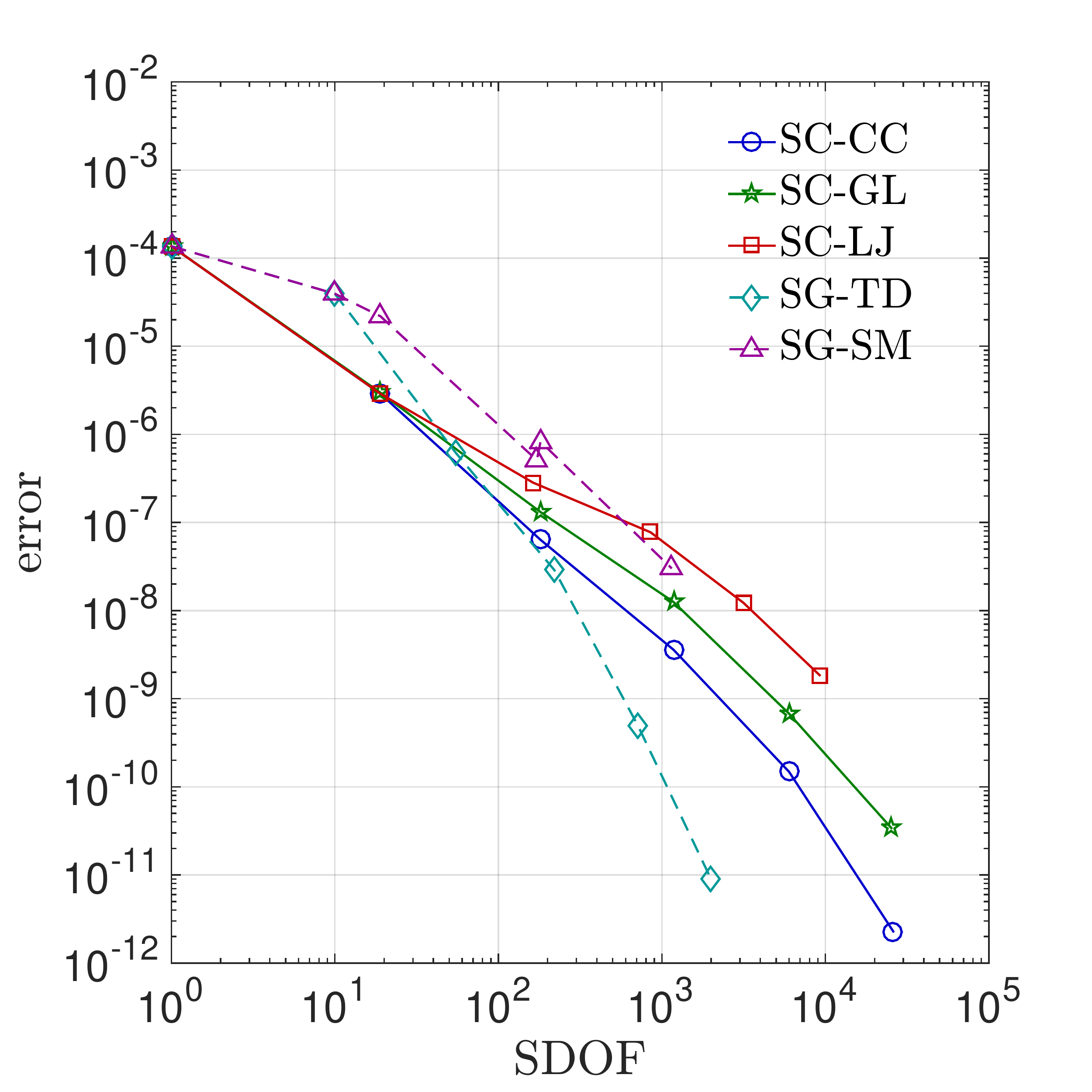}
\includegraphics[clip=true,trim=02mm 02mm 16mm 12mm,width=0.326\textwidth]{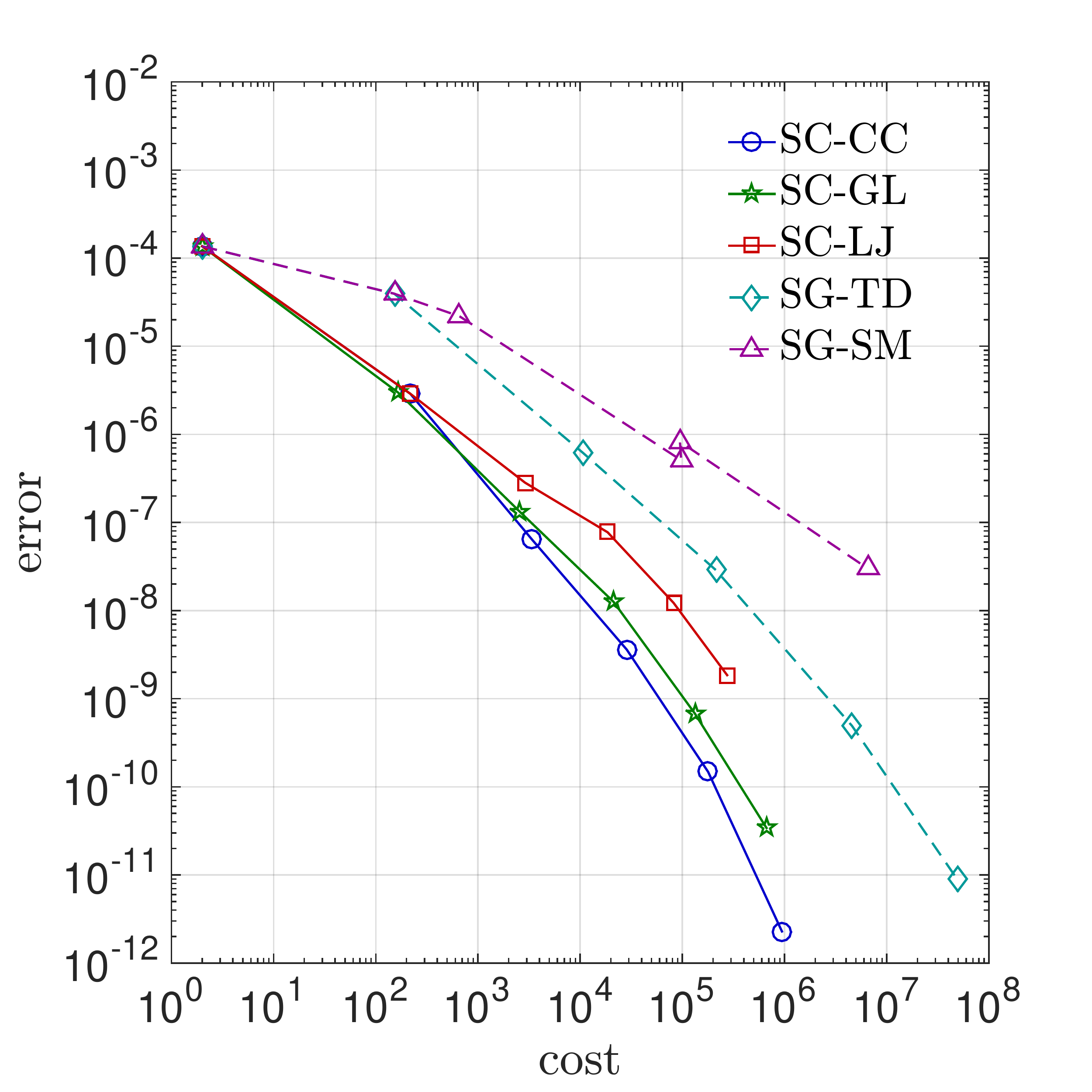}
\caption{{\bf Left:} Error versus SDOFs in solving problem 
    \eqref{eq:model_problem} with coefficient \eqref{eq:log_truncation} and forcing $f(x_1,x_2,\omega) = \cos(x_1)\sin(x_2)$.
    {\bf Right:} Error versus cost with the work estimates given in \eqref{eq:SG_cost_solve_precon} and \eqref{eq:SC_cost_solve_precon}. 
}
\label{fig:log_SDOF_and_cost}
\end{center}
\end{figure}

Figure \ref{fig:log_SDOF_and_cost} compares the error versus SDOFs. There we see that for order $p\geq 3$, the SG-TD approximation provides the best approximation with respect to SDOFs. As discussed in \S \ref{subsec:complexity_comparison}, this is to be expected since the computational complexity of solving the coupled and decoupled systems is not taken into account.
Figure \ref{fig:log_SDOF_and_cost} also displays the convergence in error versus the total computational cost of solving the system with the work estimates of \eqref{eq:SG_cost_solve_precon} and \eqref{eq:SC_cost_solve_precon}. 
Here however, the results show that the SGFEM requires significantly more work to obtain the same error than the SCFEM. 
We also observe the change in rate discussed in \S\ref{subsec:complexity_comparison} in this case, as the work required to solve \eqref{eq:model_problem} with the coefficient $a(x,\y)$ from \eqref{eq:log_truncation} now depends on the order $r$ of the projection used in the approximation of $a(x,\y)$.

For the TD-SG approximation, this can be seen as a consequence of the fact that when $r=p$, the cost of solving \eqref{eq:fully_discrete_SG_algebraic_sep_KuF} with the PCG method is of the order $\mathcal{O}(\Jh M_{\lceil p/2 \rceil}^3 \NiterSG)$, growing much more rapidly than the cost in the affine and polynomial coefficient cases, e.g., Examples \ref{exmp:affine_coefficient} and \ref{exmp:polynomial_coefficient}, as we increase the order $p$. 
Table \ref{tab:savings_log_ex} shows the amount of work required to achieve an error on the order of $10^{-k}$ for some values of $k$ in terms of the total number of matrix-vector products required by both the SC-CC and SG-TD approximations.

\begin{table}
\begin{center}
\begin{adjustbox}{max width=\textwidth}
\begin{tabular}{|c|l|r||c|l|r|} \hline
SC-CC & SC-CC                   & Mat-vec cost  & SG-TD & SG-TD                   & Mat-vec cost  \\
Level & Error                   & of SC-CC      & Order & Error                   & of SG-TD      \\ \hline
0     & $1.3626\times 10^{-4}$  & 2             & 0     & $1.3626\times 10^{-4}$  & 4             \\
1     & $2.8884\times 10^{-6}$  & 218           & 1     & $3.9444\times 10^{-5}$  & 152           \\
2     & $6.3652\times 10^{-8}$  & 3,398         & 2     & $6.1427\times 10^{-7}$  & 10,710        \\
3     & $3.6021\times 10^{-9}$  & 28,638        & 3     & $2.8851\times 10^{-8}$  & 213,010       \\
4     & $1.4794\times 10^{-10}$ & 178,894       & 4     & $4.9210\times 10^{-10}$ & 4,579,575     \\ 
5     & $2.2869\times 10^{-12}$ & 944,220       & 5     & $8.9123\times 10^{-12}$ & 49,089,051    \\ \hline
\end{tabular}
\end{adjustbox}
\end{center}
\caption{Comparison of cost in matrix-vector products for solving problem \eqref{eq:model_problem} with coefficient \eqref{eq:log_truncation} and forcing $f(x_1,x_2,\omega) = \cos(x_1)\sin(x_2)$ using the SC-CC and SG-TD approximations, with the strategy of picking the CG tolerance to be $\|\bbE[u_\text{\em ex} - \uSGCGr ]\|_{\ell^\infty}/10$ for the SGFEM and $\|\bbE[ u_\text{\em ex} - \uSCCG] \|_{\ell^\infty}/10$ for the SC method. Cost in matrix-vector products for SG-TD method is given by \eqref{eq:SG_cost_solve_precon}
    and for SC-CC is given by \eqref{eq:SC_cost_solve_precon} normalized by the cost of a finite element matrix vector product.
}
\label{tab:savings_log_ex}
\end{table}

\section{Conclusions}
\label{sec:conc}

In this work, we presented explicit cost bounds for applying the SGFEM to the solution of an elliptic PDE having both affine and non-affine random coefficients. 
To this end, we have conducted a rigorous counting argument for the sparsity of the linear system that results from the SG discretization with a global orthogonal basis defined on an isotropic total degree index set.
Our analysis shows that when the coefficient is an affine or non-affine function of the random variables having fixed polynomial order, the computational cost of solving the coupled SG system grows linearly with the dimension of the polynomial subspace. 
In these cases, the results only differ by a constant depending on the polynomial order of the random coefficient and the dimension of the parameter domain.

On the other hand, when the coefficient is a non-affine, transcendental function of the random variables 
requiring an additional orthogonal expansion, 
our analysis shows that the computational complexity,
no longer grows linearly with the polynomial subspace dimension. 
For such coefficients, we are able to provide bounds on the complexity that depend on the truncation order of the coefficient.
These estimates imply that a truncation of the expansion should be used, when possible, though attention must be paid to the well-posedness of the resulting PDE.

The analysis conducted herein motivates the study of the total computational complexity of obtaining fully discrete approximations with such methods. 
We have seen that, despite the fact that the SG method yields an approximation that is optimal in the $L^2$ sense for a given polynomial subspace, the associated computational costs of obtaining SG approximations are not optimal for all problems.
Moreover, we have observed, both through theoretical comparisons and numerical examples, that changing the underlying polynomial subspace and method used for obtaining the fully discrete approximation can often yield a solution that requires far less work to obtain, but has the same error.

\subsection*{Acknowledgements}

The first author would like to acknowledge Dr. Miroslav Stoyanov for his insightful comments and assistance in producing the stochastic collocation results with the TASMANIAN package \cite{TASMANIAN}.

\bibliographystyle{siam}
\bibliography{SGFEM}

\section{Appendix}
\label{sec:appendix}

\begin{pf}[of Corollary \ref{cor:sharp_bounds}] 
When $N=1$ we denote $\rb = r \in \N_0$ and note that from Theorem 4.1, 
$$c(r,\ell) =\begin{cases} \# \mathbf{S}(r,\ell) & r \;\;\text{even},\;\; \ell = r/2 \\ 2 \# \mathbf{S}(r,\ell) & \text{otherwise,} \end{cases}$$
with 
$\mathbf{S}(r,\ell) = \{ s \in \N_0 : s = \ell, s \leq r \} = \{\ell\}$ for $\ell \leq r$ and $\emptyset$ otherwise.
We distinguish in cases:
\begin{itemize}
\item Case $r=2k$, $k\in\N_0$,
\begin{enumerate}
\item when $0 \leq r \leq p$,
\begin{align*}
\text{nnz}(\mathbf{G}_r) = \sum_{\ell = \lceil r/2 \rceil}^r c(r, \ell) { 1 + p - \ell \choose p - \ell } 
& = {1 + p - k \choose p - k } + \sum_{\ell = k+1}^{2k} 2 { 1 + p - \ell \choose p - \ell } \\
& = (1 + p - k) - k(3 k - 2 p - 1) \\
& = 1 + p - 4 k^2 + 2kp + k^2 \\
& = (p - 2k + 1)(2k + 1) + k^2 \\
& = (p - r + 1)(r + 1) + k^2.
\end{align*}

\item when $p + 1 \leq r \leq 2 p$, we have 
$\frac{p+1}{2} \leq k \leq p$, so
\begin{align*}
\text{nnz}(\mathbf{G}_r) = \sum_{\ell = \lceil r/2 \rceil}^r c(r, \ell) { 1 + p - \ell \choose p - \ell } 
& = {1 + p - k \choose p - k } + \sum_{\ell = k+1}^{2k} 2 { 1 + p - \ell \choose p - \ell } \\
& = (1 + p - k) + \sum_{\ell = k+1}^{p} 2 { 1 + p - \ell \choose p - \ell } \\
& = (1 + p - k) + (p-k)(p-k+1) \\
& = (1 + p - k)^2.
\end{align*}

\item when $r > 2p$, then 
$k > p$, so 
\begin{align*}
\text{nnz}(\mathbf{G}_r) = \sum_{\ell = \lceil r/2 \rceil}^r c(r, \ell) { 1 + p - \ell \choose p - \ell } 
= {1 + p - k \choose p - k } + \sum_{\ell = k+1}^{2k} 2 { 1 + p - \ell \choose p - \ell } =0,
\end{align*}
since $p-k<0$ and $l>k \Rightarrow p - \ell < p - k < 0$.
\end{enumerate}

\item Case $r=2k+1$, $k\in\N_0$,
\begin{enumerate}
\item when $0 \leq r \leq p$, then $\lceil r/2 \rceil = \lceil (2k+1)/2 \rceil = \lceil k + 1/2 \rceil = k + 1$, so 
\begin{align*}
\text{nnz}(\mathbf{G}_r) = \sum_{\ell = \lceil r/2 \rceil}^r c(r, \ell) { 1 + p - \ell \choose p - \ell } 
& = 2 \sum_{\ell = k+1}^{2k+1} { 1 + p - \ell \choose p - \ell } \\ 
& = -(1+k)(3k - 2p) \\
& = -4k + 2p - 4k^2 + 2kp + k^2 + k \\
& = -2k (2k + 2) + p (2k + 2) + k^2 + k \\
& = (p - 2k)(2k + 2) + k^2 + k \\
& = (p - r + 1)(r + 1) + k^2 + k .
\end{align*}

\item when $p + 1 \leq r \leq 2 p$, then 
$p/2 \leq k \leq p - 1/2$, so 
\begin{align*}
\text{nnz}(\mathbf{G}_r) = \sum_{\ell = \lceil r/2 \rceil}^r c(r, \ell) { 1 + p - \ell \choose p - \ell } 
= 2 \sum_{\ell = k + 1}^{2k+1} { 1+ p - \ell \choose p - \ell} 
& = 2 \sum_{\ell = k+1}^p {1 + p - \ell \choose p - \ell} \\
& = (p - k) ( p - k + 1) .
\end{align*}

\item when $r > 2p$, then 
$k > p - 1/2$, so 
\begin{align*}
\text{nnz}(\mathbf{G}_r) = \sum_{\ell = \lceil r/2 \rceil}^r c(r, \ell) { 1 + p - \ell \choose p - \ell } 
& = \sum_{\ell = k+1}^{2k+1} c(r,\ell) { 1 + p - \ell \choose p - \ell} = 0,
\end{align*}
since $k > p - 1/2 \Rightarrow p - \ell \leq p - (k+1) = p - k - 1 < p - (p - 1/2) -1 = -1/2$. \qed
\end{enumerate}
\end{itemize}
\end{pf}

\end{document}